\documentclass{article}
\usepackage{amsmath,amssymb,amsfonts,amsthm}
\usepackage[all]{xy}
\usepackage{color}

\usepackage{rotating}

\begin{document}

\newtheorem{theorem}{Theorem}[subsection]
\newtheorem{proposition}[theorem]{Proposition}

\newtheorem{corollary}[theorem]{Corollary}
\newtheorem{observation}[theorem]{Observation}

\newtheorem{theoremapp}{Theorem A.\!\!}
\newtheorem{propositionapp}{Proposition A.\!\!}

\newtheorem{corollaryapp}{Corollary A.\!\!}
\newtheorem{observationapp}{Observation A.\!\!}
\newtheorem{lemmaapp}{Lemma A.\!\!}

\theoremstyle{definition}
\newtheorem{definition}[theorem]{Definition}
\newtheorem{dfn}[theorem]{definition}
\newtheorem{remark}[theorem]{Remark}
\newtheorem{example}[theorem]{Example}

\newtheorem{definitionapp}{Definition A.\!\!}
\newtheorem{dfnapp}{definitionapp}
\newtheorem{remarkapp}{Remark A.\!\!}
\newtheorem{exampleapp}{Example A.\!\!}

\newdir^{ (}{{}*!/-5pt/@^{(}}

\newcommand {\jim}[1]{{\marginpar{\textcolor{red}{\Huge{$\star$}}}\scriptsize{\bf \textcolor{red}{JIM}:}\scriptsize{\ #1 \ }}}

\newcommand {\domenico}[1]{{\marginpar{\textcolor{blue}{\Huge{$\star$}}}\scriptsize{\bf \textcolor{blue}{DOMENICO}:}\scriptsize{\ #1 \ }}}	

\newcommand {\urs}[1]{{\marginpar{\textcolor{green}{\Huge{$\star$}}}\scriptsize{\bf \textcolor{green}{URS}:}\scriptsize{\ #1 \ }}}

\def\proof {{\it Proof.}\hspace{7pt}}

\def\endofproof{\hfill{$\square$}\\}
\def\dfn{Definition}
\def\rmk{Remark}
\def\Cech{{{\v C}ech}}

\def\Loo{{$L_\infty$}}

\pagestyle{myheadings}
\markboth{Domenico-Urs-Jim}{$\infty$-Chern-Weil theory\hskip10ex \today}

\title{\Cech \  cocycles for differential characteristic classes -- An $\infty$-Lie theoretic construction}
\author{Domenico Fiorenza, Urs Schreiber  and Jim Stasheff}

\date{July 8, 2011}

\maketitle	

\begin{abstract}
  What are called \emph{secondary characteristic classes} in Chern-Weil theory
  are a refinement of ordinary
  characteristic classes of principal bundles from cohomology to differential cohomology.
  We consider the problem of refining the construction of secondary characteristic classes  
  from cohomology sets to cocycle spaces; 
  and from Lie groups to higher connected covers of Lie groups by smooth $\infty$-groups, i.e.,  by smooth groupal $A_\infty$-spaces. Namely, we realize differential characteristic
classes as morphisms from $\infty$-groupoids of smooth principal $\infty$-bundles with
connections to $\infty$-groupoids of higher $U(1)$-gerbes with connections.
  This allows us to study the homotopy fibers of the differential characteristic maps thus 
  obtained 
  and to show how these describe differential obstruction problems. This applies in particular to the
  higher twisted differential spin structures called \emph{twisted differential string structures} 
  and \emph{twisted differential fivebrane structures}.
\end{abstract}

\newpage

\centerline{\bf Summary}
\vskip .2cm

 What are called \emph{secondary characteristic classes} in Chern-Weil theory
  are a refinement of ordinary
  characteristic classes of principal bundles from cohomology to differential cohomology.
  We consider the problem of refining the construction of secondary characteristic classes  
  from cohomology sets to cocycle spaces; 
  and from Lie groups to higher connected covers of Lie groups by smooth $\infty$-groups, i.e.,  by smooth groupal $A_\infty$-spaces. Namely, we realize differential characteristic
classes as morphisms from $\infty$-groupoids of smooth principal $\infty$-bundles with
connections to $\infty$-groupoids of higher $U(1)$-gerbes with connections.
  This allows us to study the homotopy fibers of the differential characteristic maps thus 
  obtained 
  and to show how these describe differential obstruction problems. This applies in particular to the
  higher twisted differential spin structures called \emph{twisted differential string structures} 
  and \emph{twisted differential fivebrane structures}.

  To that end 
  we define for every $L_\infty$-algebra $\mathfrak{g}$ a smooth $\infty$-group $G$ integrating
  it, and define smooth $G$-principal $\infty$-bundles with connection. For every 
  $L_\infty$-algebra cocycle of suitable degree, we give a refined 
  $\infty$-Chern-Weil homomorphism that sends these $\infty$-bundles to classes in differential
  cohomology that lift the corresponding curvature characteristic classes.
  
  When applied to the canonical 3-cocycle of the Lie algebra of  
  a simple and simply connected Lie group $G$
  this construction gives a refinement of the secondary first fractional Pontryagin class
  of $G$-principal bundles
  to cocycle space. Its homotopy fiber is the 2-groupoid of smooth
  $\mathrm{String}(G)$-principal 2-bundles with 2-connection, where $\mathrm{String}(G)$
  is a smooth 2-group refinement of the topological string group. Its homotopy fibers over
  non-trivial classes we identify with the 2-groupoid of 
  \emph{twisted differential string structures} that appears in the 
  Green-Schwarz anomaly cancellation
  mechanism of heterotic string theory.

  Finally, when our construction is applied to the 
  canonical 7-cocycle on the Lie 2-algebra of the String-2-group, it produces 
  a secondary characteristic map for $\mathrm{String}$-principal 2-bundles 
  which refines the second fractional Pontryagin class.
  Its homotopy fiber is the 6-groupoid of principal 6-bundles with 6-connection
  over the \emph{Fivebrane 6-group}. Its homotopy fibers over nontrivial classes 
  are accordingly \emph{twisted differential fivebrane structures} that have beeen
  argued to control the anomaly cancellation mechanism in magnetic dual heterotic
  string theory.

\vfill

\begin{center}
  Further online resources for this document can be found at
  \\
  \medskip
  \verb"http://ncatlab.org/schreiber/show/differential+characteristic+classes"
\end{center}

\newpage

\tableofcontents

\section{Introduction}

Classical Chern-Weil theory (see for instance \cite{GHV, milnor-stasheff, hopkins-singer}) provides a toolset for refining characteristic classes of smooth principal bundles from ordinary integral cohomology to differential cohomology.

This can be described as follows. For $G$ a topological group and
$P \to X$ a $G$-principal bundle,  to any \emph{characteristic class} $[c]\in
H^{n+1}(BG,\mathbb{Z})$, there is associated
 a characteristic class of the bundle, $[c(P)] \in H^{n+1}(X,\mathbb{Z})$. This can be seen as the homotopy class of the composition
\[
 X\xrightarrow{P}BG\xrightarrow{c}K(\mathbb{Z},n+1)
 \] of the classifying map $X\xrightarrow{P}BG$ of the bundle with the \emph{characteristic map} $BG\xrightarrow{c}K(\mathbb{Z},n+1)$. If $G$ is a compact connected Lie group, and with 
real coefficients, there is a graded commutative algebra isomorphism between
$H^\bullet(BG,\mathbb{R})$ and the algebra $\mathrm{inv}(\mathfrak{g})$ of
$\mathrm{ad}_G$-invariant  polynomials on the Lie algebra $\mathfrak{g}$ of $G$. In
particular, any characteristic class $c$ will correspond to such an
invariant polynomial $\langle-\rangle$. The Chern-Weil 
homomorphism associates to a choice of \emph{connection} $\nabla$ on a  $G$-principal bundle $P$ the closed differential form $\langle F_\nabla\rangle$ on $X$, where $F_\nabla$ is the curvature of $\nabla$. The de Rham cocycle $\langle F_\nabla\rangle$ is a representative for the 
characteristic class $[c(P)]$ in $H^\bullet(X,\mathbb{R})$. This construction can be carried out at a local level: instead of considering a globally defined connection $\nabla$, one can consider an open cover $\mathcal{U}$ of $X$ and local connections $\nabla_i$ on $P\vert_{U_i}\to U_i$; then the local differential forms $\langle F_{\nabla_i}\rangle$ define a cocycle in the \v{C}ech-de Rham complex, still representing the cohomology class of $c(P)$.
\par
There is a refinement of this construction to what is sometimes called \emph{secondary characteristic classes}: 
the differential form $\langle F_\nabla \rangle$ may itself be understood as the 
higher curvature form of a higher circle-bundle-like structure $\hat{\mathbf{c}}(\nabla)$ whose higher
Chern-class is $c(P)$. In this refinement, both the original characteristic class $c(P)$ as well as its curvature differential form $\langle F_\nabla \rangle $ are unified in one single object.
This single object has originally been formalized as a \emph{Cheeger-Simons differential character}.
It may also be conceived of as a cocycle in the \v{C}ech-Deligne complex, a refinement of the \v{C}ech-de Rham complex \cite{hopkins-singer}. Equivalently, as we discuss here, these objects may naturally be described in terms of what we want to call \emph{circle 
$n$-bundles with connection}: smooth bundles whose structure group is a smooth refinement -- which we write $\mathbf{B}^{n}U(1)$ -- of the topological group $B^{n}U(1) \simeq K(\mathbb{Z},n+1)$, endowed with a smooth connection of higher order. For low $n$, such $\mathbf{B}^{n}U(1)$-principal bundles are known (more or less explicitly) as \emph{$(n-1)$-bundle gerbes}. 

The fact that we may think of  $\hat{\mathbf{c}}(\nabla)$ as being 
a smooth principal higher bundle with connection suggests that it makes sense to ask if there is a general definition of smooth $G$-principal $\infty$-bundles, for smooth $\infty$-groups $G$, and whether the Chern-Weil homomorphism extends on those to an \emph{$\infty$-Chern-Weil homomorphism}.  
Moreover, since $G$-principal bundles naturally form a parameterized groupoid -- a stack -- 
and circle $n$-bundles naturally form a 
parameterized $n$-groupoid -- an $(n-1)$-stack, an $n$-truncated $\infty$-stack, 
it is natural to ask whether we can refine the
construction of differential characteristic classes to these $\infty$-stacks.
Motivations for considering this are threefold:

\begin{enumerate}
  \item The ordinary Chern-Weil homomorphism only knows about characteristic classes of classifying spaces $B G$ for $G$ a Lie group. Already before considering the refinement to differential cohomology, this misses useful cohomological information about \emph{connected covering groups} of $G$.
  
      For instance, for $G = \mathrm{Spin},$ the Spin group, there is the second Pontryagin class represented by a map $p_2 : B \mathrm{Spin} \to B^8 \mathbb{Z}$. But on some $\mathrm{Spin}$-principal bundles $P \to X$ classified by a map $g : X \to B \mathrm{Spin}$, this class may be further divisible: there is a topological group 
      $\mathrm{String}$, called the \emph{String group}, such that we have a commuting diagram
      
      $
        \xymatrix{
           & B \mathrm{String} \ar[r]^{\frac{1}{6} p_2} \ar[d]& B^8 \mathbb{Z} \ar[d]^{\cdot 6}
           \\
           X \ar[r]^g 
             \ar@{-->}[ur]^{\tilde g}
           & B \mathrm{Spin} \ar[r]^{p_2} & B^8 \mathbb{Z}           
        }
      $
  
      of topological spaces, where the morphism on the right is given on $\mathbb{Z}$
      by multiplication with 6 \cite{SSSII}. This means that if $P$ happens to admit a  
      \emph{String structure} exhibited by a lift $\tilde g$ of its classifying map $g$
      as indicated, then its second Pontryagin class 
      $[p_2(P)] \in H^8(X,\mathbb{Z})$ is divisible by 6. But this refined information
      is invisible to the ordinary Chern-Weil homomorphism: 
      while $\mathrm{Spin}$ canonically has the structure
      of a Lie group, $\mathrm{String}$ cannot have a finite-dimensional Lie group structure
      (because it is a $B U(1)$-extension, hence has cohomology in arbitrary high degree)
      and therefore the ordinary Chern-Weil homomorphism can not model this fractional 
      characteristic class.
      
      But it turns out that $\mathrm{String}$ does have a natural smooth structure when regarded as 
      a higher group -- a 2-group in this case \cite{bcss, henriques}. We write $\mathbf{B}\mathrm{String}$ for the
      corresponding smooth refinement of the classifying space. 
      As we shall show, there is an $\infty$-Chern-Weil homomorphism
      that does apply and produces for every smooth $\mathrm{String}$-principal 2-bundle
      $\tilde g : X \to \mathbf{B}\mathrm{String}$ a smooth circle 7-bundle with connection, which 
      we write $\frac{1}{6}{\hat {\mathbf{p}}}_2(\tilde g)$. Its curvature 8-form is a 
      representative in de Rham cohomology of the fractional second Pontryagin class.

      Here and in the following
      \begin{itemize}
        \item boldface denotes a refinement from continuous (bundles) to \emph{smooth} (higher bundles);
        \item the hat denotes further \emph{differential} refinement 
           (equipping higher bundles with smooth connections).
      \end{itemize}      
  
      In this manner, the $\infty$-Chern-Weil homomorphism gives cohomological information beyond that of
      the ordinary Chern-Weil homomorphism. And this is only the beginning of a pattern:
      the sequence of smooth objects that we considered continues further as
      
      $$
        \cdots \to \mathbf{B}\mathrm{Fivebrane} \to \mathbf{B}\mathrm{String}
          \to \mathbf{B}\mathrm{Spin} \to \mathbf{B} SO \to \mathbf{B}O
      $$

      to a smooth refinement of the \emph{Whitehead tower} of $B O$. One way to think
      of $\infty$-Chern-Weil theory is as a lift of ordinary Chern-Weil theory 
      along such smooth Whitehead towers.

\item
  Traditionally the construction of secondary characteristic classes is exhibited on
  single cocycles and then shown to be independent of the representatives of the corresponding
  cohomology class. But this indicates that one is looking only at the connected components of 
  a more refined construction that explicitly sends cocycles to cocycles, and sends
  coboundaries to coboundaries such that their composition is respected up to higher
  degree coboundaries, which in turn satisfy their own coherence condition, and so forth. 
  In other words: a map between the full cocycle $\infty$-groupoids.
  
  The additional information encoded in such a refined secondary differential characteristic map 
  is equivalently found in the collection of the \emph{homotopy fibers} of the map,
  over the cocycles in the codomain. These homotopy fibers answer the question: which 
  bundles with connection have differential characteritsic class equivalent to some
  fixed class, which of their gauge transformations respect the choices of equivalences, which of 
  the higher gauge of gauge transformations respect the chosen gauge transformations, and so on. 
  This yields refined cohomological
  information whose knowledge is required in several applications of differential cohomology,
  indicated in the next item.

\item
   Much of the motivation for studies of differential cohomology originates in the 
   applications this theory has to the description of higher gauge fields in physics. 
   Notably the seminal article
   \cite{hopkins-singer} that laid the basis of generalized differential cohomology grew out of the 
   observation that this is the right machinery that describes subtle phenomena of 
   quantum anomaly cancellation in string theory, discussed by Edward Witten and others 
   in the 1990s, further spelled out 
   in \cite{freed}. In this context, the need for refined fractional characteristic classes 
   and their homotopy fibers appears.

   In higher analogy to how the quantum mechanics of a spinning particle requires
   its target space to be equipped with a $\mathrm{Spin}$-structure that is differentially
   refined to a $\mathrm{Spin}$-principal bundle with connection, the quantum dynamics
   of the (heterotic) superstring requires target space to be equipped with a
   differential refinement of  a $\mathrm{String}$-structure. Or rather, since 
   the heterotic string contains besides the gravitational $\mathrm{Spin}$-bundle
   also a $U(n)$-gauge bundle, of a \emph{twisted} $\mathrm{String}$-structure for
   a specified twist. We had argued in \cite{SSSIII} that these differentially refined
   string backgrounds are to be thought of as twisted differential structures in the 
   above sense. With the results of the present work this argument is lifted from 
   a discussion of local $\infty$-connection data to the full differential cocycles.
   We shall show that by standard homotopy theoretic arguments this allows a simple 
   derivation of the properties of untwisted
   differential string structures that have been found in \cite{waldorf}, and generalize
   these to the twisted case and all the higher analogs. 
   
   Namely, moving up along the Whitehead tower of $O(n)$, one can next ask for 
   the next higher characteristic class on $\mathrm{String}$-2-bundles and 
   its differential refinement to a secondary characteristic class. In \cite{SSSII}
   it was argued that this controls, in direct analogy to the previous case, the
   quantum super-5-brane that is expected to appear in the magnetic dual description
   of the heterotic target space theory. With the tools constructed here 
   the resulting 
   \emph{twisted differential fivebrane structures} can be analyzed in analogy
   to the case of string structures.
   
   Our results allow an analogous description of twisted differential
   structures of ever higher covering degree, but beyond the 5-brane it is currently unclear whether
   this still has applications in physics. However, there are further variants in low
   degree that do:
   
  for instance  for every $n$ there is a canonical 4-class on pairs of $n$-torus bundles and dual $n$-torus
   bundles. This has a differential refinement and thus we can apply our results to 
   this situation to produce parameterized 2-groupoids of the corresponding higher twisted
   differential torus-bundle extensions. We find that the connected components of these
   2-groupoids are precisely the \emph{differential T-duality pairs} that arise in the 
   description of differential T-duality of strings in \cite{KahleValentino}. 
   
   This suggests that there are more applications of refined higher differential
   characteristic maps in string theory, but here we shall be content with looking into 
   these three examples.
\end{enumerate}

In this article, we shall define connections on principal $\infty$-bundles and the action 
of the $\infty$-Chern-Weil homomorphism in a natural but maybe still somewhat \emph{ad hoc} 
way, which here we justify mainly by the two main theorems about two examples that we prove,
which we survey in a moment. 
The construction uses essentially standard tools of differential geometry.
The construction can be derived from \emph{first principles} 
as a model (in the precise sense of \emph{model category theory}) 
for a general abstract construction that exists in \emph{$\infty$-topos theory}. This
abstract theory is discussed in detail elsewhere \cite{survey}.

Notice that our approach goes beyond that of \cite{hopkins-singer} in two ways: the
$\infty$-stacks we consider remember \emph{smooth} gauge transformation and
thus encode smooth structure of principal $\infty$-bundles already on cocycles
and not just in cohomology; secondly, we describe \emph{non-abelian} phenomena, such as connections on 
principal bundles for non-abelian structure groups, and more in general $\infty$-connections
for non-abelian structure smooth $\infty$-groups, such as the $\mathrm{String}$-2-group
and the $\mathrm{Fivebrane}$-6-group. This is the very essence of (higher) Chern-Weil
theory: to characterize non-abelian cohomology by abelian characteristic classes. 
Since \cite{hopkins-singer} work with spectra, nothing non-abelian is directly available there.
On the other hand, the construction we describe does not as easily allow differential 
refinements of cohomology theories represented by non-connective spectra.

\medskip
We now briefly indicate the means by which we will approach these issues in the following.

The construction that we discuss is the result of applying a refinement of the machine
of \emph{$\infty$-Lie integration} \cite{henriques, getzler} to the $L_\infty$-algebraic 
structures discussed in \cite{SSSI, SSSIII}:

For $\mathfrak{g}$ an $L_\infty$-algebra, 
its \emph{Lie integration} to a Lie $\infty$-group $G$ with smooth classifying object
$\mathbf{B}G$ turns out to be encoded in the simplicial presheaf given by the assignment 
to each smooth test manifold $U$ of the simplicial set
$$
  \exp_\Delta(\mathfrak{g}) : 
  (U,[k]) \mapsto \mathrm{Hom}_{\mathrm{dgAlg}}(\mathrm{CE}(\mathfrak{g}), 
   \Omega^\bullet(U \times \Delta^k)_{\mathrm{vert}})
  \,,
$$
where $\mathrm{CE}(\mathfrak{g})$ is the Chevalley-Eilenberg algebra of $\mathfrak{g}$
and `vert' denotes forms which see only vector fields along $\Delta^k$.
This has a canonical projection $\exp_\Delta(\mathfrak{g}) \to \mathbf{B}G$, hence the name 
$\exp_\Delta(\mathfrak{g}).$ 
One can think of this as saying that a $U$-parameterized smooth family of $k$-simplices
in $G$ is given by the \emph{parallel transport} over the $k$-simplex 
of a flat $\mathfrak{g}$-valued vertical 
differential form on the trivial simplex bundle $U \times \Delta^k \to U$. 
This we discuss  in detail in section \ref{section.principal_infinity-bundles}. 

The central step of our construction is
a \emph{differential refinement}
$\mathbf{B}G_{\mathrm{diff}}$ of $\mathbf{B}G$, where the above is enhanced to
$$
  \exp_\Delta(\mathfrak{g})_{\mathrm{diff}} : 
  (U,[k]) \mapsto 
  \left\{
   \raisebox{23pt}{
  \xymatrix{
    \Omega^\bullet(U \times \Delta^k)_{\mathrm{vert}}
    &
    \mathrm{CE}(\mathfrak{g})
    \ar[l]
    \\
    \Omega^\bullet(U \times \Delta^k)
    \ar[u]
    &
    \mathrm{W}(\mathfrak{g})
    \ar[l]\ar[u]
  }}
  \right\}
  \,,
$$
with $\mathrm{W}(\mathfrak{g})$ the Weil algebra of $\mathfrak{g}$.
We also consider a simplicial sub-presheaf 
$\exp_\Delta(\mathfrak{g})_{\mathrm{conn}} \hookrightarrow \exp_\Delta(\mathfrak{g})_{\mathrm{diff}}$
defined by a certain horizontality constraint. This may be thought of as assigning non-flat
$\mathfrak{g}$-valued forms on the total space of the trivial simplex bundle $U \times \Delta^k$. 
The horizontality constraint generalizes one of the conditions of an
\emph{Ehresmann connection} \cite{ehresmann} on an ordinary $G$-principal bundle. This we discuss
in detail in \ref{section.infinity-stack_of_principal_bundles_with_connection}.

We observe that an $L_\infty$-algebra cocycle  $\mu \in \mathrm{CE}(\mathfrak{g})$
in degree $n$, when we equivalently regard it as a morphism of $L_\infty$-algebras 
$\mu : \mathfrak{g} \to b^{n-1}\mathbb{R}$ to the Eilenberg-MacLane object $b^{n-1}\mathbb{R}$,
tautologically \emph{integrates} to a morphism
$$
  \exp_\Delta(\mu) :  \exp_\Delta(\mathfrak{g}) \to \exp_\Delta(b^{n-1}\mathbb{R})
$$
of the above structures. What we identify as the $\infty$-Chern-Weil homomorphism is 
obtained by first extending this to the differential refinement
$$
  \exp_\Delta(\mu)_{\mathrm{diff}} : \exp_\Delta(\mathfrak{g})_{\mathrm{diff}} \to \exp_\Delta(b^{n-1}\mathbb{R})_{\mathrm{diff}}
$$
in a canonical way -- this we shall see introduces \emph{Chern-Simons elements} -- and then 
descending the construction along the projection 
$\exp(\mathfrak{g})_{\mathrm{diff}} \to \mathbf{B}G_{\mathrm{diff}}$.
This quotients out a lattice $\Gamma \subset \mathbb{R}$ and makes the resulting
higher bundles with connection be circle $n$-bundles with connection which represent
classes in differential cohomology. This we discuss in section 
\ref{section.infinity_chern_weil_homomorphism}.

Finally, in the last part of section \ref{section.infinity_chern_weil_homomorphism} 
we discuss two classes of applications and obtain the following
statements.

\begin{theorem} 
  \label{theorem.main}
  Let $X$ be a paracompact smooth manifold and choose a good open cover $\mathcal{U}$.

  \noindent Let $\mathfrak{g}$ be a semisimple Lie algebra with normalized 
  binary Killing form $\langle -,-\rangle$ in transgression with the 3-cocycle
  $\mu_3 = \frac{1}{2}\langle -,[-,-]\rangle$. Let $G$ be the corresponding simply
  connected Lie group.
  
  \begin{itemize}
    \item {\bf 1.}
    Applied to this $\mu_3$, the $\infty$-Chern-Weil homomorphism 
  $$
    \exp(\mu)_{\mathrm{conn}} 
     :  
     \check{C}(\mathcal{U},\mathbf{B}G_{\mathrm{conn}}) \to 
    \check{C}(\mathcal{U}, \mathbf{B}^3 U(1)_{\mathrm{conn}})
  $$
  from \v{C}ech cocycles with coefficients in the complex that classifies $G$-principal bundles
  with connection
  to \v{C}ech-Deligne cohomology in degree 4 is a fractional multiple of the Brylinski-McLaughlin
  construction \cite{brylinski-mclaughlin} of {\Cech}-Deligne cocycles representing the 
  differential refinement of the characteristic class corresponding to $\langle - , -\rangle$.
    
  In particular, in cohomology it represents the 
  refined Chern-Weil homomorphims
  $$
     \frac{1}{2}\hat {\mathbf{p}}_1 : H^1(X,G)_{\mathrm{conn}} \to \hat H^4(X,\mathbb{Z})
  $$ 
  induced by the Killing form and with coefficients in 
  degree 4 differential cohomology. For 
  $\mathfrak{g} = \mathfrak{so}(n),$ this is the differential refinement 
  of the first fractional Pontryagin class.
  
  \end{itemize}

  \noindent Next let $\mu_7 \in \mathrm{CE}(\mathfrak{g})$ be a 7-cocycle 
  on the semisimple Lie algebra $\mathfrak{g}$ (this is unique up to a scalar factor). Let
  $\mathfrak{g}_{\mu_3} \to \mathfrak{g}$ be the $L_\infty$-algebra-extension of 
  $\mathfrak{g}$ classified by $\mu_3$ (the \emph{string Lie 2-algebra}. Then  
  $\mu_7$ can be seen as a 7-cocycle also on $\mathfrak{g}_{\mu_3}$.

  \begin{itemize}
    \item {\bf 2.}
     Applied to $\mu_7$ regarded as a cocycle on $\mathfrak{g}_\mu$, the 
     $\infty$-Chern-Weil homomorphism produces a map
     $$
       \check{C}(\mathcal{U}, \mathbf{B}\mathrm{String}(G)_{\mathrm{conn}})
       \to 
       \check{C}(\mathcal{U}, \mathbf{B}^7 U(1)_{\mathrm{conn}})
     $$
     from \v{C}ech cocycles with coefficients in the complex that classifies 
     $\mathrm{String}(G)$-2-bundles with connection to degree 8 \v{C}ech-Deligne cohomology. 
     For $\mathfrak{g} = \mathfrak{so}(n)$ this gives a fractional refinement
     of the ordinary refined Chern-Weil homomorphism
     $$
       \frac{1}{6} {\hat {\mathbf{p}}_2} : H^1(X,\mathrm{String})_{\mathrm{conn}} \to
       \hat H^8(X,\mathbb{Z})
     $$
     that represents the differential refinement of the second fractional 
     Pontryagin class on $\mathrm{Spin}$-bundles with $\mathrm{String}$-structure.
  \end{itemize}
\end{theorem}

These are only the first two instances of a more general statement. 
But this will be discussed elsewhere.

\def\dfn{Definition}

\section{A review of ordinary Chern-Weil theory}

We briefly review standard aspects of ordinary Chern-Weil theory whose generalization we consider later on. In this section we assume the reader is familiar with basic properties of Chevalley-Eilenberg and of Weil algebras; the unfamiliar reader can find a concise account at the
beginning of Section \ref{ooLieIntegration}.

\subsection{The Chern-Weil homomorphism}\label{chern-weil-review1}

For $G$ a Lie group and $X$ a smooth manifold, 
the idea of a \emph{connection} on a smooth $G$-principal bundle $P\to X$ 
can be expressed in a variety of equivalent ways:
as a distribution of horizontal spaces on the tangent bundle total space $T P$,
as the corresponding family of projection operators in terms of local connection 1-forms on $X$ or, more generally, as defined by Ehresmann \cite{ehresmann}, and, ultimately, purely algebraically, by H. Cartan \cite{cartan:g-alg-sem-I,cartan:g-alg-sem-II}.
\par
Here, following this last approach, we review how the Weil algebra can be used to give an algebraic description of connections on principal bundles, and of the Chern-Weil homomorphism.
\par
We begin by recalling the classical definition of connection on a $G$-principal bundle $P\to X$ as a $\mathfrak{g}$-valued 1-form $A$ on $P$ which is $G$-equivariant and induces the Maurer-Cartan form of $G$ on the fibers (these are known as the Cartan-Ehresmann conditions). The key insight is then the identification of $A\in \Omega^1(P,\mathfrak{g})$ with a differential graded algebra morphism; this is where the Weil algebra $\mathrm{W}(\mathfrak{g})$ comes in. We will introduce Weil algebras in a precise and intrinsic way in the wider context of Lie $\infty$-algebroids in Section \ref{section.Lie_infinity-algebroids}, so we will here 
content ourselves with  thinking of the Weil algebra of $\mathfrak{g}$
as a perturbation of the Chevalley-Eilenberg cochain complex $\mathrm{CE}(\mathfrak{g})$ for $\mathfrak g$ with coefficients in 
the polynomial algebra generated by the dual $\mathfrak g^*$. More precisely, 
the Weil algebra $\mathrm{W}(\mathfrak g)$ is a commutative dg-algebra freely generated by two copies 
of $\mathfrak g^*$, one in degree 1 and one in degree 2;
the differential $d_{\mathrm{W}}$ is the sum of the Chevalley-Eilenberg differential plus $\sigma$, the shift isomorphism from
$\mathfrak g^*$ in degree 1 to $\mathfrak g^*$ in degree 2, extended as a derivation. 

A crucial property of the Weil algebra is its freeness: dgca morphisms out of the Weil algebra are uniquely and freely determined by graded vector space morphism out of the copy of $\mathfrak{g}^*$ in degree 1. This means that a $\mathfrak{g}$-valued 1-form $A$ on $P$ can be equivalently seen as a dgca morphism 
\[
A:\mathrm{W}(\mathfrak{g})\to \Omega^\bullet(P)
\]
to the de Rham dg-algebra of differential forms on $P$.
Now we can read the Cartan-Ehresmann conditions on a $\mathfrak{g}$-connection as properties of this dgca morphism. First, the Maurer-Cartan form on $G$, i.e., the left-invariant $\mathfrak{g}$-valued form $\theta_G$ on $G$ induced by the identity on $\mathfrak{g}$ seen as a linear morphism $T_eG\to \mathfrak{g}$, is an element of $\Omega^1(G,\mathfrak{g})$, and so it defines a dgca morphism $\mathrm{W}(\mathfrak{g})\to \Omega^\bullet(G)$. This morphism actually factors through the Chevalley-Eilenberg algebra of $\mathfrak{g}$; this is the algebraic counterpart of the fact that the curvature 2-form of $\theta_G$ vanishes. Therefore, the first Cartan-Ehresmann condition on the behaviour of the connection form $A$ on the fibres of $P\to X$ is encoded in the commutativity of the following diagram of differential graded commutative algebras:
\[
\xymatrix{
\,\phantom{mm}\Omega^\bullet(P)_\mathrm{vert}&\mathrm{CE}(\mathfrak{g})\ar[l]_{\phantom{mm}A_\mathrm{vert}}\\
\Omega^\bullet(P)\ar[u]&\mathrm{W}(\mathfrak{g})\ar[l]_{A}\ar[u]
}\,.
\]
In the upper left corner, $\Omega^\bullet(P)_\mathrm{vert}$ is the dgca of \emph{vertical} differential forms on $P$, i.e. the quotient of $\Omega^\bullet(P)$ by the differential ideal consisting of differential forms on $P$ which vanish  when evaluated on a vertical multivector field.
\par
Now we turn to the second Cartan-Ehresmann condition. The symmetric algebra $\mathrm{Sym}^\bullet(\mathfrak{g}^*[-2])$ on $\mathfrak{g}^*$ placed in degree 2 is a graded commutative subalgebra of the Weil algebra $\mathrm{W}(\mathfrak{g})$, but it is not a dg-subalgebra. However, the subalgebra $\mathrm{inv}(\mathfrak{g})$ of $\mathrm{Sym}^\bullet(\mathfrak{g}^*[-2])$ consisting of $\mathrm{ad}_\mathfrak{g}$-invariant polynomials is a dg-subalgebra of $\mathrm{W}(\mathfrak{g})$. The composite morphism of dg-algebras
\[
\mathrm{inv}(\mathfrak{g})\to \mathrm{W}(\mathfrak{g})\xrightarrow{A}\Omega^\bullet(P)
\]
is the evaluation of  invariant polynomials on the \emph{curvature} 2-form of $A$, i.e. on the $\mathfrak{g}$-valued 2-form $F_A=dA+\frac{1}{2}[A,A]$. 
Invariant polynomials are $d_{\mathrm{W}}$-closed as elements in the Weil algebra, therefore, their images in $\Omega^\bullet(P)$ are closed differential forms. Assume now $G$ is connected. Then, if $\langle-\rangle$ is an $\mathrm{ad}_\mathfrak{g}$-invariant polynomial, by the $G$-equivariance of $A$ it follows that the closed differential form $\langle F_A\rangle$ descends to a closed differential form on the base $X$ of the principal bundle. Thus, the second Cartan-Ehresmann condition on $A$ implies the commutativity of the diagram 
\[
\xymatrix{
\Omega^\bullet(P)&\mathrm{W}(\mathfrak{g})\ar[l]_{A}\\
\Omega^\bullet(X)\ar[u]& \mathrm{inv}(\mathfrak{g})\ar[l]_{F_A}\ar[u]
}\,.
\]
Since the image of $\mathrm{inv}(\mathfrak{g})$ in $\Omega^\bullet(X)$ consists of closed forms, we have an induced graded commutative algebras morphism
\[
\mathrm{inv}(\mathfrak{g})\to H^\bullet(X,\mathbb{R}),
\]
the \emph{Chern-Weil homomorphism}. This morphism is independent of the particular connection form chosen and natural in $X$. Therefore we can think of elements of $\mathrm{inv}(\mathfrak{g})$ as representing universal cohomology classes, hence as \emph{characteristic classes}, of $G$-principal bundles. And indeed, if $G$ is a compact connected finite dimensional Lie group, then we have an isomorphism of graded commutative algebras
$\mathrm{inv}(\mathfrak{g}) \cong H^\bullet(BG, \mathbb{R})$,
corresponding to $H^\bullet(G,\mathbb{R})$ being isomorphic to  an exterior algebra  on odd dimensional generators \cite{chern}, the indecomposable Lie algebra cohomology classes of $\mathfrak{g}$. The isomorphism $\mathrm{inv}(\mathfrak{g}) \cong H^\bullet(BG, \mathbb{R})$ is to be thought as the universal Chern-Weil homomorphism. Traditionally this is conceived of in terms
of a smooth manifold version of the universal $G$-principal bundle on $B G$. 
We will here instead refine $B G$ to a smooth $\infty$-groupoid $\mathbf{B}G$. This classfies not just
equivalence classes of $G$-principal bundles but also their automorphisms. We shall argue that
the context of smooth $\infty$-groupoids is the natural place (and \emph{place} translates to  \emph{topos}) in which to conceive of the Chern-Weil homomorphism.

\subsection{Local curvature 1-forms}
\label{OrdinaryConnections}

Next we focus on the description of $\mathfrak{g}$-connections in terms of local $\mathfrak{g}$-valued 1-forms and gauge transformations. We discuss this in terms of the local transition function data from which the
total space of the bundle may be reconstructed. It is this local point
of view that we will explicitly generalize in section \ref{ooLieIntegration}. 
More precisely, in section \ref{section.infinity-stack_of_principal_bundles_with_connection} we will present algebraic data which encode an $\infty$-connection on a trivial higher bundle on a Cartesian space $\mathbb{R}^n$, and will then globalize this local picture by \emph{descent/stackification}.\par
To prepare this general construction, let us show how it works in the case of ordinary $\mathfrak{g}$-connections on $G$-principal bundles. 
For that purpose, consider a Cartesian space $U = \mathbb{R}^n$. Every $G$-principal
bundle on $U$ is equivalent to the trivial $G$-bundle $U \times G$ equipped with the evident
action of $G$ on the second factor, and under stackification this completely characterizes
$G$-principal bundles on general spaces. A connection on this
trivial $G$-bundle is given by a $\mathfrak{g}$-valued 1-form
$A \in \Omega^1(U,\mathfrak{g})$. An isomorphism  
$\xymatrix{A \ar[r]^{g}& A'}$ from the trivial bundle with 
connection $A$ to that with connection $A'$ is given by a function
$g \in C^\infty(U,G)$ such that the equation

\begin{eqnarray}
  \label{gaugetransformation}
  A' = g^{-1} A g + g^{-1} d g
\end{eqnarray}
holds. Here the first term on the right denotes the adjoint action of the Lie group
on its Lie algebra, while the second term denotes the pullback of the Maurer-Cartan
form on $G$ along $g$ to $U$.

We wish to amplify a specific way to understand this formula as the Lie integration of a path of infinitesimal 
gauge transformations: 
write $\Delta^1 = [0,1]$ for the standard interval regarded as a smooth manifold (with boundary)
and consider a smooth 1-form $A \in \Omega^1(U \times \Delta^1, \mathfrak{g})$ on the 
product of $U$ with $\Delta^1$. If we think of this as the trivial interval bundle 
$U \times \Delta^1 \to U$  and  are inspired by the discussion in section
\ref{chern-weil-review1}, we can equivalently conceive of $A$
as a morphism of dg-algebras
$$
  A: \mathrm{W}(\mathfrak{g}) \to \Omega^\bullet(U \times \Delta^1)
$$
from the Weil algebra of $\mathfrak{g}$
into the de Rham algebra of differential forms on the total space of the
interval bundle. It makes sense to decompose $A$ as the sum of a horizontal
1-form $A_U$ and a vertical 1-form $\lambda \, d t$, where $t\colon \Delta^1\to\mathbb{R}$ is
the canonical coordinate on $\Delta^1$: 
$$
  A = A_U + \lambda\, dt
  \,.
 $$ 
 The vertical part $A_\mathrm{vert}=\lambda\, dt$ of $A$ is an element of the completed tensor product $C^\infty(U)\hat{\otimes}\Omega^1(\Delta^1,\mathfrak{g})$ and can be seen as a family of $\mathfrak{g}$-connections on a trivial $G$-principal bundle on $\Delta^1$, parametrized by $U$. At any fixed $u_0\in U,$ the 1-form $\lambda(u_0,t)\,dt\in \Omega^1(\Delta^1,\mathfrak{g})$ satisfies the Maurer-Cartan equation by trivial dimensional reasons, and so we have a commutative diagram
\[
\xymatrix{
\,\phantom{mm}\Omega^\bullet(U\times\Delta^1)_\mathrm{vert}&\mathrm{CE}(\mathfrak{g})\ar[l]_{\phantom{mm}A_\mathrm{vert}}\\
\Omega^\bullet(U\times \Delta^1)\ar[u]&\mathrm{W}(\mathfrak{g})\ar[l]_{A}\ar[u]
}\,
\]
By the discussion in
section \ref{chern-weil-review1}, this can be seen as a \emph{first Cartan-Ehresmann
condition in the $\Delta^1$-direction}; it precisely encodes the fact that the 1-form $A$ on the total space of $U\times \Delta^1\to U$ is flat in the vertical direction.
\par
The curvature 2-form of $A$ decomposes as
$$
  F_{A} = F_{A_U} + F_{\Delta^1}
  \,,
$$
where the first term is at each point $t \in \Delta^1$ the ordinary curvature 
$F_{A_U} = d_U A_U + \frac{1}{2}[A_U,A_U]$ of $A_U$ at fixed $t\in\Delta^1$ and where the second term is 
$$
   F_{\Delta^1} = \left( d_U \lambda + [A_U,\lambda] -\frac{\partial}{\partial t} A_U\right)\wedge dt
  \,.
$$
We shall require that $F_{\Delta^1} = 0$; this is the \emph{second Ehresmann condition in the $\Delta^1$-direction} .
It implies that we have a commutative diagram
\[
\xymatrix{
\Omega^\bullet(\Delta^1\times U)&\mathrm{W}(\mathfrak{g})\ar[l]_{A}\\
\Omega^\bullet(U)\ar[u]& \mathrm{inv}(\mathfrak{g})\ar[l]_{F_A}\ar[u]
}\,.
\]
The condition $F_{\Delta^1}=0$ is equivalent to the differential equation
\[
\frac{\partial}{\partial t} A_U=d_U \lambda + [A_U,\lambda],
\]
 whose unique solution for given boundary condition $A_U\vert_{t = 0}$
specifies $A_U\vert_{t = 1}$ by the formula 
$$
  A_U(1) = g^{-1} A_U(0) g + g^{-1} d g
  \,, 
$$
where 
$$
  g := \mathcal{P} \exp(\int_{\Delta^1} \lambda d t) : U \to G
$$
is, pointwise in $U,$ the parallel transport of $\lambda dt$ along the interval. 
We may think of this as exhibiting formula
(\ref{gaugetransformation}) for gauge transformations
as arising from \emph{Lie integration} of infinitesimal data.
\par
Globalizing this local picture of connections on trivial bundles and gauge transformations between them now amounts to the following. For any (smooth, paracompact) manifold $X,$ we may find a \emph{good} open cover
$\{U_i \to X\}$, i.e., an open cover such that every non-empty $n$-fold intersection
$U_{i_1} \cap \cdots \cap U_{i_n}$ for all $n \in \mathbb{N}$ is diffeomorphic to
a Cartesian space. The cocycle data for a $G$-bundle with connection relative to this cover
is in degree 0 and 1 given by diagrams
\begin{equation}\label{first-diagram}
\raisebox{38pt}{\xymatrix{
0&\mathrm{CE}(\mathfrak{g})\ar[l]_{\phantom{mm}A_\mathrm{vert}}\\
\Omega^\bullet(U_{i})\ar[u]&\mathrm{W}(\mathfrak{g})\ar[l]_{A}\ar[u]\\
\Omega^\bullet(U_{i})\ar[u]& \mathrm{inv}(\mathfrak{g})\ar[l]_{F_A}\ar[u]
}}
\qquad\text{ and }\qquad 
\raisebox{38pt}{\xymatrix{
\,\phantom{m}\Omega^\bullet(\Delta^1\times U_{ij})_\mathrm{vert}&\mathrm{CE}(\mathfrak{g})\ar[l]_{\phantom{mmmm}A_\mathrm{vert}}\\
\Omega^\bullet(\Delta^1\times U_{ij})\ar[u]&\mathrm{W}(\mathfrak{g})\ar[l]_{A}\ar[u]\\
\Omega^\bullet(U_{ij})\ar[u]& \mathrm{inv}(\mathfrak{g})\ar[l]_{F_A}\ar[u]
}}
\,,
\end{equation}
where the latter restricts to the former after pullback along the two inclusions $U_{ij}\to U_i,U_j$ 
and along the face maps $\Delta^0=\{*\} \rightrightarrows \Delta^1$. This gives  
a collection of 1-forms $\{A_i \in \Omega^1(U_i,\mathfrak{g})\}_i$ and of smooth function $\{g_{i j} \in C^\infty(U_i \cap U_j, G)\}$, such that 
the formula 
$$
  A_j = g_{i j}^{-1} A_i g_{i j} + g_{i j}^{-1} d g_{i j}
$$
for gauge transformation holds on each double intersection
$U_i \cap U_j$.
This is almost the data defining a $\mathfrak{g}$-connecton on a $G$-principal bundle $P\to X$, but not quite yet, since it does not yet constrain the transition functions $g_{ij}$ on the triple intersections $U_i \cap U_j \cap U_k$ to obey the cocycle relation $g_{i j} g_{j k} = g_{i k}$. But since each $g_{ij}$ is the parallel transport of our connection along a vertical 1-simplex, the cocycle condition precisely says that parallel transport along the three edges of a vertical 2-simplex is trivial, i.e. that the  vertical parts of our connection forms on $U_{ijk}\times \Delta^1$ are the boundary data of a connection form on $U_{ijk}\times \Delta^2$ which is \emph{flat} in the vertical direction. In other words, the collection of commutative diagrams (\ref{first-diagram}) is to be seen as the 0 and 1-simplices of a simplicial set whose 2-simplices are the commutative diagrams
\[
\xymatrix{
\,\phantom{m}\Omega^\bullet(\Delta^2\times U_{ijk})_\mathrm{vert}&\mathrm{CE}(\mathfrak{g})\ar[l]_{\phantom{mmmm}A_\mathrm{vert}}\\
\Omega^\bullet(\Delta^2\times U_{ijk})\ar[u]&\mathrm{W}(\mathfrak{g})\ar[l]_{A}\ar[u]\\
\Omega^\bullet(U_{ijk})\ar[u]& \mathrm{inv}(\mathfrak{g})\ar[l]_{F_A}\ar[u]
}\,.
\]
Having added 2-simplices to our picture, we have finally recovered 
the standard description of connections in terms of local differential form
data. By suitably replacing Lie algebras with $L_\infty$-algebras in this derivation,
we will obtain a definition of connections on higher bundles in Section \ref{section.principal_infinity-bundles}. As one can expect, in the simplicial description of connections on higher bundles, simplices of arbitrarily high dimension will appear.

\section{Smooth $\infty$-groupoids}
\label{section.Lie_infinity-groupoids}

In this section we introduce a central concept that we will be dealing with in this paper, smooth $\infty$-groupoids, as a natural generalization of the classical notion of Lie groups. 

A Lie groupoid is, by definition,  a groupoid internal to the category of smooth spaces and smooth maps.
It is a widely appreciated fact in Lie groupoid theory that many features of Lie groupoids can be usefully  thought of in terms of their associated groupoid-valued presheaves on the category of manifolds, called the \emph{differentiable stack} represented by the Lie groupoid. This is the perspective that immediately generalizes to higher groupoids. 

Since many naturally appearing smooth spaces  are not manifolds 
-- particularly the spaces $[\Sigma,X]$ of smooth maps $\Sigma \to X$ between two manifolds -- for the development
of the general theory it is convenient to adopt a not too strict notion of `smooth space' . This generalized notion will have to be more flexible than the notion of manifold but at the same time not too far from that. 
The basic example to have in mind is the following: every smooth manifold $X$ of course represents a sheaf 
\begin{align*}
X:\text{SmoothManifolds}^{\mathrm {op}}&\to \text{Sets}\\
U &\mapsto C^\infty(U,X).
\end{align*}
on the category of smooth manifolds. But since manifolds themselves are by definition glued from Cartesian spaces $\mathbb{R}^n$, all the information about $X$ is in fact already encoded in the restriction of this sheaf to the category of Cartesian spaces and smooth maps between them:
$$
  X : \mathrm{CartSp}^\mathrm{op} \to \mathrm{Sets}
  \,.
$$
Now notice that also the spaces $[\Sigma,X]$ of smooth maps 
$\Sigma \to X$ between two manifolds naturally exist as sheaves on $\mathrm{CartSp}$, 
given by the assignment
$$
  [\Sigma,X] : U \mapsto C^\infty(\Sigma \times U, X)
  \,.
$$
Sheaves of this form are examples of generalized smooth spaces that are known as \emph{diffeological spaces} or \emph{Chen smooth spaces}. While not manifolds, these smooth spaces do have an underlying topological space and behave like smooth manifolds in many essential ways.

Even more generally, we will need to consider also `smooth spaces'  that do not  have even an underlying topological space. The central example of such is the sheaf of (real valued) closed differential $n$-forms
$$
  U \mapsto \Omega^n_{cl}(U)
  \,,
$$
\noindent
which we will need to consider later in the paper. We may think of these as modelling a kind of smooth Eilenberg-MacLane space that support a single (up to scalar multiple) smooth closed $n$-form. A precise version of this statement will play a central role later in the theory of Lie $\infty$-integration that we will describe in Section \ref{ooLieIntegration}.

Thus we see  that the common feature of generalized smooth spaces is not that they are \emph{representable} in one way or other. Rather, the common feature is that they all define sheaves on the category of the archetypical smooth spaces: the Cartesian spaces. This is a special case of an old insight going back to Grothendieck, Lawvere and others: with a category $\mathcal{C}$ of test spaces fixed, the correct context in which to consider generalized spaces modeled on $\mathcal{C}$ is the category $\mathrm{Sh}(\mathcal{C})$ of \emph{all} sheaves on $\mathcal{C}$:
the sheaf topos \cite{johnstone}. In there we may find a hierachy of types of generalized spaces ranging from ones that are very close to being like these test spaces, to ones that are quite a bit more general. In applications, it is good to find models as close as possible to the test spaces, but for the development of the theory it is better to admit them all.

Now if the manifold $X$ happens, in addition, to  be equipped with the structure of a Lie group $G$, then it represents more than just an ordinary sheaf of sets: from each group we obtain a simplicial set, its \emph{nerve}, whose set of $k$-cells is the set of $k$-tuples of elements in the group, and whose face and degeneracy maps are built from the product operation and the neutral element in the group. Since, for every $U \in \mathrm{CartSp}$, also the set of functions $C^\infty(U,G)$ forms a group, this means that from a Lie group we obtain a \emph{simplicial presheaf}
$$
  \mathbf{B}G : U \mapsto 
   \left\{
  \xymatrix{\cdots\ar@<6pt>[r] \ar@<2pt>[r] \ar@<-2pt>[r] \ar@<-6pt>[r]&
   C^\infty(U, G \times G)\ar@<4pt>[r] \ar[r] \ar@<-4pt>[r] &
   C^\infty(U,G)\ar@<2pt>[r]\ar@<-2pt>[r]& {*}
   },
  \right\}
  \,,
$$
where the degeneracy maps have not been displayed in order to make the diagram more readable.
\par
The simplicial presheaves arising this way are, in fact, special examples of presheaves taking values in \emph{Kan complexes}, i.e., in simplicial sets in which every \emph{horn} -- 
a simplex minus its interior and minus one face -- has a completion to a simplex; see for instance \cite{goerss-jardine} for a review. It turns out 
(see section 1.2.5 of \cite{lurie}) that Kan complexes may be thought of as modelling 
$\infty$-groupoids: the generalization of groupoids where one has not only morphisms between objects,
but also 2-morphisms between morphisms and generally $(k+1)$-morphisms between $k$-morphisms for all
$k \in \mathbb{N}$.
 The traditional theory of Lie groupoids may be thought of as dealing with those simplicial presheaves on $\mathrm{CartSp}$ that arise from nerves of Lie groupoids in the above manner

This motivates the definitions that we now turn to.

\subsection{Presentation by simplicial presheaves}
\label{section.Presentation_by_simplicial_presheaves}

\begin{definition} 
  A \emph{smooth  $\infty$-groupoid} $A$ is a simplicial  presheaf on the category
  $\mathrm{CartSp}$ of Cartesian spaces and smooth maps between them such that, over each $U \in \mathrm{CartSp},$
 $A$  is a Kan complex.
\end{definition} 
Much of ordinary Lie theory lifts from Lie groups to this context.
The reader is asked to keep in mind that smooth $\infty$-groupoids are objects whose smooth structure may be considerably more general than that of a Kan complex internal to smooth manifolds, i.e., of a simplicial smooth manifold satisfying a horn filling condition. Kan complexes internal to smooth manifolds, such as for instance nerves of ordinary Lie groupoids, can be thought of as representable smooth $\infty$-groupoids.
\begin{example}
The basic example of a representable smooth $\infty$-groupoids are ordinary Lie groupoids; in particular smooth manifolds and Lie groups are smooth $\infty$-groupoids. A particularly important example of representable smooth $\infty$-groupoid
is the  \emph{{\Cech} $\infty$-groupoid}: for $X$ a smooth manifold and $\mathcal{U}=\{U_i \to X\}$ an open cover, there is the simplicial manifold
\[
  \check{C}(\mathcal{U})
   :=
   \left\{
\xymatrix{\cdots\ar@<6pt>[r] \ar@<2pt>[r] \ar@<-2pt>[r] \ar@<-6pt>[r]&
\{U_{ijk}\}\ar@<4pt>[r] \ar[r] \ar@<-4pt>[r] &
\{U_{ij}\}\ar@<2pt>[r]\ar@<-2pt>[r]& \{U_i\}
}
  \right\}
\]
which in degree $k$ is the disjoint union of the $k$-fold intersections $U_i \cap U_j \cap \cdots$ of open subsets (the degeneracy maps are not depicted). This is a Kan complex internal to smooth manifolds in the evident way. 
\end{example}
While the notion of simplicial presheaf itself is straightforward, the correct concept of morphism between them is more subtle: we need a notion of morphisms such that the resulting category 
--or $\infty$-category as it were --
of our smooth $\infty$-groupoids reflects the prescribed notion of gluing of test objects. In fancier words,
we want simplicial presheaves to be equivalent to a higher analog of a sheaf topos: 
an $\infty$-topos \cite{lurie}. This may be achieved by equipping the naive category of simplicial presheaves with a \emph{model category structure} \cite{hovey}. This provides the information as to which objects in the category are to be regarded as equivalent, and how to resolve objects by equivalent objects for purposes of mapping between them. 
\noindent
There are some technical aspects to this this that we have relegated to the appendix. 
For all details and proofs of the definitions and propositions, respectively, 
in the remainder of this section see there.
\medskip
\par
\begin{definition}
Write $[\mathrm{CartSp}^{\mathrm{op}}, \mathrm{sSet}]_{\mathrm{proj}}$ for the global
projective model category structure on simplicial presheaves: weak equivalences and
fibrations are objectwise those of simplicial sets.
\end{definition}
This model structure presents the $\infty$-category of \emph{$\infty$-presheaves} on $\mathrm{CartSp}$. 
We impose now an $\infty$-sheaf condition.
\begin{definition}
\label{TheInfinityTopos}
Write $[\mathrm{CartSp}^{\mathrm{op}}, \mathrm{sSet}]_{\mathrm{proj},\mathrm{loc}}$ 
for the left Bousfield localization (see for instance section A.3 of \cite{lurie})
of $[\mathrm{CartSp}^{\mathrm{op}}, \mathrm{sSet}]_{\mathrm{proj}}$ at 
the set of all {\Cech} nerve projections $\check{C}(\mathcal{U}) \to U$ for $\mathcal{U}$
a differentiably good open cover of $U$, i.e., an open cover $\mathcal{U} = \{U_i \to U\}_{i \in I}$ of $U$ such that for all $n \in \mathbb{N}$ every $n$-fold intersection $U_{i_1} \cap \cdots \cap U_{i_n}$
is either empty or \emph{diffeomorphic} to $\mathbb{R}^{\dim U}$.
\end{definition}
This is the model structure that presents the $\infty$-category of \emph{$\infty$-sheaves}
or \emph{$\infty$-stacks} on $\mathrm{CartSp}$.
By standard results, it 
is a simplicial model category with respect to the canonical simplicial enrichment of 
simplicial presheaves, see \cite{dugger}. For $X, A$ two simplicial presheaves, we write
\begin{itemize}
  \item 
   $[\mathrm{CartSp}^{\mathrm{op}}, \mathrm{sSet}](X,A) \in \mathrm{sSet}$ 
      for the simplicial hom-complex of morphisms;
  \item
    $\mathbf{H}(X,A) := [\mathrm{CartSp}^{\mathrm{op}}, \mathrm{sSet}](Q(X), P(A))$ 
    for the right derived hom-complex (well defined up to equivalence) where $Q(X)$ 
    is any local cofibrant
    resolution of $X$ and $P(A)$ any local fibrant resolution of $A$.
\end{itemize}
\medskip
\par
Notice some standard facts about left Bousfield localization:
\begin{itemize}
  \item every weak equivalence in $[\mathrm{CartSp}^{\mathrm{op}}, \mathrm{sSet}]_{\mathrm{proj}}$
    is also a weak equivalence in $[\mathrm{CartSp}^{\mathrm{op}}, \mathrm{sSet}]_{\mathrm{proj}, \mathrm{loc}}$;
  \item
    the classes of cofibrations in both model structures coincide.
  \item 
    the fibrant objects of the local structure are precisely the objects that are fibrant in the 
    global structure and in addition satisfy \emph{descent} over all differentiably good open covers
    of Cartesian spaces. What this means precisely is stated in corollary
    A.\ref{CharacterizationOfLocalFibrancy} in the appendix.
   \item
     the localization right Quillen functor 
     $$\mathrm{Id} : [\mathrm{CartSp}^{\mathrm{op}}, \mathrm{sSet}]_{\mathrm{proj}}
      \to  [\mathrm{CartSp}^{\mathrm{op}}, \mathrm{sSet}]_{\mathrm{proj},\mathrm{loc}}$$
     presents 
     $\infty$-sheafification, 
     which is a left adjoint left exact $\infty$-functor \cite{lurie}, therefore
     all homotopy colimits and all finite homotopy limits in the local 
     model structure can be computed in the global model structure.
\end{itemize}
In particular, notice that an acyclic fibration in the global model structure will not, in general, be an acyclic fibration in the local model structure; nevertheless, it will be a weak equivalence in the local model structure.
\begin{definition}
We write
\begin{itemize}
  \item $\xymatrix{\ar[r]^\simeq &}$ for isomorphisms of simplicial presheaves;
  \item $\xymatrix{\ar[r]^\sim &}$ for weak equivalences in the global model structure;
  \item $\xymatrix{\ar[r]^{\sim_{\mathrm{loc}}} &}$ for weak equivalences in the local model structure;
\end{itemize}
(Notice that each of these generalizes the previous.)
\begin{itemize}
  \item $\xymatrix{\ar@{->>}[r] &}$ for fibrations in the global model structure.
\end{itemize}
We do not use notation for fibrations in the local model structure.
\end{definition}
Since the category $\mathrm{CartSp}$ has fewer objects than the category of all manifolds, we have that 
the conditions for simplicial presheaves to be fibrant in $[\mathrm{CartSp}^{\mathrm{op}}, \mathrm{sSet}]_{\mathrm{proj},\mathrm{loc}}$
are comparatively weak. For instance 
$$
\mathbf{B}G : U \mapsto N((C^\infty(U,G) \rightrightarrows *)
$$
is locally fibrant over $\mathrm{CartSp}$ but not over the site of all manifolds. This is discused below in section
\ref{section.Lie_infinity-groupoids_Examples}.
Conversely, the condition to be cofibrant is stronger over $\mathrm{CartSp}$ than it is over all manifolds. 
But by a central result by \cite{dugger}, we have fairly good control over cofibrant resolutions: these include 
notably {\Cech} nerves $\check{C}(\mathcal{U})$ of \emph{differentiably good} open covers, i.e.,  the {\Cech} nerve
$\check{C}(\mathcal{U}) \to X$ of a differentiably good open cover over a paracompact smooth manifold
$X$ is a cofibrant resolution of $X$ in $[\mathrm{CartSp}^{\mathrm{op}}, \mathrm{sSet}]_{\mathrm{proj},\mathrm{loc}}$,
and so we write
$$
  \check{C}(\mathcal{U}) \xrightarrow{\sim_{\mathrm{loc}}}X
  \,.
$$
Notice that in the present article these will be the only local weak equivalences that are not global
weak equivalences that we need to consider. 

In the practice of our applications, all this means that much of the technology hidden in definition 
\ref{TheInfinityTopos}
boils down to a simple  algorithm: after solving the comparatively easy tasks of finding a 
version $A$ of a given smooth $\infty$-groupoid that is fibrant over $\mathrm{CartSp}$, 
for describing morphisms of smooth $\infty$-groupoids from 
a manifold $X$ to  $A,$ we are to choose a differentiably good open cover $\mathcal{U} = \{U_i \to X\}$, 
form the {\Cech} nerve simplicial presheaf
$\check{C}(\mathcal{U})$ and then consider spans of ordinary morphisms of simplicial presheaves of the form
$$
  \xymatrix{
     \check{C}(\mathcal{U}) \ar[r]^g \ar[d]^{\begin{turn}{270}$\scriptstyle{\sim_\mathrm{loc}}$\end{turn}}& A
     \\
     X     
  }
  \,.
$$
Such a diagram of simplicial presheaves presents an object in $\mathbf{H}(X,A)$, in the hom-space of
the $\infty$-topos of smooth $\infty$-groupoids. As discussed below in section 
\ref{section.Lie_infinity-groupoids_Examples}, 
here the morphisms $g$ are naturally identified with cocycles in nonabelian 
{\Cech} cohomology on $X$ with coefficients in $A$. In section \ref{section.principal_infinity-bundles},
we discus that we may also think of these cocycles as transition data for $A$-principal $\infty$-bundles on $X$
\cite{SSSIII}. 

For discussing the $\infty$-Chern-Weil homomorphism, we are crucially interested in
composites of such spans: a \emph{characteristic map} on a coefficient object $A$ is 
nothing but a morphism $\mathbf{c} : A\to B$ in the $\infty$-topos, presented itself by a span
$$
  \xymatrix{
    \hat A \ar@{->>}[d]^\wr \ar[r] & B
    \\
    A
  }
  \,.
$$
The evaluation of this characteristic map on the $A$-principal bundle on $X$ encoded by a cocycle $g :\check{C}(\mathcal{U}) \to A$ is the composite morphism $X\to A \to B$ in the $\infty$-topos, which is presented by the composite span of simplicial presheaves
$$
  \xymatrix{
    Q X \ar@{->>}[d]^{\wr}\ar[r] 
      & \hat A \ar@{->>}[d]^\wr \ar[r] & B
    \\
    \check{C}(\mathcal{U})
    \ar[r]
        \ar[d]^{\begin{turn}{270}$\scriptstyle{\sim_\mathrm{loc}}$\end{turn}}
    &
    A
    \\
    X
  }
  \,.
$$
Here $Q X \to \check{C}(\mathcal{U})$ is the pullback of the acyclic fibration $\hat A \to A$, hence itself
an acyclic fibration; moreover, since $\check{C}(\mathcal{U})$ is cofibrant, we are guaranteed that a section $\check{C}(\mathcal{U}) \to Q X$ exists and is unique up to homotopy. Therefore the composite morphism $X\to A\to B$ is encoded in a cocylce $\check{C}(\mathcal{U}) \to B$ as in the diagram below:
$$
  \xymatrix{
    Q X \ar@{->>}[d]^{\wr}\ar[r] 
      & \hat A \ar@{->>}[d]^\wr \ar[r] & B
    \\
    \check{C}(\mathcal{U})
    \ar[r]
    \ar@/^1pc/[u]
    \ar[d]^{\begin{turn}{270}$\scriptstyle{\sim_\mathrm{loc}}$\end{turn}}
    &
    A
    \\
    X
  }
  \,.
$$
\par
Our main theorems will involve the construction of such span composites.

\medskip

\subsection{Examples}
\label{section.Lie_infinity-groupoids_Examples}

In this paper we will consider three main sources of smooth $\infty$-groupoids
\begin{itemize}
\item Lie groups and Lie groupoids, leading to Kan complexes via their nerves; examples of this kind will be the smooth $\infty$-groupoid $\mathbf{B}G$ associated with a Lie group $G$, and its refinements 
$\mathbf{B}G_{\mathrm{diff}}$ and 
$\mathbf{B}G_{\mathrm{conn}}$;
\item complexes of abelian groups concentrated in nonnegative degrees, leading to Kan complexes via the Dold-Kan correspondence; examples of this kind will be the smooth $\infty$-groupoid $\mathbf{B}^nU(1)$ associated with the chain complex of abelian groups consisting in $U(1)$ concentrated in degree $n$, and its refinements $\mathbf{B}^nU(1)_{\mathrm{diff}}$ and $\mathbf{B}^nU(1)_{\mathrm{conn}}$;
\item Lie algebras and $L_\infty$-algebras, via flat connections over simplices; this construction will produce, for any Lie or $L_\infty$-algebra $\mathfrak{g}$, a smooth $\infty$-groupoid $\exp_\Delta(\mathfrak{g})$ integrating $\mathfrak{g}$; other examples of this kind are the refinements $\exp_\Delta(\mathfrak{g})_{\mathrm{diff}}$ and $\exp_\Delta(\mathfrak{g})_{\mathrm{conn}}$ of $\exp_\Delta(\mathfrak{g})$. 
\end{itemize}
In the following sections, we will investigate these examples and show how they naturally combine in $\infty$-Chern-Weil theory.

\noindent{\bf Smooth $\infty$-groups.} With a useful notion of \emph{smooth $\infty$-groupoids} and their morphisms thus established, we automatically obtain a good notion of \emph{smooth $\infty$-groups}. This is accomplished simply by following the general principle by which essentially all basic constructions and results familiar from classical homotopy theory lift from the 
archetypical $\infty$-topos
$\mathrm{Top}$ of (compactly generated) topological spaces (or, equivalently, of \emph{discrete} $\infty$-groupoids) to any other $\infty$-topos, such as our
$\infty$-topos $\mathbf{H}$ of smooth $\infty$-groupoids.
\par
Namely, in classical homotopy theory a \emph{monoid} up to higher coherent homotopy 
is a topological space $X \in \mathrm{Top} \simeq \infty \mathrm{Grpd}$ equipped with 
 $A_\infty$-structure \cite{stasheff-H-spaces} or, equivalently, an $E_1$-structure, i.e., a
homotopical action of the little 1-cubes operad \cite{May}. A \emph{groupal} $A_\infty$-space -- 
an \emph{$\infty$-group} --
is one where this homotopy-associative product is invertible, up to homotopy. Famously, \emph{May's recognition theorem}
identifies such $\infty$-groups as being precisely, up to weak homotopy equivalence,
loop spaces. This establishes an equivalence of \emph{pointed} connected spaces with 
$\infty$-groups, given by looping $\Omega$ and delooping $B$:
$$
   \xymatrix{
      \infty \mathrm{Grp}
        \ar@<-4pt>[r]_{B}^\simeq
        \ar@{<-}@<+4pt>[r]^{\Omega}
      &
      \infty \mathrm{Grpd}_*
   }
  \,.
$$
Lurie shows in section 6.1.2 of \cite{lurie} (for $\infty$-groups) and in 
theorem 5.1.3.6 of \cite{LurieAlgebra}  
that these classical statements have 
direct analogs in any $\infty$-topos. We are thus entitled to think of any (pointed) connected smooth $\infty$-groupoid $X$
as the delooping $\mathbf{B}G$ of a smooth $\infty$-group $G \simeq \Omega X$
$$
   \xymatrix{
      \mathrm{Smooth}\infty \mathrm{Grp}
        \ar@<-4pt>[r]_<<<<{\mathbf{B}}^<<<<\simeq
        \ar@{<-}@<+4pt>[r]^<<<<<{\Omega}
      &
      \mathbf{H}_* = \mathrm{Smooth}\infty \mathrm{Grpd}_*
   }
  \,,
$$
where we use boldface $\mathbf{B}$ to indicate that the
delooping takes place in the $\infty$-topos $\mathbf{H}$ of smooth $\infty$-groupoids.

The most basic example for this we have already seen above: for $G$ any Lie group
the Lie groupoid $\mathbf{B}G$ described above is precisely the delooping of $G$ 
-- not in $\mathrm{Top}$ but in our $\mathbf{H}$.

In this article most smooth $\infty$-groups $G$ appear in the form of their smooth delooping
$\infty$-groupoids $\mathbf{B}G$. Apart from Lie groups, the main examples that we consider
will be the higher line and circle Lie groups $\mathbf{B}^n U(1)$ and $\mathbf{B}^n \mathbb{R}$
that have arbitrary many delooping, as well as the nonabelian smooth 2-group
$\mathrm{String}$ and the nonabelian smooth 6-group $\mathrm{Fivebrane}$, which are smooth refinements of the higher connected covers of the $\mathrm{Spin}$-group.

\subsubsection{$\mathbf{B}G$, $\mathbf{B}G_{\mathrm{conn}}$ and principal $G$-bundles with connection}
\label{section.BG_BG-conn_and_principal_G-bundles_with_connection}

The standard example of a stack on manifolds is the classifying stack $\mathbf{B}G$ 
for $G$-principal bundles with $G$ a Lie group. As an illustration of our setup,
we describe what this looks like in terms of simplicial presheaves over the
site $\mathrm{CartSp}$. Then we discuss its differential refinements $\mathbf{B}G_{\mathrm{diff}}$ and $\mathbf{B}G_{\mathrm{conn}}$.
\begin{definition}
Let $G$ be a Lie group. The smooth $\infty$-groupoid $\mathbf{B}G$ 
is defined to associate to a Cartesian space $U$ the nerve of the action groupoid 
$*// C^\infty(U,G)$, i.e.,  of the one-object groupoid with $C^\infty(U,G)$ as its set of morphisms
and composition given by the product of $G$-valued functions.
\end{definition}
\begin{remark}
Often this object is regarded over the site of all manifolds, where it is just a pre-stack, hence not fibrant.
Its fibrant replacement over that site is the stack $G \mathrm{Bund} : \mathrm{Manfd}^{\mathrm{op}} \to 
\mathrm{Grpd}$ that sends a manifold to the groupoid of $G$-principal bundles over it. We may think instead
of $\mathbf{B}G$ as sending a space to just the \emph{trivial} $G$-principal bundle and its automorphisms.
But since the site of Cartesian spaces is smaller, we have:
\end{remark}
\begin{proposition}
  The object $\mathbf{B}G \in [\mathrm{CartSp}^{\mathrm{op}}, \mathrm{sSet}]_{\mathrm{proj},\mathrm{loc}}$ is fibrant.
\end{proposition}
On the other hand, over the site of manifolds, every manifold itself is cofibrant. 
This means that to compute the groupoid of $G$-bundles on a manifold $X$ in terms of morphisms of stacks
over all manifolds, one usually passes to the fibrant replacement $G \mathrm{Bund}$ of $\mathbf{B}G$, then
considers $\mathrm{Hom}(X,G \mathrm{Bund})$ and uses the 2-Yoneda lemma to
identify this with the groupoid $G \mathrm{Bund}(X)$ of principal
$G$-bundles on $X$.
When working over $\mathrm{CartSp}$ instead, the situation is the opposite: here $\mathbf{B}G$ is already
fibrant, but the manifold $X$ is in general no longer cofibrant! To compute the groupoid of $G$-bundles
on $X$, we pass to a cofibrant replacement of $X$ given according to 
proposition \ref{CechResolution} by the \v{C}ech nerve $\check{C}(\mathcal{U})$ of a differentiably good open cover and then
compute $\mathrm{Hom}_{[\mathrm{CartSp}^{\mathrm{op}}, \mathrm{sSet}]}(\check{C}(\mathcal{U}), \mathbf{B}G)$. 
To see that the resulting groupoid is again equivalent to $G \mathrm{Bund}(X)$ (and hence to prove
the above proposition by taking $X = \mathbb{R}^n$) one proceeds as follows:

The object $\check{C}(\mathcal{U})$ is equivalent to the homotopy colimit in 
$[\mathrm{CartSp}^{\mathrm{op}}, \mathrm{sSet}]_{\mathrm{proj}}$ over the simplicial diagram of its components
\begin{align*}
  \check{C}(\mathcal{U}) &\simeq
  \mathrm{hocolim}
  \left(
  \xymatrix{
    \cdots\ar@<6pt>[r] \ar@<2pt>[r] \ar@<-2pt>[r] \ar@<-6pt>[r]
    &
    \coprod_{i,j,k}U_{ijk}\ar@<4pt>[r] \ar[r] \ar@<-4pt>[r] 
    &
    \coprod_{i j}U_{ij}\ar@<2pt>[r]\ar@<-2pt>[r]& \coprod_i U_i
  }
  \right)\\
  &
  \simeq \int^{[k] \in \Delta} \Delta[k] \cdot \coprod_{i_0, \cdots, i_k} U_{i_0,\dots, i_k}
  \,.
\end{align*}
(Here in the middle we are notationally suppressing the degeneracy maps for readability and on the 
second line we display for the inclined reader the formal coend expression
that computes this homotopy colimit as a weighted colimit \cite{hovey}. The dot denotes the tensoring 
of simplicial presheaves over simplicial sets).
Accordingly $\mathrm{Hom}(\check{C}(\mathcal{U}), \mathbf{B}G)$ is the homotopy limit
\[
  \begin{aligned}
  \mathrm{Hom}(\check{C}(\mathcal{U}), \mathbf{B}G)
  &\simeq
  \mathrm{holim}
  \left(
  \xymatrix{
    \cdots&\{\mathbf{B}G(U_{ijk})\}\ar@<6pt>[l] \ar@<2pt>[l] \ar@<-2pt>[l] \ar@<-6pt>[l]&
   \{\mathbf{B}G(U_{ij})\}\ar@<4pt>[l] \ar[l] \ar@<-4pt>[l] &
   \{\mathbf{B}G(U_i)\}\ar@<2pt>[l]\ar@<-2pt>[l]
 }
 \right)
 \\
 &
 \simeq\int_{[k] \in \Delta} \mathrm{Hom}(\Delta[k], \prod_{i_0,\cdots, i_k} \mathbf{B}G(U_{i_0, \dots, i_k}))
 \end{aligned}
 \,.
\]
The last line tells us that 
an element $g : \check{C}(\mathcal{U}) \to \mathbf{B}G$ in this Kan complex is a diagram
\[
  \xymatrix{
    \vdots & \vdots
    \\
      \ar@<+10pt>[u]
      \ar@<+5pt>[u]
      \ar@<+0pt>[u]
      \ar@<-5pt>[u]
      \ar@<-10pt>[u]
    \Delta[2] 
      \ar[r]^<<<<<{g^{(2)}} 
      & 
      \ar@<+10pt>[u]
      \ar@<+5pt>[u]
      \ar@<+0pt>[u]
      \ar@<-5pt>[u]
      \ar@<-10pt>[u]
    \prod_{i,j,k}\mathbf{B}G(U_{ijk})
    \\
    \Delta[1] 
      \ar[r]^<<<<<{g^{(1)}} 
      \ar@<+5pt>[u]
      \ar@<+0pt>[u]
      \ar@<-5pt>[u]
      & 
    \prod_{i,j}\mathbf{B}G(U_{ij})
      \ar@<+5pt>[u]
      \ar@<+0pt>[u]
      \ar@<-5pt>[u]
   \\
    \Delta[0] 
      \ar[r]^<<<<<{g^{(0)}} 
      \ar@<-3pt>[u]
      \ar@<3pt>	[u]
      & \prod_{i}\mathbf{B}G(U_i)
      \ar@<-3pt>[u]
      \ar@<3pt>	[u]
  }
\]
of simplicial sets. This is a collection $(\{g_i\}, \{g_{i j}\}, \{g_{i j k}\}, \cdots)$, where
\begin{itemize}
  \item $g_i$ is a vertex in $\mathbf{B}G(U_i)$;
  \item $g_{i j}$ is an edge in $\mathbf{B}G(U_{i j})$;
  \item $g_{i j k}$ is a 2-simplex in $\mathbf{B}G(U_{i j k})$
  \item etc.
\end{itemize}
such that the $k$-th face of the $n$-simplex in $\mathbf{B}G(U_{i_0, \cdots, i_n})$ is the 
image of the \hbox{$(n-1)$}-simplex under the $k$-th face inclusion $\mathbf{B}G(U_{i_0, \cdots, \hat i_k , \cdots, i_n})
\to \mathbf{B}G(U_{i_0, \cdots, i_n})$. (And similarly for the coface maps, which we continue to disregard
for brevity.) This means that an element $g : \check{C}(\mathcal{U}) \to \mathbf{B}G$ is precisely an element of the set $\check{C}(\mathcal{U},\mathbf{B}G)$ of \emph{nonabelian {\Cech} cocycles} with coefficients in $\mathbf{B}G$.
Specifically, by definition of $\mathbf{B}G$, 
this reduces to
\begin{itemize}
\item a collection of smooth maps $g_{ij}:U_{ij}\to G$, for every pair of indices $i,j$;
\item the constraint $g_{ij}g_{jk}g_{ki}=1_G$ on $U_{ijk}$, for every $i,j,k$ (the \emph{cocycle} constraint).
\end{itemize}
These are manifestly the data of transition functions defining a principal $G$-bundle over $X$. 
\par
Similarly working out the morphisms (i.e., the 1-simplices) in $\mathrm{Hom}(\check{C}(\mathcal{U}),\mathbf{B}G)$, we find that
their components are collections $h_i : U_i \to G$ of smooth functions, such that 
$g'_{i j} = h_i^{-1} g_{i j} h_j$. These are precisely the gauge transformations between the $G$-principal
bundles given by the transition functions $(\{g_{i j}\})$ and $(\{g'_{i j}\})$. Since the
cover $\{U_i \to X\}$ is \emph{good}, it follows that 
we have indeed reproduced the groupoid of $G$-principal bundles
$$
  \mathrm{Hom}(\check{C}(\mathcal{U}), \mathbf{B}G)=\check{C}(\mathcal{U},
\mathbf{B}G)
  \simeq G \mathrm{Bund}(X)
  \,.
$$
Two cocycles define isomorphic principal $G$-bundles precisely when they define the same element in 
{\Cech} cohomology with coefficients in the sheaf of smooth functions with values in $G$. Thus we recover the standard fact that isomorphism classes of principal $G$-bundles are in natural bijection with $H^1(X,G)$.

We now consider a differential refinement of $\mathbf{B}G$.
\begin{definition}
Let $G$ be a Lie group with Lie algebra $\mathfrak{g}$. The smooth $\infty$-groupoid 
$\mathbf{B}G_{\mathrm{conn}}$ is defined to associate with a Cartesian space $U$ the nerve of the action groupoid $\Omega^1(U,\mathfrak{g})// C^\infty(U,G)$.
\end{definition}
This is over $U$ the groupoid $\mathbf{B}G_{\mathrm{conn}}(U)$
\begin{itemize}
  \item
    whose set of objects is the set of smooth $\mathfrak{g}$-valued 1-forms 
    $A \in \Omega^1(U,\mathfrak{g})$;
  \item
    whose morphisms $g: A \to A'$ are labeled by smooth functions $g \in C^\infty(U,G)$ 
    such that they relate the source and target by a \emph{gauge transformation}
   \[
     A' = g^{-1} A g + g^{-1}d g
     \,,
   \]
   where $g^{-1}Ag$ denotes pointwise the adjoint action of $G$ on $\mathfrak{g}$ 
   and $g^{-1}dg$ is the pull-back $g^*(\theta)$ of the Maurer-Cartan form $\theta\in\Omega^1(G,\mathfrak{g})$.    
\end{itemize}
With $X$ and $\check{C}(\mathcal{U})$ as before we now have:
\begin{proposition} The smooth $\infty$-groupoid 
$\mathbf{B}G_{\mathrm{conn}}$ is fibrant and there is a natural equivalence of groupoids
  $$
    \mathbf{H}(X, \mathbf{B}G_{\mathrm{conn}}) 
    \simeq G \mathrm{Bund}_{\mathrm{conn}}(X)
    \,,
  $$
  where on the right we have the groupoid of $G$-principal bundles on $X$ equipped with connection.
\end{proposition}
This follows along the above lines, by unwinding the nature of the simplicial hom-set $\check{C}(\mathcal{U}, \mathbf{B}G_{\mathrm{conn}}):=\mathrm{Hom}(\check{C}(\mathcal{U}), \mathbf{B}G_{\mathrm{conn}})$ of nonabelian {\Cech} cocycles
with coefficients in $\mathbf{B}G_{\mathrm{conn}}$. Such a cocycle is a collection 
$(\{A_i\}, \{g_{i j}\})$ consisting of
\begin{itemize}
\item a 1-form $A_i\in\Omega^1(U_i,\mathfrak{g})$ for each index $i$;
\item a smooth function $g_{ij}:U_{ij}\to G$, for all indices $i,j$;
\item the gauge action constraint $A_j=g_{ij}^{-1}A_ig_{ij}+g_{ij}^{-1}dg_{ij}$ on $U_{ij}$, for all indices $i,j$;
\item the cocycle constraint $g_{ij}g_{jk}g_{ki}=1_G$ on $U_{ijk}$, for all indices $i,j,k$.
\end{itemize}
These are readily seen to be the data defining a $\mathfrak{g}$-connection on a principal $G$-bundle over $X$. 

Notice that there is an evident ``forget the connection''-morphism 
$\mathbf{B}G_{\mathrm{conn}}\to \mathbf{B}G$, given over $U \in \mathrm{CartSp}$ by
$$
  (A \stackrel{g}{\to} A') \mapsto (\bullet \stackrel{g}{\to} \bullet)
  \,. 
$$
We denote the set of isomorphism classes of principal $G$-bundles with connection by the symbol $H^1(X,G)_{\mathrm{conn}}$. Thus we obtain a morphism
\[
H^1(X,G)_{\mathrm{conn}}\to H^1(X,G).
\]
Finally, we introduce a smooth $\infty$-groupoid $\mathbf{B}G_{\mathrm{diff}}$ in between 
$\mathbf{B}G$ and $\mathbf{B}G_{\mathrm{conn}}$. This may seem a bit curious, but we'll 
see in section \ref{section.infinity-stack_of_principal_bundles_with_connection} how it is the degree one case of a completely natural and noteworthy
general construction. Informally, $\mathbf{B}G_{\mathrm{diff}}$ is obtained from $\mathbf{B}G$ by freely decorating the vertices of the simplices in $\mathbf{B}G$ by elements in $\Omega^1(U,\mathfrak{g})$. More formally, we have the following definition.
\begin{definition}
Let $G$ be a Lie group with Lie algebra $\mathfrak{g}$. The smooth $\infty$-groupoid 
$\mathbf{B}G_{\mathrm{diff}}$ is defined to associate with a Cartesian space $U$ the 
nerve of the groupoid 
\begin{enumerate}
\item whose set of objects is $\Omega^1(U,\mathfrak{g})$;
\item a morphism $A \stackrel{(g,a)}{\to} A'$ is labeled by 
      $g\in C^\infty(U,G)$ and $a \in \Omega^1(U,\mathfrak{g})$ such that
      $$
        A = g^{-1} A' g + g^{-1} d  g + a
      $$
\item composition of morphisms is given by 
$$
(g,a)\circ(h,b)=(gh, h^{-1}ah+h^{-1}dh+b).
$$
\end{enumerate}
\end{definition}
\begin{remark}\label{remark.BG_diff} This definition intentionally carries an evident redundancy: given any $A$, $A'$ and
$g$ the element $a$ that makes the above equation hold does exist uniquely; the 1-form $a$ measures
the failure of $g$ to constitute a morphism from $A$ to $A'$ in $\mathbf{B}G_{\mathrm{conn}}$. 
We can equivalently express the redundancy of $a$ by saying that there is a natural isomorphism between $\mathbf{B}G_{\mathrm{diff}}$ and the direct product of $\mathbf{B}G$ with the codiscrete groupoid on the sheaf of sets $\Omega^1(-;\mathfrak{g})$.
\end{remark}
\begin{proposition}
The evident forgetful morphism $\mathbf{B}G_{\mathrm{conn}} \to \mathbf{B}G$ factors
through $\mathbf{B}G_{\mathrm{diff}}$ by a monomorphism followed by an acyclic fibration (in the global model structure)
\[
\mathbf{B}G_{\mathrm{conn}}\hookrightarrow\mathbf{B}G_{\mathrm{diff}}
\stackrel{\sim}{\twoheadrightarrow}\mathbf{B}G
  \,.
\]
\end{proposition}

\subsubsection{$\mathbf{B}G_2$, and nonabelian gerbes and principal 2-bundles}
\label{section.Principal_2_Bundles}

We now briefly dicuss the first case of $G$-principal $\infty$-bundles after ordinary 
principal bundles, the case where $G$ is a \emph{Lie 2-group}:
$G$-principal 2-bundles. 

When $G = \mathrm{AUT}(H)$ the \emph{automorphism 2-group} of a Lie group $H$ (see below)
these structures have the same classification (though are conceptually somewhat different from) the smooth version 
of the $H$-banded \emph{gerbe}s of \cite{Giraud} (see around def. 7.2.2.20 in \cite{lurie} for a conceptually
clean account in the modern context of higher toposes):
both are classified by the \emph{nonabelian cohomology} $H^1_{\mathrm{Smooth}}(-,\mathrm{AUT}(H))$ 
with coefficients in that 2-group. But the main examples of 2-groups that we shall be
interested in, namely \emph{string 2-groups}, are not equivalent to $\mathrm{AUT}(H)$ for
any $H$, hence the 2-bundles considered here are strictly more general than 
Giraud's gerbes. The literature knows what has been called \emph{nonabelian bundle gerbes},
but despite their name these are not Giraud's gerbes, but are instead models
for the total spaces of what we call here \emph{principal 2-bundles}. A good 
discussion of the various equivalent incarnations of principal 2-bundles is in 
\cite{NikolausWaldorf}.

To start with, note the general abstract notion of smooth 2-groups:
\begin{definition}
   A \emph{smooth 2-group} is a 1-truncated group object in 
   $\mathbf{H} = \mathrm{Sh}_{\infty}(\mathrm{CartSp})$. These are equivalently given by
   their (canonically pointed) delooping 2-groupoids $\mathbf{B}G \in \mathbf{H}$, which are
   precisely, up to equivalence, the connected 2-truncated objects of $\mathbf{H}$.
   
   For $X \in \mathbf{H}$ any object, $G 2\mathrm{Bund}_{\mathrm{smooth}}(X) :=
   \mathbf{H}(X, \mathbf{B}G)$ is the 2-groupoid of smooth $G$-principal 2-bundles on $G$.
\end{definition}
While nice and abstract, in applications one often has -- or can get -- hold of a 
\emph{strict model} of a given smooth 
2-group. The following definitions can be found recalled in any reference on these matters, for instance in 
\cite{NikolausWaldorf}. 
\begin{definition}
  \label{CrossedModuleAndStrict2Group}
  \begin{enumerate}
  \item A smooth \emph{crossed module} of Lie groups is a pair of homomorphisms 
  $\partial : G_1 \to G_0$ and $\rho : G_0 \to \mathrm{Aut}(G_1)$ of 
  Lie groups, such that for all $g \in G_0$ and $h,h_1, h_2 \in G_1$
  we have $\rho(\partial h_1)(h_2) = h_1 h_2 h_1^{-1}$ and
  $\partial \rho(g)(h) = g \partial(h) g^{-1}$.

  \item For $(G_1 \to G_0)$ a smooth crossed module, the corresponding \emph{strict Lie 2-group} is 
  the smooth groupoid 
  $\xymatrix{
    G_0 \times G_1 \ar@<+4pt>[r]_{}\ar@<-4pt>[r]^{} & G_0}$, whose source map is given by projection on 
    $G_0$, whose target map is given by applying $\partial$ to the second factor and then multiplying
    with the first in $G_0$, and whose composition is given by multiplying in $G_1$.
 
   This groupoid has a strict monoidal structure with strict inverses given by equipping 
   $G_0 \times G_1$ with the semidirect product group structure $G_0 \ltimes G_1$ induced
   by the action $\rho$ of $G_0$ on $G_1$.
 
  \item 
   The corresponding one-object strict smooth 2-groupoid we write $\mathbf{B}(G_1 \to G_0)$.
   As a simplicial object (under Duskin nerve of 2-categories)
  this is of the form
  $$
    \mathbf{B}(G_1 \to G_0) =
    \mathrm{cosk}_3
    \left(
     \xymatrix{
       G_0^{\times 3} \times G_1^{\times 3}\ar@<-6pt>[r]\ar[r]\ar@<+6pt>[r] & G_0^{\times 2} \times G_1 \ar@<-4pt>[r]\ar@<+4pt>[r] & G_0 \ar[r] & {*}
       }
    \right)
    \,.
  $$
  \end{enumerate}
\end{definition}
{\bf Examples.}
\begin{enumerate}
  \item
    For $A$ any abelian Lie group, $A \to 1$ is a crossed module. Conversely, 
    for $A$ any Lie group $A \to 1$ is a crossed module precisely if $A$ is abelian.
    We write $\mathbf{B}^2 A = \mathbf{B}(A \to 1)$. This case and its 
    generalizations is discussed below in \ref{section.Deligne_cohomology}.
  \item
    For $H$ any Lie group with automorphism Lie group $\mathrm{Aut}(H)$, the morphism
    $H \stackrel{\mathrm{Ad}}{\to} \mathrm{Aut}(H)$ that sends group elements to inner
    automorphisms, together with $\rho = \mathrm{id}$, is a crossed module.
    We write $\mathrm{AUT}(H) := (H \to \mathrm{Aut}(H))$ and speak of the 
    \emph{automorphism 2-group} of $H$, because this is $\simeq \mathrm{Aut}_{\mathbf{H}}(\mathbf{B}H)$.
  \item
    For $G$ an ordinary Lie group and $\mathbf{c} : \mathbf{B}G \to \mathbf{B}^3 U(1)$
    a morphism in $\mathbf{H}$ (see the \ref{section.Deligne_cohomology} for a discussion of 
    $\mathbf{B}^n U(1)$), its homotopy fiber $\mathbf{B}\hat G \to \mathbf{B}G$ is the
    delooping of a smooth 2-group $\hat G$. If $G$ is compact, simple and simply connected, then 
   this is equivalent (\cite{survey}, section 4.1) to a strict 2-group 
   $({\hat \Omega} G \to P G)$ given by a $U(1)$-central extension of the loop group of $G$,
   as described in \cite{bcss}. This is called the \emph{string 2-group} extension of $G$ by 
   $\mathbf{c}$. We come back to this in \ref{section.characteristic_maps_by_Lie_integration}.
\end{enumerate}
\begin{observation}
  \label{GlobalFibrancyOfBOf2Group}
  For every smooth crossed module, its delooping object $\mathbf{B}(G_1 \to G_0)$ is fibrant in $[\mathrm{CartSp}^{\mathrm{op}}, \mathrm{sSet}]$.
\end{observation}
\proof
  Since $(G_1 \to G_0)$ induces a strict 2-group, there are horn fillers defined by the smooth operations
  in the 2-group: we can always solve for the missing face in a horn in terms of an expression 
  involving the smooth composite-operations and inverse-operations in the 2-group.
\endofproof
\begin{proposition}
  \label{LocalFibrancyOfBOf2Group}
  Suppose that the smooth crossed module $(G_1 \to G_0)$ 
  is such that the quotient $\pi_0 G = G_0/G_1$ is a smooth manifold and the
  projection $G_0 \to G_0/G_1$ is a submersion.
  
  Then $\mathbf{B}(G_1 \to G_0)$ is fibrant also in 
  $[\mathrm{CartSp}^{\mathrm{op}}, \mathrm{sSet}]_{\mathrm{proj}, \mathrm{loc}}$.
\end{proposition}
\proof
  We need to show that for $\{U_i \to \mathbb{R}^n\}$ a good open cover, the canonical 
  descent morphism
  $$
    B(C^\infty(\mathbb{R}^n, G_1) \to C^\infty(\mathbb{R}^n, G_0))
    \to 
    [\mathrm{CartSp}^{\mathrm{op}}, \mathrm{sSet}](\check{C}(\mathcal{U}), \mathbf{B}(G_1 \to G_0))
  $$
  is a weak homotopy equivalence. 
  The main point to show is that, since the Kan complex on the left is connected by 
  construction, also the Kan complex on the right is. 

  To that end, notice that the category $\mathrm{CartSp}$ equipped with the open cover topology is a 
  \emph{Verdier site} in the sense of section 8 of \cite{dugger-hollander-isaksen}.
  By the discussion there it follows that every hypercover over $\mathbb{R}^n$ can be refined
  by a split hypercover, and these are cofibrant resolutions of $\mathbb{R}^n$ in both the 
  global and the local model structure
  $[\mathrm{CartSp}^{\mathrm{op}}, \mathrm{sSet}]_{\mathrm{proj}, \mathrm{loc}}$. 
  Since also $\check{C}(\mathcal{U}) \to \mathbb{R}^n$ is a cofibrant resolution and  since
  $\mathbf{B}G$ is fibrant in the \emph{global} structure by observation \ref{GlobalFibrancyOfBOf2Group}, 
  it follows from the
  existence of the global model structure that morphisms out of $\check{C}(\mathcal{U})$ into 
  $\mathbf{B}(G_1 \to G_0)$ capture all cocycles over any hypercover over $\mathbb{R}^n$,
  hence  that 
  $$
    \pi_0 [\mathrm{CartSp}^{\mathrm{op}}, \mathrm{sSet}](\check{C}(\mathcal{U}), \mathbf{B}(G_1 \to G_0))
    \simeq
    H^1_{\mathrm{smooth}}(\mathbb{R}^n, (G_1 \to G_0))
  $$
  is the standard {\v C}ech cohomology of $\mathbb{R}^n$, defined as a colimit over refinements of
  covers of equivalence classes of {\v C}ech cocycles.
  
  Now by prop. 4.1 of \cite{NikolausWaldorf} (which is the smooth refinement of the 
  statement of \cite{baezstevenson} in the continuous context)
  we have that under our assumptions on $(G_1 \to G_0)$ 
  there is a
  topological classifying space for this smooth {\v C}ech cohomology set. Since $\mathbb{R}^n$
  is topologically contractible, it follows that this is the
  singleton set and hence the above descent morphism is indeed an isomorphism on $\pi_0$.
  
  Next we can argue that it is also an isomorphism on $\pi_1$, by reducing to the analogous
  local trivialization statement for ordinary principal bundles: a loop in 
  $[\mathrm{CartSp}^{\mathrm{op}}, \mathrm{sSet}](\check{C}(\mathcal{U}), \mathbf{B}(G_1 \to G_0))$
  on the trivial cocycle
  is readily seen to be a $G_0 // (G_0 \ltimes G_1)$-principal groupoid bundle, over the
  action groupoid as indicated. The underlying $G_0 \ltimes G_1$-principal bundle
  has a trivialization on the contractible $\mathbb{R}^n$ (by classical results or, in fact,
  as a special case of the previous argument), and so equivalence classes of such loops are 
  given by $G_0$-valued smooth functions on $\mathbb{R}^n$. The descent morphism
  exhibits an isomorphism on these classes. 
  
  Finally the equivalence classes of spheres on both sides are directly seen to be 
  smooth $\mathrm{ker}(G_1 \to G_0)$-valued
  functions on both sides, identified by the descent morphism.
\endofproof
\begin{corollary}
  For $X \in \mathrm{SmoothMfd} \subset \mathbf{H}$ a paracompact smooth manifold, and
  $(G_1 \to G_0)$ as above, we have for any good open cover $\{U_i \to X\}$ that
  the 2-groupoid of smooth $(G_1 \to G_0)$-principal 2-bundles is
  $$
    (G_1 \to G_0)\mathrm{Bund}(X) := \mathbf{H}(X, \mathbf{B}(G_1))
    \simeq 
    [\mathrm{CartSp}^{\mathrm{op}}, \mathrm{sSet}](\check{C}(\mathcal{U}), \mathbf{B}(G_1 \to G_0))
  $$
  and its set of connected components is naturally isomorphic to the nonabelian {\v C}ech
  cohomology
  $$
    \pi_0 \mathbf{H}(X, \mathbf{B}(G_1 \to G_0)) \simeq H^1_{\mathrm{smooth}}(X, (G_1 \to G_0))
    \,.
  $$
\end{corollary}

\subsubsection{$\mathbf{B}^nU(1)$, $\mathbf{B}^nU(1)_{\mathrm{conn}}$, circle $n$-bundles and Deligne cohomology}
\label{section.Deligne_cohomology}

A large class of examples of smooth $\infty$-groupoids is induced from chain complexes of sheaves
of abelian groups by the Dold-Kan correspondence \cite{goerss-jardine}.
\begin{proposition}
The \emph{Dold-Kan correspondence} is an equivalence of categories
$$
  \xymatrix{
    \mathrm{Ch}_\bullet^+
    \ar@{<-}@<-3pt>[r]_{N_\bullet}
    \ar@<+3pt>[r]^{\mathrm{DK}}
    &
    \mathrm{sAb}    
  }  
  \,,
$$
between non-negatively graded chain complexes and simplicial abelian groups, where 
$N_\bullet$ forms the \emph{normalized chains complex} of a simplicial abelian group $A_\Delta$. 
Composed with the forgetful functor $\mathrm{sAb} \to \mathrm{sSet}$ and prolonged to a functor
on sheaves of chain complexes, the functor
$$
  \mathrm{DK} : [\mathrm{CartSp}^{op}, \mathrm{Ch}_\bullet^+] \to [\mathrm{CartSp}^{\mathrm{op}}, \mathrm{sSet}]
$$
takes degreewise surjections to fibrations and degreewise quasi-isomorphisms to weak equivalences
in $[\mathrm{CartSp}^{\mathrm{op}}, \mathrm{sSet}]_{\mathrm{proj}}$.
\end{proposition}
We will write an element $(A_\bullet,\partial)$ of $\mathrm{Ch}_\bullet^+$ as
\[
\cdots\to A_{k}\to A_{k-1}\to\cdots \to A_2\to A_1\to A_0
\]
and will denote by $[1]$ the ``shift on the left'' functor on chain complexes defined by $(A_\bullet[1])_k=A_{k-1}$, i.e., $A_\bullet[1]$ is the chain complex
\[
\cdots\to A_{k-1}\to \cdots \to A_2\to A_1\to A_0\to 0.
\]
\begin{remark}The reader used to cochain complexes, and so to the shift functor $(A^\bullet[1])^k=A^{k+1}$ could at first be surprised by the minus sign in the shift functor on chain complexes; but the shift rule is actually the same in both contexts, as it is evident by writing it as $(A_\bullet[1])_k=A_{k+\deg(\partial)}$.
\end{remark}
For $A$ any abelian group, we can consider $A$ as a chain complex concentrated in degree zero, and so $A[n]$ will be the chain complex consisting of $A$ concentrated in degree $n$.
\begin{definition}
Let $A$ be an abelian Lie group. Define the simplicial presheaf $\mathbf{B}^n A$ 
to be the image under $\mathrm{DK}$ of the sheaf of complexes $C^\infty(-,A)[n]$:
$$
  \mathbf{B}^n A
     :
     U 
   \mapsto
   \mathrm{DK}(C^\infty(U,A) \to 0 \to \cdots \to 0)
  \,,
$$
with $C^\infty(U,A)$ in degree $n$. Similarly,
for $K \to A$ a morphism of abelian groups, 
write $\mathbf{B}^n (K \to A)$ for the image under $\mathrm{DK}$ of the complex of sheaves of abelian groups
$$
  (C^\infty(-,K) \to C^\infty(-,A) \to 0 \to \cdots \to 0)
$$
with $C^\infty(-,A)$ in degree $n$; for $n \geq 1$ we write $\mathbf{E} \mathbf{B}^{n-1}A$ for $\mathbf{B}^{n-1}(A \stackrel{\mathrm{Id}}{\to} A)$. 
\end{definition}
\begin{proposition}
  For $n \geq 1$ the object $\mathbf{B}^n A$ is indeed the delooping of the object $\mathbf{B}^{n-1}A$.
\end{proposition}
\proof
  This means that there is an \emph{$\infty$-pullback} diagram \cite{lurie}
  $$
    \xymatrix{
      \mathbf{B}^{n-1}A \ar[d]\ar[r] & {*} \ar[d]
      \\
      {*} \ar[r] & \mathbf{B}^n A
    }
    \,.
  $$
  This is presented by the corresponding homotopy pullback in $[\mathrm{CartSp}^{\mathrm{op}},\mathrm{sSet}]$.
  Consider the diagram
  $$
    \xymatrix{
      \mathbf{B}^{n-1}A \ar[d]\ar[r]&  \mathbf{E}\mathbf{B}^{n-1}A \ar@{->>}[d]\ar[r]^{\phantom{mm}\sim}  & {*}
      \\
      {*}\ar[r] & \mathbf{B}^n A
    }
    \,,
  $$
  The right vertical morphism is a replacement of the point inclusion by a fibration and the 
  square is a pullback in $[\mathrm{CartSp}^{\mathrm{op}}, \mathrm{sSet}]$
  (the pullback of presheaves is computed objectwise and under the DK-correspondence
  may be computed in $\mathrm{Ch}_\bullet^+$, where it is evident). Therefore this exhibits
  $\mathbf{B}^{n-1}A$ as the homotopy pullback, as claimed.
\endofproof
\begin{proposition} \label{FibrancyOfBnA}
	For $A = \mathbb{Z}, \mathbb{R}, U(1)$ and all $n\geq 1$ we have that 
	$\mathbf{B}^n A$ satisfies descent over $\mathrm{CartSp}$ in that it is fibrant
	in $[\mathrm{CartSp}^{\mathrm{op}}, \mathrm{sSet}]_{\mathrm{proj}, \mathrm{loc}}$.
\end{proposition}
\proof
  One sees directly in terms of {\Cech} cocycles that the homotopy groups based at the trivial
  cocycle in the simplicial hom-sets $[\mathrm{CartSp}^{\mathrm{op}}, \mathrm{sSet}](\check{C}(\mathcal{U}), \mathbf{B}^nA)$
  and $[\mathrm{CartSp}^{\mathrm{op}}, \mathrm{sSet}](U, \mathbf{B}^nA)$ are naturally identified.
  Therefore it is sufficient to show that
  $$
    * \simeq 
    \pi_0[\mathrm{CartSp}^{\mathrm{op}}, \mathrm{sSet}](U, \mathbf{B}^nA)
      \to 
    \pi_0[\mathrm{CartSp}^{\mathrm{op}}, \mathrm{sSet}](\check{C}(\mathcal{U}), \mathbf{B}^nA)        
  $$
  is an isomorphism. This amounts to proving that the $n$-th \v{C}ech cohomology group of $U$ with coefficients in $\mathbb{Z}$, $\mathbb{R}$ or $U(1)$ is trivial, which is immediate since $U$ is contractible (for $U(1)$ one uses the isomophism $H^n(U,U(1))\simeq H^{n+1}(U,\mathbb{Z})$ in \v{C}ech cohomology).
\endofproof
\begin{definition}
	\label{definition.AbelianNBundles}
  For $X$ a smooth $\infty$-groupoid and $QX \to X$ a cofibrant replacement
 we say that 
  \begin{itemize}
  \item
    for $X \xleftarrow{\sim_\mathrm{loc}} QX \stackrel{g}{\to} \mathbf{B}^n A$
    a span in $[\mathrm{CartSp}^{\mathrm{op}}, \mathrm{sSet}]$, the corresponding 
    \emph{$(\mathbf{B}^{n-1})A$-principal $n$-bundle} is the $\infty$-pullback
    $$
      \xymatrix{
         P \ar[d]\ar[r] & {*}\ar[d]
         \\
         X \ar[r]^g & \mathbf{B}^n A
      }
    $$
    hence the ordinary pullback in $[\mathrm{CartSp}^{\mathrm{op}}, \mathrm{sSet}]$
    $$
      \xymatrix{
         P \ar[r] \ar[d] & \mathbf{E}\mathbf{B}^{n-1} A \ar[d]
         \\
         QX \ar[d]^{\begin{turn}{270}$\scriptstyle{\sim_\mathrm{loc}}$\end{turn}} \ar[r]^g & \mathbf{B}^n A
         \\
         X
      }
      \,.
    $$
  \item
  the Kan complex
  $$
    (\mathbf{B}^{n-1}A) \mathrm{Bund}(X)
    :=
    \mathbf{H}(X, \mathbf{B}^n A)    
  $$
  is the \emph{$n$-groupoid} of smooth $\mathbf{B}^{n-1}A$-principal $n$-bundles on $X$.
  \end{itemize}
\end{definition}

\begin{proposition}
  For $X$ a smooth paracompact manifold,
  the $n$-groupoid $(\mathbf{B}^{n-1} A)\mathrm{Bund}(X)$ is equivalent to the $n$-groupoid $\check{C}(\mathcal{U},  \mathbf{B}^n A)$ of degree $n$
  {\Cech} cocycles on $X$ with coefficients in the sheaf of smooth functions with values in $A$. In particular
  $$
    \pi_0 (\mathbf{B}^{n-1}A)\mathrm{Bund}(X)=\pi_0\mathbf{H}(X, \mathbf{B}^n A) \simeq H^n(X,A)
  $$
  is the {\Cech} cohomology of $X$ in degree $n$ with coefficients in $A$. \end{proposition}
\proof
  This follows from the same arguments as in the previous section given for 
  the more general nonabelian {\Cech} cohomology.
\endofproof

We will be  interested mainly in the abelian Lie group $A = U(1)$. The exponential exact sequence
$0\to \mathbb{Z}\to \mathbb{R}\to U(1)\to 1$ induces an acyclic fibration (in the global model structure) $\mathbf{B}^n(\mathbb{Z}\hookrightarrow\mathbb{R})\stackrel{\sim}{\twoheadrightarrow}\mathbf{B}^nU(1)$, and one has the
long fibration sequence 
\[
\xymatrix{&&\mathbf{B}^n(\mathbb{Z}\hookrightarrow\mathbb{R})\ar@{->>}^{\wr}[d]\ar[r] &\mathbf{B}^{n+1}\mathbb{Z}\cdots\\
\cdots\to\mathbf{B}^{n}\mathbb{Z}\ar[r]& \mathbf{B}^n\mathbb{R}\ar[r]& \mathbf{B}^n U(1)&
 }
\]
from which one recovers the classical isomorphism $H^n(X,U(1))\simeq H^{n+1}(X,\mathbb{Z})$. Next, we consider differential refinements of these cohomology groups.
\begin{definition}
The smooth $\infty$-groupoid $\mathbf{B}^nU(1)_{\mathrm{conn}}$ is the image via the Dold-Kan correspondence of the \emph{Deligne complex} $U(1)[n]_D^\infty$, i.e., of the chain complex of sheaves of abelian groups
\[
U(1)[n]_D^\infty:=\left(C^\infty(-,U(1))\xrightarrow{d\mathrm{log}}\Omega^1(-,\mathbb{R})\xrightarrow{d}\cdots\xrightarrow{d} \Omega^n(-,\mathbb{R})\right)
\]
concentrated in degrees $[0,n]$.
Similarly, the smooth $\infty$-groupoid $\mathbf{B}^{n}(\mathbb{Z}\hookrightarrow\mathbb{R})_{\mathrm{conn}}$ is the image via the Dold-Kan correspondence of the complex of sheaves of abelian groups
\[
\mathbb{Z}[n+1]^{\infty}_D:=\left(\mathbb{Z}\hookrightarrow C^\infty(-,\mathbb{R})\xrightarrow{d}\Omega^1(-,\mathbb{R})\xrightarrow{d}\cdots\xrightarrow{d} \Omega^n(-,\mathbb{R})\right),
\]
concentrated in degrees $[0,n+1]$.
\end{definition}
The natural morphism of sheaves of complexes $\mathbb{Z}[n+1]^{\infty}_D\to U(1)[n]_D^\infty$ is an acyclic fibration and so we have an induced acyclic fibration (in the global model structure) $\mathbf{B}^{n}(\mathbb{Z}\hookrightarrow\mathbb{R})_{\mathrm{conn}}\stackrel{\sim}{\twoheadrightarrow}
\mathbf{B}^{n}U(1)_{\mathrm{conn}}$. Therefore, 
we find a natural isomorphism 
\[
H^0(X,U(1)[n]^\infty_D)\simeq H^0(X,\mathbb{Z}[n+1]^{\infty}_D)
\]
and a commutative diagram
\[
\xymatrix{
H^0(X,U(1)[n]^\infty_D)\ar[d]^{\begin{turn}{270}$\scriptstyle{\simeq}$\end{turn}}\ar[r]&H^n(X,U(1))\ar[d]^{\begin{turn}{270}$\scriptstyle{\simeq}$\end{turn}}
\\
H^0(X,\mathbb{Z}[n+1]^{\infty}_D)\ar[r]& H^{n+1}(X,\mathbb{Z})\,.
}
\]
\begin{definition}
We denote the cohomology group $H^0(X,\mathbb{Z}[n]^{\infty}_D)$ by the symbol $\hat{H}^n(X,\mathbb{Z})$, and call it the $n$-th \emph{differential cohomology} group of $X$ (with integer coefficients). The natural morphism $\hat{H}^n(X,\mathbb{Z})\to H^n(X,\mathbb{Z})$ will be called the \emph{differential refinement} of ordinary cohomology.
\end{definition}
\begin{remark}The reader experienced 
with gerbes and higher gerbes  will have recognized that
$H^n(X,U(1))\simeq H^{n+1}(X,\mathbb{Z})$ is the set of isomorphim classes of $U(1)$-$(n-1)$-gerbes on a manifold $X$, whereas $H^0(X,U(1)[n]^\infty_D)\simeq \hat{H}^{n+1}(X,\mathbb{Z})$ is the set of isomorphim classes of $U(1)$-$(n-1)$-gerbes with connection on $X$, and that the natural morphism $\hat{H}^{n+1}(X,\mathbb{Z})\to H^{n+1}(X,\mathbb{Z})$ is `forgetting the connection', see, e.g., \cite{gajer}.
\end{remark}
The natural projection $U(1)[n]_D^\infty\to C^\infty(-,U(1)[n]$ is a fibration, so we have a natural fibration $\mathbf{B}^nU(1)_{\mathrm{conn}}\twoheadrightarrow \mathbf{B}^nU(1)$, and, as for the case of Lie groups, we have a natural factorization 
\[
\mathbf{B}^nU(1)_{\mathrm{conn}}\hookrightarrow \mathbf{B}^nU(1)_{\mathrm{diff}}\stackrel{\sim}{\twoheadrightarrow} \mathbf{B}^nU(1)
\]
 into a monomorphism followed by an acyclic fibration (in the global model structure).
 
 The smooth $\infty$-groupoid $ \mathbf{B}^nU(1)_{\mathrm{diff}}$ is best defined at the level of chain complexes, where we have the well known ``cone trick'' from homological algebra to get the desired factorization.  In the case at hand, it works as follows: let $\mathrm{cone}(\ker\pi\hookrightarrow U(1)[n]_D^\infty(U))$ be  the mapping cone of the inclusion of the kernel of $\pi:U(1)[n]_D^\infty\to C^\infty(-,U(1)[n]$ into $U(1)[n]_D^\infty$, i.e., the 
chain complex
\[
\xymatrix@R=1pt{
       C^\infty(-,U(1)) 
        \ar[r]^{d \mathrm{log}} & 
       \Omega^1(-)
       \ar[r]^d & \Omega^2(-) 
       \ar[r]^d & \cdots
        \ar[r]^d
        &
        \Omega^n(-)
       \\
       \oplus & \oplus & \oplus&\oplus
       \\
       \Omega^1(-) \ar[uur]|{\mathrm{Id}} \ar[r]_d
       & \Omega^2(-)
         \ar@{--}[uur]
         \ar[r]
         & \cdots
         \ar[r]
         &
         \Omega^n(-)
         \ar[uur]|{\mathrm{Id}}
         \ar[r]
         &0
    }
\]
 Then $U(1)[n]_D^\infty$ naturally injects into $\mathrm{cone}(\ker\pi\hookrightarrow U(1)[n]_D^\infty)$, and $\pi$ induces a morphism of complexes $\pi:\mathrm{cone}(\ker\pi\hookrightarrow U(1)[n]_D^\infty)\to C^\infty(-,U(1)[n]$ which is an acyclic fibration; the composition
\[
U(1)[n]_D^\infty\hookrightarrow \mathrm{cone}(\ker\pi\hookrightarrow U(1)[n]_D^\infty)\xrightarrow{\pi}C^\infty(-,U(1)[n]
\]
is the sought for  factorization.
\begin{definition}
Define the simplicial presheaf
\[
  \mathbf{B}^n U(1)_{\mathrm{diff}}
  = 
  \mathrm{DK}
  \left(\mathrm{cone}(\ker\pi\hookrightarrow U(1)[n]_D^\infty)
  \right)
\]
  to be the image under the Dold-Kan equivalence of the chain complex of sheaves of abelian groups $\mathrm{cone}(\ker\pi\hookrightarrow U(1)[n]_D^\infty)$.
\end{definition}  
\par
The last smooth $\infty$-groupoid we introduce in this section is the natural ambient for curvature forms to live in. As above, we work at the level of sheaves of chain complexes first. So, let $\flat\mathbb{R}[n]_{\mathrm{dR}}^\infty$ be the truncated de Rham complex
\[
\flat\mathbb{R}[n+1]_{\mathrm{dR}}^\infty:=\left(\Omega^1(-) \stackrel{d}{\to} \Omega^2(-) \stackrel{d}{\to} \cdots \stackrel{d}{\to} \Omega^{n+1}_{\mathrm{cl}}(-)\right)
\]
seen as a chain complex concentrated in degrees $[0,n]$. 

There is a natural morphism of complexes of sheaves, which we call the \emph{curvature map},
\[
\mathrm{curv}:\mathrm{cone}(\ker\pi\hookrightarrow U(1)[n]_D^\infty)\to \flat\mathbb{R}[n+1]_{\mathrm{dR}}^\infty 
\]
given by the projection $\mathrm{cone}(\ker\pi\hookrightarrow U(1)[n]_D^\infty)\to\ker\pi[1]$ in degrees $[1,n]$ and given by the de Rham differential $d:\Omega^n(-)\to \Omega^{n+1}_{\mathrm{cl}}(-)$ in degree zero.
Note that the preimage of $(0 \to 0 \to \cdots \to \Omega^{n+1}_{\mathrm{cl}}(-))$ via $\mathrm{curv}$ is precisely the complex $U(1)[n]_D^\infty$, and that for $n=1$ the induced morphism
\[
\mathrm{curv}:U(1)[1]_D^\infty\to \Omega^2_{\mathrm{cl}}(-)
\]
is the map sending a connection on a principal $U(1)$-bundle to its curvature 2-form.

\begin{definition} The smooth $\infty$-groupoid $\mathbf{\flat}_{\mathrm{dR}} \mathbf{B}^{n+1}\mathbb{R}$ is
  $$
    \mathbf{\flat}_{\mathrm{dR}} \mathbf{B}^{n+1}\mathbb{R} = 
    \mathrm{DK}\left(
      \Omega^1(-) \stackrel{d}{\to} \Omega^2(-) \stackrel{d}{\to} \cdots \stackrel{d}{\to} \Omega^{n+1}_{\mathrm{cl}}(-)
    \right),
  $$
   the image under DK of the truncated de Rham complex.
\end{definition}
The above discussion can be summarized as
\begin{proposition}\label{prop.curv}
In $[\mathrm{CartSp^{\mathrm{op}}}, \mathrm{sSet}]_{\mathrm{proj},\mathrm{loc}}$ we have a natural commutative diagram
\[
  \xymatrix{
    \mathbf{B}^n U(1)_{\mathrm{conn}}
    \ar[d]
    \ar[r]
    & \Omega^{n+1}_{\mathrm{cl}}(-)
    \ar[d]
    \\
    \mathbf{B}^n U(1)_{\mathrm{diff}}
    \ar[r]^{\mathrm{curv}}
    \ar@{->>}[d]^\wr
    &
    \mathbf{\flat}_{\mathrm{dR}} \mathbf{B}^{n+1}\mathbb{R} 
    \\
    \mathbf{B}^n U(1)
  }
    \,.
\]
whose upper square is a pullback and whose lower part presents a morphism of smooth $\infty$-groupoids from $\mathbf{B}^n U(1)$ to 
$\mathbf{\flat}_{\mathrm{dR}}\mathbf{B}^{n+1}\mathbb{R}$. We call this morphism the 
\emph{curvature characteristic map}.
\end{proposition}

\begin{remark}
One also has a natural $(\mathbb{Z}\hookrightarrow \mathbb{R})$ version of $\mathbf{B}^n U(1)_{\mathrm{diff}}$, i.e., we have a smooth $\infty$-groupoid $\mathbf{B}^{n}(\mathbb{Z}\hookrightarrow\mathbb{R})_\mathrm{diff}$
with a natural morphism
\[
\mathbf{B}^{n}(\mathbb{Z}\hookrightarrow\mathbb{R})_{\mathrm{diff}}\stackrel{\sim}{\twoheadrightarrow}
 \mathbf{B}^nU(1)_{\mathrm{diff}}
\]   
which is an acyclic fibration in the global model structure.
\end{remark}

\section{Differential $\infty$-Lie integration}
\label{ooLieIntegration}

As the notion of $L_\infty$-algebra 
 generalizes
that of Lie algebra so that of  \emph{Lie $\infty$-group} generalizes that of  Lie group. 
We describe a way to \emph{integrate}
an $L_\infty$-algebra $\mathfrak{g}$ to the smooth delooping $\mathbf{B}G$ of the corresponding
Lie $\infty$-group by a slight variant of the construction of \cite{henriques}.
(Recall from the introduction that we use ``Lie'' to indicate generalized 
smooth structure which may or may not be represented by smooth manifolds). 
Then we generalize this to a \emph{differential} integration: an integration of  an $L_\infty$-algebroid
$\mathfrak{g}$
to smooth $\infty$-groupoids $\mathbf{B}G_{\mathrm{diff}}$ and $\mathbf{B}G_{\mathrm{conn}}$. 
Cocycles with coefficients in $\mathbf{B}G$ give $G$-principal $\infty$-bundles; those with coefficients
in $\mathbf{B}G_{\mathrm{diff}}$ support the $\infty$-Chern-Weil homomorphism,  those with coefficients
in $\mathbf{B}G_{\mathrm{conn}}$ give $G$-principal $\infty$-bundles with connection.

\subsection{Lie $\infty$-Algebroids: cocycles, invariant polynomials and CS-elements}
\label{ooLieAlgebroid}
\label{section.Lie_infinity-algebroids}

We summarize the main definitions and properties of $L_\infty$-algebroids from \cite{ls, lada-markl}, and their
cocycles, invariant polynomials and Chern-Simons elements from \cite{SSSI, SSSIII}.

\begin{definition}
Let $R$ be a commutative $\mathbb{R}$-algebra, and let $\mathfrak{g}$ be a chain complex of finitely generated (in each degree) $R$-modules, concentrated in non-negative degree. Then a (reduced)  
$L_\infty$-algebroid (or Lie $\infty$-algebroid) structure on $\mathfrak{g}$ is the datum of a degree 1 
$\mathbb{R}$-derivation $d_{\mathrm{CE}(\mathfrak{g})}$ on the exterior algebra
\[
  \wedge_R^\bullet\,\mathfrak{g}^*:=\mathrm{Sym}_R^\bullet(\mathfrak{g}^*[-1])
\]
(the free graded commutative algebra on the shifted dual of $\mathfrak{g}$), which is a differential (i.e., squares to zero) compatible with the differential of $\mathfrak{g}$.

A chain complex $\mathfrak{g}$ endowed with an $L_\infty$-algebroid structure will be called a \emph{$L_\infty$-algebroid}. The differential graded commutative algebra
\[
\mathrm{CE}(\mathfrak{g}):=(\wedge_R^\bullet\,\mathfrak{g}^*,d_{\mathrm{CE}(\mathfrak{g}}))
\]
will be called the \emph{Chevalley-Eilenberg algebra} of the 
$L_\infty$-algebroid $\mathfrak{g}$. A morphism of $L_\infty$-algebroids $\mathfrak{g}_1\to\mathfrak{g}_2$ is defined to be a dgca morphism $\mathrm{CE}(\mathfrak{g}_2)\to \mathrm{CE}(\mathfrak{g}_1)$
\end{definition}
Since all  $L_\infty$-algebroids which will be met in this  paper will be reduced, we will just say Lie $\infty$-algebroid to mean reduced Lie $\infty$-algebroid in what follows.
\begin{remark}\label{Loormk}
$\,$
\begin{itemize}

\item Lie $\infty$-algebroids could be more intrinsically defined as follows: the category  $L_\infty \mathrm{Algd} \subset dg\mathrm{Alg}^{op}$ of \emph{ $L_\infty$-algebroids} is the full subcategory of
 the opposite of that of  differential graded commutative $\mathbb{R}$-algebras
   on those dg-algebras whose underlying graded-commutative algebra is free on a finitely generated graded module concentrated in positive
 degree over the commutative algebra in degree 0.
\item
The dual $\mathfrak{g}^*=\mathrm{Hom}^{-\bullet}_R(\mathfrak{g},R)$ is a \emph{cochain} complex concentrated in nonnegative degrees. In particular the shift to the right functor $[-1]$ changes it into a cochain complex concentrated in strictly positive degrees.
\item
The restriction to finite generation is an artifact of dualizing $\mathfrak{g}$ rather than working with graded alternating multilinear functions on $\mathfrak{g}$ as the masters (Chevalley-Eilenberg-Koszul) did in the original ungraded case.  In particular, the direct generalization of their approach consists in working with the cofree connected cocommutative coalgebra cogenerated by 
$\mathfrak{g}[1]$, see, e.g., \cite{jds:intrinsic}. At least for $L_\infty$-algebras,  there are alternate definitions and conventions as to bounds on the grading, signs, etc. cf. \cite{ls, lada-markl} among others.

   \item Given an $L_\infty$-algebroid $\mathfrak{g}$, the degree 0 part $\mathrm{CE}(\mathfrak{g})_0$ of $\mathrm{CE}(\mathfrak{g})$ is is a commutative $\mathbb{R}$-algebra  which we think
of as the formal dual to the space of objects over which the  Lie $\infty$-algebroid is
defined. If $\mathrm{CE}(\mathfrak{g})_0 = \mathbb{R}$ equals the ground field, we say we have an $\infty$  algebroid
over the point, or equivalently that we have an \emph{$L_\infty$-algebra}.

  \item The underlying algebra in degree 0 can be generalized to an algebra over 
    some Lawvere theory. In particular in a proper setup of higher differential geometry,
    we would demand $\mathrm{CE}(\mathfrak{g})_0$ to be equipped with the structure of a 
    $C^\infty$-ring.    
\end{itemize}
\end{remark}

\begin{example}
\,
\begin{itemize}
  \item
    For $\mathfrak{g}$ an ordinary (finite-dimensional) Lie algebra, $\mathrm{CE}(\mathfrak{g})$
    is the ordinary Chevalley-Eilenberg algebra with coefficients in $\mathbb{R}$. The differential
    is given by the dual of the Lie bracket, 
    $$
      d_{\mathrm{CE}(\mathfrak{g})} = [-,-]^*
    $$  
    extended uniquely as a graded derivation.
    
    \item
    For a dg-Lie algebra $\mathfrak{g} = (\mathfrak{g}_\bullet, \partial)$, the differential is
    $$
      d_{\mathrm{CE}(\mathfrak{g})} = [-,-]^* + \partial^*
      \,.
    $$  
 \item   
     In the general case, 
      the total differential is further determined by (and is equivalent to) a sequence of higher 
      multilinear brackets \cite{ls}.

  \item
    For $n \in \mathbb{N},$ the $L_\infty$-algebra $b^{n-1} \mathbb{R}$ is defined
    in terms of $\mathrm{CE}(b^{n-1}\mathbb{R})$ which is the dgc-algebra on a single generator in degree $n$
    with vanishing differential.
  \item
    For $X$ a smooth manifold, its \emph{tangent Lie algebroid} is defined to have
    $\mathrm{CE}(T X) = (\Omega^\bullet(X), d_{dR})$ the de Rham algebra of $X$. Notice that
    $\Omega^{\bullet}(X) = \wedge^\bullet_{C^\infty(X)} \Gamma(T^* X)$.
    
    We shall extensively use the tangent Lie algebroid $T (U \times \Delta^k)$ where $U \in \mathrm{CartSp}$
    and $\Delta^k$ is the standard $k$-simplex.
\end{itemize}
\end{example} 

\begin{definition}
  For $\mathfrak{g}$ a  Lie $\infty$-algebroid and $n \in \mathbb{N}$, a 
  \emph{cocycle} in degree $n$ on $\mathfrak{g}$ is, equivalently
  \begin{itemize}
    \item an element $\mu \in \mathrm{CE}(\mathfrak{g})$ in degree $n$, such that 
     $d_{\mathrm{CE}(\mathfrak{g})}\; \mu = 0$;
    \item a morphism of dg-algebras $\mu: \mathrm{CE}(b^{n-1}\mathbb{R})\to \mathrm{CE}(\mathfrak{g})$;
    
    \item a morphism of Lie $\infty$-algebroids $\mu : \mathfrak{g} \to b^{n-1}\mathbb{R}$.
  \end{itemize}

\end{definition}
\begin{example}
\begin{itemize}
  \item
    For $\mathfrak{g}$ an ordinary Lie algebra, a cocycle in the above sense is the same as
    a Lie algebra cocycle in the ordinary sense (with values in the trivial module).
  \item
    For $X$ a smooth manifold, a cocycle in degree $n$ on the tangent Lie algebroid $T X$
    is precisely a closed $n$-form on $X$.
\end{itemize}
\end{example}
For our purposes, a particularly important Chevalley-Eilenberg algebra is the Weil algebra.

\begin{definition}
\label{weilalg}
 The \emph{Weil algebra} of an $L_\infty$-algebra $\mathfrak{g}$ is
 the dg-algebra
 $$ 
   \mathrm{W}(\mathfrak{g}) := (\mathrm{Sym}^\bullet (\mathfrak{g}^*[-1] \oplus \mathfrak{g}^*[-2]), d_{W(\mathfrak{g})})
   \,,
 $$ 
 where the differential  on the copy $\mathfrak{g}^*[-1]$ is the sum
 $$
   d_{W(\mathfrak{g})}|_{\mathfrak{g}^*} = d_{CE(\mathfrak{g})} + \sigma
   \,,
 $$
 with $\sigma : \mathfrak{g}^* \to \mathfrak{g}^*[-1]$ is the grade-shifting isomorphism, i.e. it is the identity of $\mathfrak{g}^*$ seen as a degree 1 map $\mathfrak{g}^*[-1]\to \mathfrak{g}^*[-2]$, extended
 as a graded derivation, and where
 $$
   d_{W(\mathfrak{g})} \circ \sigma = - \sigma \circ d_{W(\mathfrak{g})}
   \,.
 $$
\end{definition}
\begin{proposition}\label{weil-is-free}
  The Weil algebra is a representative of the free differential graded commutative algebra 
  on the graded vector space $\mathfrak{g}^*[-1]$ in that there exist a natural isomorphism
  \[
  \mathrm{Hom}_{\rm{dgca}}(W(\mathfrak{g}),\Omega^\bullet)\simeq 
   \mathrm{Hom}_{\rm{gr-vect}}(\mathfrak{g}^*[-1],\Omega^\bullet),
  \]
  for $\Omega^\bullet$ an arbitrary dgca. Moreover, the Weil algebra is precisely that algebra
  with this property for which the projection morphism 
  $i^*  : \mathfrak{g}^*[-1]\oplus \mathfrak{g}^*[-2] \to \mathfrak{g}^*[-1]$ of graded vector spaces
  extends to a dg-algebra homomorphism 
  $$
    i^* : \mathrm{W}(\mathfrak{g}) \to \mathrm{CE}(\mathfrak{g})
    \,.
  $$
\end{proposition}
Notice that the free dgca on a graded vector space is defined only up to isomorphism. 
The condition on $i^*$ is what picks the Weil algebra among all free dg-algebras. A proof of the above proposition can be found, e.g., in \cite{SSSI}.

Equivalently, one can state the freeness of the Weil algebra by saying that the dgca-morphisms 
$A:\mathrm{W}(\mathfrak{g})\to \Omega^\bullet$ are in natural bijection with the degree 1 elements in the 
graded vector space $\Omega^\bullet\otimes\mathfrak{g}$. 
\begin{example}\label{example.weil-algebras}

\begin{itemize}
  \item For $\mathfrak{g}$ an ordinary Lie algebra, $\mathrm{W}(\mathfrak{g})$ is the ordinary Weil algebra 
   \cite{cartan:g-alg-sem-I}. In that paper, H. Cartan defines a $\mathfrak{g}$-algebra  as an analog of the dg-algebra $\Omega^\bullet(P)$ of differential forms on a principal bundle, i.e. as a dg-algebra equipped with operations $i_\xi$ and $L_\xi$ for all $\xi\in  \mathfrak{g}$ satisfying the usual relations, including
   $$L_\xi= di_\xi + i_\xi d.$$
 Next, Cartan introduces the Weil algebra $\mathrm{W}(\mathfrak{g})$ as the universal $\mathfrak{g}$-algebra and identifies a $\mathfrak{g}$-connection $A$ on a principal bundle  $P$ as a morphism of $\mathfrak{g}$-algebras
\[
A:\mathrm{W}(\mathfrak{g})\to \Omega^\bullet(P).
\]
This can in turn be seen as a dgca morphism satisfying the 
Cartan-Ehresmann conditions, and it is this latter point of view that we generalize to an arbitrary $L_\infty$-algebra.
\item
     The dg-algebra $\mathrm{W}(b^{n-1}\mathbb{R})$ of $b^{n-1}\mathbb{R}$ 
     is the free dg-algebra on a single generator
     in degree $n$. As a graded algebra, it has a generator $b$ in degree $n$ and a generator
     $c$ in degree $(n+1)$ and the differential acts as $d_\mathrm{W} : b \mapsto c$. Note that, since $d_{\mathrm{CE}}b=0$, this is equivalent to $c=\sigma b$.
          
\end{itemize}
\end{example}
\begin{remark}\label{remark.curvature-forms}
Since the Weil algebra is itself a dg-algebra whose underlying graded algebra is
a graded symmetric algebra, it is itself the CE-algebra of an $L_\infty$-algebra. The $L_\infty$-algebra
thus defined we denote $\mathrm{inn}(\mathfrak{g})$:
$$
  \mathrm{CE}(\mathrm{inn}(\mathfrak{g})) = \mathrm{W}(\mathfrak{g})
  \,.
$$
Note that the underlying graded vector space of $\mathrm{inn}(\mathfrak{g})$ is $\mathfrak{g}\oplus\mathfrak{g}[1]$.  Looking at $\mathrm{W}(\mathfrak{g})$ as the Chevalley-Eilenberg algebra of $\mathrm{inn}(\mathfrak{g})$ we therefore obtain the following descrition of morphisms out of $\mathrm{W}(\mathfrak{g})$: for any dgca $\Omega^\bullet$, a dgca morphism $\mathrm{W}(\mathfrak{g})\to \Omega^\bullet$ is the datum of a pair $(A,F_A)$, where $A$ and $F_A$ are a degree 1 and a degree 2 element in $\Omega^\bullet\otimes \mathfrak{g}$, respectively, such that $(A,F_A)$ satisfies the Maurer-Cartan equation in the $L_\infty$-algebra $\Omega^\bullet\otimes \mathrm{inn}(\mathfrak{g})$. The Maurer-Cartan equation actually completely determines $F_A$ in terms of $A$; this is an instance of the freeness property of the Weil algebra stated in Proposition \ref{weil-is-free}.
\end{remark}

\begin{definition}
  For $X$ a smooth manifold, a \emph{$\mathfrak{g}$-valued connection form} on $X$ is a morpism of Lie $\infty$-algebroids
  $A : T X \to \mathrm{inn}(\mathfrak{g})$, hence a morphism of dg-algebras
  $$
   A:\mathrm{W}(\mathfrak{g})\to  \Omega^\bullet(X)
    \,.
  $$  
 \end{definition}
 \begin{remark}
A $\mathfrak{g}$-valued connection form on $X$ can be equivalently seen as an element $A$ in the set $\Omega^1(X,\mathfrak{g})$ of degree 1 elements in $\Omega^\bullet(X)\otimes\mathfrak{g}$, or as a pair $(A,F_A)$, where $A\in\Omega^1(U,\mathfrak{g})$,  $F_A\in\Omega^2(U,\mathfrak{g})$, and $A$ and $F_A$ are related by the Maurer-Cartan equation in $\Omega^\bullet(X,\mathrm{inn}(\mathfrak{g}))$. The element $F_A$ is called the \emph{curvature form} of $A$.
\end{remark}
\begin{example}
If $\mathfrak{g}$ is an ordinary Lie algebra, then a $\mathfrak{g}$-valued connection form $A$ on $X$ is a 1-form on $X$ with coefficients in $\mathfrak{g}$, i.e. it is naturally a connection 1-form on a trivial principal $G$-bundle on $X$. The element $F_A$ in $\Omega^2(X,\mathfrak{g})$ is then given by equation
\[
F_A=dA+\frac{1}{2}[A,A],
\]
so it is precislely the usual curvature form of $A$.
\end{example}
\bigskip

The last ingredient we need to generalize 
from Lie algebras to $L_\infty$-algebroids is the algebra 
$\mathrm{inv}(\mathfrak{g})$ of invariant polynomials. 
\begin{definition}
An \emph{invariant polynomial} on $\mathfrak{g}$ is a  $d_{\mathrm{W}(\mathfrak{g})}$-closed  element 
  $\langle - \rangle$ in $\mathrm{Sym}^\bullet (\mathfrak{g}^*[-2]) \subset \mathrm{W}(\mathfrak{g})$.
 \end{definition} 
To see how this definition encodes the classical definition of invariant polynomials on a Lie algebra, notice that invariant polynomials are elements of $\mathrm{Sym}^\bullet (\mathfrak{g}^*[-2])$ that are both horizontal and $\mathrm{ad}$-invariant (``\emph{basic} forms''). Namely, for any $v \in \mathfrak{g}$ we have, for an invariant polynomial $\langle - \rangle$, the identities
$$
    \iota_v \langle - \rangle = 0 \;\;\; \mbox{(horizontality)}
  $$
  and
  $$
     \mathcal{L}_v \langle - \rangle  = 0
     \;\;\;\; \mbox{($\mathrm{ad}$-invariance)}
    \,,
    $$
where    $\iota_v : \mathrm{W}(\mathfrak{g}) \to \mathrm{W}(\mathfrak{g})$
 the contraction derivation defined by $v$ and  $\mathcal{L}_v:= [d_{\mathrm{W}(\mathfrak{g})}, \iota_v]$ is the corresponding Lie derivative.  
\bigskip

We want to identify two indecomposable invariant polynomials which differ
by a ``horizontal shift''. A systematic way of doing this is to introduce the
following equivalence relation on the dgca of all invariant polynomials: we say that two invariant polynomials $\langle -\rangle_1, \langle -\rangle_2$ are \emph{horizontally equivalent}
  if there exists $\omega$ in $\mathrm{ker}(\mathrm{W}(\mathfrak{g}) \to \mathrm{CE}(\mathfrak{g}))$ such that
  $$
    \langle - \rangle_1 = \langle -\rangle_2 + d_{\mathrm{W}} \omega
    \,.
  $$
 Write $\mathrm{inv}(\mathfrak{g})_V$ for the quotient graded vector space of horizontal equivalence classes of invariant polynomials. 
\begin{definition}
  \label{InvariantPolynomials}
The dgca $\mathrm{inv}(\mathfrak{g})$ is defined as the free polynomial algebra on the graded vector space $\mathrm{inv}_V(\mathfrak{g})$, endowed with the trivial differential. 
\end{definition}
\begin{remark}
 A choice of a linear section to the projection
 \[
 \{\text{invariant polynomials}\}\to \mathrm{inv}(\mathfrak{g})_V
 \]
 gives a
morphism of graded vector spaces $\mathrm{inv}(\mathfrak{g})_V \to \mathrm{W}(\mathfrak{g})$,
canonical up to horizontal homotopy, that sends each equivalence class to a representative. This linear morphism uniquely extends to a dg-algebra homomorphism
$$
  \mathrm{inv}(\mathfrak{g}) \to \mathrm{W}(\mathfrak{g})
  \,.
$$  
 \end{remark}
 \begin{remark}
 The algebra $\mathrm{inv}(\mathfrak{g})$ is at first sight a quite abstract construction which is apparently unrelated to an equivalence relation on indecomposable invariant polynomials. A closer look shows that it is actually not so. Namely, only indecomposable invariant polynomials can be representatives for the nonzero equivalence classes. Indeed, if $\langle- \rangle_1$
  and $\langle -\rangle_2$ are two nontrivial invariant polynomials, then since the cohomology of 
  $\mathrm{W}(\mathfrak{g})$ is trivial in positive degree, there is $\mathrm{cs}_1$ in $\mathrm{W}(\mathfrak{g})$ (not necessarily in $\mathrm{ker}(\mathrm{W}(\mathfrak{g}) \to \mathrm{CE}(\mathfrak{g}))$ ) such that 
  $d_{\mathrm{W}} \mathrm{cs}_1 = \langle - \rangle_1$, but then $\mathrm{cs}_1 \wedge \langle - \rangle_2$ is a horizontal trivialization of $\langle - \rangle_1 \wedge \langle - \rangle_2$.
One therefore obtains a very concrete description of the algebra $\mathrm{inv}(\mathfrak{g})$ as follows: one picks a representative indecomposable invariant polynomial for each horizontal equivalence class and considers the subalgebra of $\mathrm{W}(\mathfrak{g})$ generated by these representatives. The morphism $\mathrm{inv}(\mathfrak{g})\to \mathrm{W}(\mathfrak{g})$ is then realized as the inclusion of this subalgebra into the Weil algebra. Different choices of representative generators lead to distinct but equivalent subalgebras: each one is isomorphic to the others via an horizontal shift in the generators.
\end{remark}
\begin{remark}
  For $\mathfrak{g}$ an ordinary reductive Lie algebra, definition \ref{InvariantPolynomials} 
  reproduces the traditional
  definition of the algebra of $\mathrm{ad}_{\mathfrak{g}}$-invariant polynomials. Indeed, for a Lie algebra $\mathfrak{g}$, the condition $d_{\mathrm{W}(\mathfrak{g})} \langle - \rangle = 0$ 
  is precisely the usual $\mathrm{ad}_{\mathfrak{g}}$-invariance of an element $\langle - \rangle$  in $\mathrm{Sym}^\bullet \mathfrak{g}^*[-2]$. Morover, the horizonal equivalence on indecomposables is trivial in this case and it 
  is a classical fact (for instance theorem I on page 242 in volume III of \cite{GHV}) that the graded algebra of $\mathrm{ad}_{\mathfrak{g}}$-invariant polynomials
  is indeed free on the space of indecomposables.
\end{remark}

\begin{definition}
For any dgca morphism $A:\mathrm{W}(\mathfrak{g})\to \Omega^\bullet$, the composite morphism
$\mathrm{inv}(\mathfrak{g})\to \mathrm{W}(\mathfrak{g})\to \Omega^\bullet$ is the evaluation of invariant polynomials on the element $F_A$. In particular, if $X$ is a smooth manifold and $A$ is a $\mathfrak{g}$-valued connection form on $X$, then the image of $F_A:\mathrm{inv}(\mathfrak{g})\to \Omega^\bullet(X)$ is a collection of differential forms on $X$, to be called the \emph{curvature characteristic forms} of $A$. 
\end{definition}
\par
\begin{example}
For the $L_\infty$-algbera $b^{n-1}\mathbb{R}$, in the notations of Example \ref{example.weil-algebras}, one has $\mathrm{inv}(b^{n-1}\mathbb{R})=\mathbb{R}[c]$.
\end{example}
\begin{definition}
  We say an invariant polynomial $\langle - \rangle$ on $\mathfrak{g}$ is \emph{in transgression} 
  with a cocycle $\mu$
  if there exists an element $\mathrm{cs} \in \mathrm{W}(\mathfrak{g})$ such that
  \begin{enumerate}
    \item
      $i^* \mathrm{cs} = \mu$;
    \item
      $d_{\mathrm{W}(\mathfrak{g})}\mathrm{cs}   = \langle - \rangle$.
  \end{enumerate}
  We call $\mathrm{cs}$ a \emph{Chern-Simons element} for $\mu$ and $\langle - \rangle$.
\end{definition}
For ordinary Lie algebras this reduces to the classical notion, for instance 6.13 in vol. III of \cite{GHV}.
\begin{remark}
If we think of $\mathrm{inv}(\mathfrak{g}) \subset \mathrm{ker}\, i^*$ as a subcomplex of the kernel
of $i^*$, then this transgression exhibits the connecting homomorphism 
$H^{n-1}(\mathrm{CE}(\mathfrak{g})) \to H^n(\mathrm{ker}\,i^*)$ of the long sequence in cohomology
induced from the short exact sequence 
$\mathrm{ker}\,i^* \to \mathrm{W}(\mathfrak{g}) \stackrel{i^*}{\to} \mathrm{CE}(\mathfrak{g})$.

If we think of 
\begin{itemize}
  \item $\mathrm{W}(\mathfrak{g})$ as differential forms on the total space of a universal principal bundle,
  \item $\mathrm{CE}(\mathfrak{g})$ as differential forms on the fiber,
  \item $\mathrm{inv}(\mathfrak{g})$ as forms on the base,
\end{itemize}
then the above notion of transgression is precisely the classical one of transgression of forms
in the setting of fiber bundles (for instance section 9 of \cite{borel}).

\end{remark}
\begin{example} \label{ExamplesCSElements}
\hspace{1pt}
\begin{itemize}
 \item
   For $\mathfrak{g}$ a semisimple Lie algebra with $\langle -,-\rangle$ the Killing form
   invariant polynomial, the corresponding cocycle in transgression is 
   $\mu_3 = \frac{1}{2} \langle -,[-,-]\rangle$. The Chern-Simons element witnessing this transgression
   is $\mathrm{cs} = \langle \sigma(-),-\rangle  + \frac{1}{2} \langle -,[-,-]\rangle$.
 \item
 For the Weil algebra $\mathrm{W}(b^{n-1}\mathbb{R})$ of Example \ref{example.weil-algebras}, the element $b$ (as element of the Weil algebra) is a Chern-Simons element transgressing the cocycle $b$ (as element of the Chevalley-Eilenberg algebra $\mathrm{CE}(b^{n-1}\mathbb{R})$) to the invariant polynomial $c$.
  \item
    For $\mathfrak{g}$ a semisimple Lie algebra, $\mu_3 = \frac{1}{2}\langle -,[-,-]\rangle$ the canonical Lie algebra 3-cocycle in transgression with the Killing form, let $\mathfrak{g}_{\mu_3}$ be the corresponding \emph{string Lie 2-algebra} given by the next definition, \ref{LooExtension},
  and discussed below in \ref{section.string-2_group}. Its Weil algebra is given by
 $$
    d_{W} t^a = - \frac{1}{2} C^a{}_{b c} t^b \wedge t^c + r^a
 $$
 $$
   d_W b = h - \mu_3 
 $$
 and the corresponding Bianchi identities,
 with $\{t^a\}$ a dual basis for $\mathfrak{g}$ in degree 1, with $b$ a generator in degree 2 and $h$ its curvature generator in degree 3. 
 We see that every invariant polynomial of $\mathfrak{g}$ is also an invariant polynomial of $\mathfrak{g}_{\mu_3}$. But the Killing form $\langle -,-\rangle$ is now horizontally trivial: let $\mathrm{cs}_3$ be any Chern-Simons element for $\langle -,-\rangle$ in $\mathrm{W}(\mathfrak{g})$. This is not horizontal. But  the element
$$
  \tilde {\mathrm{cs}}_3 := \mathrm{cs}_3 - \mu_3 + h
$$
is in $\mathrm{ker}(W(\mathfrak{g}_{\mu_3}) \to \mathrm{CE}(\mathfrak{g}_{\mu_3}))$ and
$$
  d_W \tilde {\mathrm{cs}_3} = \langle -,- \rangle
  \,.
$$
  Therefore $\mathrm{inv}(\mathfrak{g}_{\mu_3})$ has the same generators as $\mathrm{inv}(\mathfrak{g})$ except the Killing form, which  is discarded.
\end{itemize}
\end{example}
\begin{definition}
  \label{LooExtension}
  For $\mu : \mathfrak{g} \to b^{n-1}\mathbb{R}$ a cocycle in degree $n \geq 1$, 
  the \emph{extension} that it classifies
  is the $L_\infty$-algebra given by the pullback
  $$
    \xymatrix{
      \mathfrak{g}_\mu \ar[r] \ar[d] & \mathrm{inn}(b^{n-2}\mathbb{R}) \ar[d]
      \\
      \mathfrak{g} \ar[r]^\mu & b^{n-1}\mathbb{R}
    }
    \,.
  $$
\end{definition}
\begin{remark} 
Dually, the $L_\infty$-algebra $\mathfrak{g}_\mu$ is the pushout
  $$
    \xymatrix{
      \mathrm{CE}(\mathfrak{g}_\mu) \ar@{<-}[r] \ar@{<-}[d] 
      & 
       \mathrm{W}(b^{n-2}\mathbb{R}) \ar@{<-}[d]
      \\
      \mathrm{CE}(\mathfrak{g}) \ar@{<-}[r]^\mu & \mathrm{CE}(b^{n-1}\mathbb{R})
    }
  $$
  in the category $\mathrm{dgcAlg}$.
  This means that $\mathrm{CE}(\mathfrak{g}_\mu)$ is obtained from $\mathrm{CE}(\mathfrak{g})$
  by adding one more generator $b$ in degree $(n-1)$ and setting 
  $$
    d_{\mathrm{CE}(\mathfrak{g}_\mu)} : b \mapsto -\mu
    \,.
  $$
These are standard constructions on dgc-algebras familiar from rational homotopy theory,
realizing $\mathrm{CE}(\mathfrak{g}) \to \mathrm{CE}(\mathfrak{g}_\mu)$ as a relative Sullivan algebra. Yet, it is still worthwhile to make the $\infty$-Lie theoretic meaning in terms of $L_\infty$-algebra extensions manifest: we may think of $\mathfrak{g}_\mu$ as the homotopy fiber of $\mu$ or equivalently as the extension of $\mathfrak{g}$ classified by $\mu$.  In Section \ref{section.infinity_chern_weil_homomorphism} we discuss how these $L_\infty$-algebra extensions  are integrated to extensions of smooth $\infty$-groups; the homotopy fiber point of view will be emphasized in Section \ref{HomotopyFibers}.
\end{remark}
\begin{example}
  For $\mathfrak{g}$ a semisimple Lie algebra and $\mu = \frac{1}{2}\langle -,[-,-]\rangle$ the
  cocycle in transgression with the Killing form, the corresponding extension
  is the \emph{string Lie 2-algebra} $\mathfrak{g}_\mu$ discussed in \ref{section.string-2_group}, \cite{bcss, henriques}.
\end{example}
We may summarize the situation as follows:
for  $\mu$ a degree $n$ cocycle which is in transgression with  an invariant polynomial $\langle - \rangle$ via a Chern-Simons
  element $\mathrm{cs}$, the corresponding morphisms of dg-algebras fit into a commutative diagram
  $$
    \xymatrix{
      \mathrm{CE}(\mathfrak{g}) 
      \ar@{<-}[r]^{\mu}
      \ar@{<-}[d] 
       & \mathrm{CE}(b^{n-1}) \ar@{<-}[d]
      \\
      \mathrm{W}(\mathfrak{g})
      \ar@{<-}[r]^{\mathrm{cs}}
      \ar@{<-}[d]
      &
      \mathrm{W}(b^{n-1}\mathbb{R})
      \ar@{<-}[d]
      \\
      \mathrm{inv}(\mathfrak{g})
      \ar@{<-}[r]^{\langle- \rangle}
      &
      \mathrm{inv}(b^{n-1}\mathbb{R})
    }
  $$
In section \ref{section.infinity-stack_of_principal_bundles_with_connection} we will see that 
under \emph{$\infty$-Lie integration} this diagram corresponds to a universal circle $n$-bundle connection
on $\mathbf{B}G$. The  composition of the diagrams defining the cells in 
$\exp_\Delta(\mathfrak{g})_{\mathrm{diff}}$ (section \ref{section.infinity-stack_of_principal_bundles_with_connection})
with this diagram models the $\infty$-Chern-Weil homomorphism for the characteristic class given by $\langle-\rangle$.

\subsection{Principal $\infty$-bundles}
\label{section.principal_infinity-bundles}

We describe the integration of Lie $\infty$-algebras $\mathfrak{g}$ to smooth $\infty$-groupoids
$\mathbf{B}G$ in the sense of section \ref{section.Lie_infinity-groupoids}. 

The basic idea is Sullivan's old construction \cite{sullivan:inf} in rational homotopy theory 
of a simplicial set from a dg-algebra.
It was essentially noticed by Getzler \cite{getzler}, following Hinich \cite{hinich}, that this construction
may be interpreted in $\infty$-Lie  theory as forming the smooth $\infty$-groupoid underlying the Lie
integration of an $L_\infty$ -algebra. Henriques \cite{henriques} refined the construction to land
in $\infty$-groupoids internal to Banach spaces. Here we observe that the construction has an evident
refinement to yield genuine smooth $\infty$-groupoids in the sense of section 
\ref{section.Lie_infinity-groupoids} (this refinement has independently also been considered 
by Roytenberg in \cite{Roytenberg}): the integrated smooth $\infty$-groupoid sends each Cartesian space
$U$ to a Kan complex which in degree $k$ is the set of smoothly $U$-parameterized families of
smooth flat $\mathfrak{g}$-valued differential forms on 
the standard $k$-simplex $\Delta^k \subset \mathbb{R}^k$ regarded as a smooth manifold (with boundary and corners).

To make this precise we need a suitable notion of smooth differential forms on the $k$-simplex.
Recall that an ordinary smooth form on $\Delta^k$ is a smooth form on an open neighbourhood of 
$\Delta^n$ in $\mathbb{R}^n$. This says that the derivatives are well behaved at the boundary.
The following technical definition imposes even more restrictive conditions on the behaviour at the
boundary.
\begin{definition}
For any point $p$ in $\Delta^k$, let $\Delta_p$ be the lowest dimensional subsimplex of $\Delta^k$ the point $p$ belongs to, and let $\pi_p$ the orthogonal projection on the affine subspace spanned by $\Delta_p$.  A smooth differential form $\omega$ on $\Delta^k$ is said to have \emph{sitting instants} 
along the boundary if for any point $p$ in $\Delta^k$ there is a neighborhood $V_p$ of $p$ such that $\omega=\pi_p^*(\omega\vert_{\Delta_p})$ on $V_p$.

For any $U \in \mathrm{CartSp},$ a smooth differential form $\omega$ on $U \times\Delta^k$ is said to have
sitting instants if for all points $u : * \to U$ the pullback along 
$(u, \mathrm{Id}) : \Delta^k \to U \times \Delta^k$ has sitting instants.

Smooth forms with sitting instants clearly form a sub-dg-algebra of all smooth forms. We shall write
$\Omega^\bullet_{\mathrm{si}}(U \times \Delta^k)$ to denote this sub-dg-algebra.
\end{definition}
\begin{remark}
The inclusion $\Omega_{\mathrm{si}}^\bullet(\Delta^k)\hookrightarrow
\Omega^\bullet(\Delta^k)$ is a quasi-isomorphism. Indeed, by using bump
functions with sitting instants one sees that the sheaf of differential
forms with sitting instants is fine, and it is imemdiate to show that the
stalkwise Poincar\'e lemma holds for this sheaf. Hence the usual
hypercohomology argument applies. We thank Tom Goodwillie for having
suggested a sheaf-theoretic proof of this result. \end{remark}

\begin{remark}
For a point $p$ in the interior of the simplex $\Delta^k$ the sitting instants condition is clearly empty; this justifies the name ``sitting instants along the boundary''. Also note that the dimension of the normal direction to the boundary depends on the dimension of the boundary stratum: 
there is one perpendicular direction to a codimension-1 face, and there are $k$ perpendicular directions to a vertex. 
\end{remark}
\begin{definition}
For a Cartesian space $U$, we denote by the symbol 
$$\Omega_{\mathrm{si}}^\bullet(U\times \Delta^k)_{\mathrm{vert}} \subset \Omega^\bullet(U \times \Delta^k)$$  
the sub-dg-algebra on forms that are \emph{vertical} with respect to the projection
$U \times \Delta^k \to U$.

Equivalently this is the completed tensor product, 
$$\Omega_{\mathrm{si}}^\bullet(U\times \Delta^k)_{\mathrm{vert}}=C^\infty(U;\mathbb{R})\hat{\otimes}\Omega^\bullet_{\mathrm{si}}(\Delta^k),$$ where
$C^\infty(U;\mathbb{R})$ is regarded as a dg-algebra concentrated in degree zero.
\end{definition}

\begin{example}
\begin{itemize}
  \item
    A 0-form (a smooth function) has sitting instants on $\Delta^1$ if in a neighbourhood of the endpoints it is constant.
    A smooth function $f : U \times \Delta^1 \to \mathbb{R}$ is in $\Omega^0_{\mathrm{si}}(U \times \Delta^1)_{\mathrm{vert}}$
    if for each $u \in U$ it is constant in a neighbourhood of the endpoints of $\Delta^1$.
  \item
    A 1-form has sitting instants on $\Delta^1$ if in a neighbourhood of the endpoints it vanishes.
  \item
    Let $X$ be a smooth manifold and $\omega \in \Omega^\bullet(X)$ be a smooth form on $X$.
    Let $\phi : \Delta^k \to X$ be a smooth map with sitting instants in the ordinary sense:
    for every $r$-face of $\Delta^k$ there is a neighbourhood such that $\phi$ is perpendicularly
    constant on that neighbourhood. Then the pullback form $\phi^* \omega$ is a form with sitting
    instants on $\Delta^k$. 
\end{itemize}
\end{example}
\begin{remark}
The point of the definition of sitting instants, clearly reminiscent of the use of normal cylindrical collars in cobordism theory, is that, when glueing compatible forms on simplices along faces, the resulting differential form is smooth.
\end{remark}
\begin{proposition}
  \label{SmoothGluing}
  Let $\Lambda^k_i \subset \Delta^k$ be the $i$th horn of $\Delta^k$, regarded naturally as a 
  closed subset of $\mathbb{R}^{k-1}$. If $\{\omega_j \in \Omega_{\mathrm{si}}^\bullet(\Delta^{k-1})\}$ is a collection
  of smooth forms with sitting instants on the $(k-1)$-simplices of $\Lambda^k_i$ that match on their
  coinciding faces, then there is a unique smooth form $\omega$ on $\Lambda^k_i$ that restricts to 
  $\omega_j$ on the $j$th face.
\end{proposition} 
\proof
  By the condition that forms with sitting instants are constant perpendicular to their value on a face in 
  a neighbourhood of any face it follows that if two agree on an adjacent face then all derivatives 
  at that face of the unique form that extends both exist in all directions. Hence that unique form
  extending both is smooth.
\endofproof
\begin{definition}
 \label{LieIntegrationOfLInfinityAlgebra}
For $\mathfrak{g}$ an $L_\infty$-algebra, the simplicial presheaf $\exp_{\Delta}(\mathfrak{g})$ on the site of Cartesian spaces is defined as
\[
  \exp_\Delta(\mathfrak{g})
    :
(U,[k])\mapsto {\rm Hom}_{\mathrm{dgAlg}}(\mathrm{CE}(\mathfrak{g}),\Omega_{\mathrm{si}}^\bullet(U\times \Delta^k)_{\mathrm{vert}}).
\]
\end{definition}
Note that the construction of $\exp_\Delta(\mathfrak{g})$ is functorial in $\mathfrak{g}$: a morphism of $L_\infty$-algebras $\mathfrak{g_1}\to \mathfrak{g_2}$, i.e., a dg-algebra morphism 
$\mathrm{CE}(\mathfrak{g_2})\to \mathrm{CE}(\mathfrak{g}_1)$, induces a morphism of simplicial presheaves $\exp_\Delta(\mathfrak{g}_1)\to \exp_\Delta(\mathfrak{g}_2)$.
\begin{remark}
A $k$-simplex in $\exp_\Delta(\mathfrak{g})(U)$ may be thought of as a smooth family of \emph{flat} 
$\mathfrak{g}$-valued forms on $\Delta^n$, parametrized by $U$.
We write $\exp_\Delta(\mathfrak{g})$ for this simplicial presheaf to indicate that it plays a role analogous to the formal exponentiation of a Lie algebra to a Lie group. 
\end{remark}
\begin{proposition} 
  \label{proposition.exp.is.infty.groupoid}
  The simplicial presheaf $\exp_\Delta(\mathfrak{g})$ is a smooth $\infty$-groupoid in that 
  it is fibrant in $[\mathrm{CartSp}^{\mathrm{op}}, \mathrm{sSet}]_{\mathrm{proj}}$: it 
  takes values in Kan complexes. We say that the smooth $\infty$-groupoid $\exp_\Delta(\mathfrak{g})$ \emph{integrates} the $L_\infty$-algebra $\mathfrak{g}$.
\end{proposition}
\proof 
  Since our forms have sitting instants, this follows in direct analogy to the standard
  proof that the singluar simplicial complex of any topological space is a Kan complex: 
 we may use the standard retracts of simplices onto their horns to pull back forms from horns to
 simplices. The retraction maps are smooth except where they cross faces, but since the forms have
 sitting instants there, their smooth pullback exists nevertheless. 
 
  Let $\pi:\Delta^k\to\Lambda^k_i$ be the standard retraction map of a $k$-simplex on 
  its $i$-th horn. Since $\pi$ is smooth away from the primages of the faces, the commutative diagram
\[
\xymatrix{
U\times \Lambda^k_i\ar[r]^{\mathrm{id}\times i}\ar[dr]_{\mathrm{id}}&U\times \Delta^k\ar[d]^{\mathrm{id}\times \pi}\\
&U\times \Lambda^k_i
}
\]
induces a commutative diagram of dgcas
\[
\xymatrix{
\mathrm{Hom}_{\mathrm{dgAlg}}(\mathrm{CE}(\mathfrak{g}),\Omega_{\mathrm{si}}^\bullet(U\times \Lambda^k_i)_{\mathrm{vert}})&&\mathrm{Hom}_{\mathrm{dgAlg}}(\mathrm{CE}(\mathfrak{g}),\Omega_{\mathrm{si}}^\bullet(U\times \Delta^k)_{\mathrm{vert}})\ar[ll]_{(\mathrm{id}\times i)^*\circ-}\\
&&\mathrm{Hom}_{\mathrm{dgAlg}}(\mathrm{CE}(\mathfrak{g}),\Omega_{\mathrm{si}}^\bullet(U\times \Lambda^n_i)_{\mathrm{vert}})\ar[u]_{(\mathrm{id}\times \pi)^*\circ-}\ar[ull]^{\mathrm{id}}
}
 \,,
\]
so that, in particular, the horn-filling map $\mathrm{Hom}_{\mathrm{dgAlg}}(\mathrm{CE}(\mathfrak{g}),\Omega_{\mathrm{si}}^\bullet(U\times \Lambda^k_i)_{\mathrm{vert}})\to \mathrm{Hom}_{\mathrm{dgAlg}}(\mathrm{CE}(\mathfrak{g}),\Omega_{\mathrm{si}}^\bullet(U\times \Delta^k)_{\mathrm{vert}})$ is surjective.
\endofproof
\begin{example}
We may parameterize the 2-simplex as
$$
  \Delta^2 = \{ (x,y) \in \mathbb{R}^2 | |x| \leq 1\,,\; 0 \leq  y \leq 1-|x| \}
  \,.
$$
The retraction map $\Delta^2 \to \Lambda^2_1$ in this parameterization is
$$
  (x,y) \mapsto (x,1-|x|)
  \,.
$$
This is smooth away from $x = 0$. A 1-form with sitting instants on $\Lambda^1_1$  vanishes in 
a neighbourhood of $x = 0$, hence its pullback along this map exists and is smooth.
\end{example}

Typically one is interested not in $\exp_\Delta(\mathfrak{g})$ itself, but in a \emph{truncation} thereof.
For our purposes truncation is best modeled by the coskeleton operation.

 Write $\Delta_{\leq n} \hookrightarrow \Delta$ for the full subcategory of the simplex category on the first $n$
  objects $[k]$, with $0 \leq k \leq n$. Write $\mathrm{sSet}_{\leq n}$  for the category of presheaves
  on $\Delta_{\leq n}$. By general abstract reasoning 
  the canonical projection $\mathrm{tr}_n : \mathrm{sSet} \to \mathrm{sSet}_{\leq n}$
  has a left adjoint $\mathrm{sk}_n : \mathrm{sSet}_{\leq n} \to \mathrm{sSet}$ 
  and a right adjoint $\mathrm{cosk}_n : \mathrm{sSet}_{\leq n} \to \mathrm{sSet}$.
  $$
    (\mathrm{sk}_n \dashv \mathrm{tr}_n \dashv \mathrm{cosk}_n) :
    \xymatrix{
       \mathrm{sSet}
       \ar@<+7pt>@{<-}[r]^{\mathrm{sk}_n}
       \ar[r]|{\mathrm{tr}_n}
       \ar@<-7pt>@{<-}[r]_{\mathrm{cosk}_n}
       &
       \mathrm{sSet}_{\leq n}
    }
    \,.
  $$
  The \emph{coskeleton} operation on a simplicial set is the composite
  $$
    \mathbf{cosk}_n := \mathrm{cosk}_n \circ \mathrm{tr}_n : \mathrm{sSet} \to \mathrm{sSet}
    \,.
  $$
  Since $\mathbf{cosk}_n$ is a functor, it extends to an operation of simplicial presheaves,
  which we shall denote by the same symbol
  $$
    \mathbf{cosk}_n : [\mathrm{CartSp}^{\mathrm{op}}, \mathrm{sSet}] \to 
      [\mathrm{CartSp}^{\mathrm{op}}, \mathrm{sSet}]
  $$
  For $X \in \mathrm{sSet}$ or $X \in [\mathrm{CartSp^{\mathrm{op}, \mathrm{sSet}}}]$ 
  we say $\mathbf{cosk}_n X$ is its $n$-\emph{coskeleton}.   
\begin{remark}
Using the adjunction relations, we have that the $k$-cells of $\mathbf{cosk}_n X$ are
images of the $n$-truncation of $\Delta[k]$ in the $n$-truncation of $X$:
$$
  (\mathbf{cosk}_n X)_k 
   = 
  \mathrm{Hom}_{\mathrm{sSet}}(\Delta[k], \mathbf{cosk}_n X)
  =
  \mathrm{Hom}_{\mathrm{sSet}_{\leq n}}(\mathrm{tr}_n \Delta[k], \mathrm{tr}_n X) 
  \,.
$$
\end{remark}
A standard fact (e.g. \cite{dwyer-kan}, \cite{goerss-jardine}) is
\begin{proposition}
  For $X$ a Kan complex
  \begin{itemize}
    \item
      the simplicial homotopy groups $\pi_k$ of $\mathbf{cosk}_n X$ vanishing in degree $k \geq n$;
    \item the canonical morphism $X \to \mathbf{cosk}_n X$ (the unit of the adjunction)
    is an isomorphism on all $\pi_k$ in degree $k < n$;
    \item
      in fact, the sequence
    $$
      X \to \cdots \to 
      \mathbf{cosk}_k X \to \mathbf{cosk}_{k-1} X \to \cdots \to 
      \mathbf{cosk}_1 X \to \mathbf{cosk}_0 X \simeq *
    $$
    is a model for the Postnikov tower of $X$.
   \end{itemize}
\end{proposition}
\begin{example}
For $\mathcal{G}$ a 
groupoid and $N \mathcal{G}$ its simplicial nerve, the canonical morphism
$N \mathcal{G} \to \mathbf{cosk}_2 N \mathcal{G}$ is an isomorphism. 
\end{example}
\begin{definition}
  We say a Kan complex or $L_\infty$-groupoid $X$ is an \emph{$n$-groupoid} if the canonical morphism
  $$
    X \to \mathbf{cosk}_{n+1} X
  $$
  is an isomorphism. If this morphism is just a weak equivalence, we say $X$ is an \emph{$n$-type}.
\end{definition}

We now spell out details of the Lie $\infty$-integration for 
\begin{enumerate}
  \item an ordinary Lie algebra
  \item the \emph{string Lie 2-algebra} 
  \item the line Lie $n$-algebras $b^{n-1}\mathbb{R}$. 
\end{enumerate}
The basic mechanism is is that discused in \cite{henriques}, there for Banach $\infty$-groupoids.
We present now analogous discussions for the context of smooth $\infty$-groupoids that we need 
for the differential refinement in \ref{section.infinity-stack_of_principal_bundles_with_connection} 
and then for the construction of the $\infty$-Chern-Weil homomorphism 
in \ref{section.infinity_chern_weil_homomorphism}

\subsubsection{Ordinary Lie group}

Let $G$ be a Lie group with Lie algebra $\mathfrak{g}$. Then every smooth $\mathfrak{g}$-valued 1-form on the 1-simplex defines an element of $G$ by parallel transport:
\begin{align*}
\mathrm{tra}:\Omega^1_{\mathrm{si}}([0,1],\mathfrak{g})&\to G\\
\omega&\mapsto \mathcal{P}\exp\left(\int_{[0,1]}\omega\right)
\,,
\end{align*}
where the right hand $\mathcal{P} \exp(\cdots)$ is notation
defined to be the endpoint evaluation $f(1)$ of the unique solution $f : [0,1] \to G$ 
to the differential equation
$$
  d f + {r_f}_*(\omega) = 0
$$
with initial condition $f(0) = e$, where $r_g : G \to G$ denotes the right action of $g \in G$ on $G$ itself. 
In the special case that $G$ is simply connected, there is a unique smooth path $\gamma : [0,1] \to G$
starting at the neutral element $e$ 
such that $\omega$ equals the pullback $\gamma^* \theta$ of the Maurer-Cartan form on $G$.
The value of the parallel transport is then the endpoint of this path in $G$.

More generally, this construction works in families and produces 
for every Cartesian space $U$, a parallel transport map
\[
\mathrm{tra}:\Omega^1_{\mathrm{si}}(U\times[0,1],\mathfrak{g})_{\mathrm{vert}}\to C^\infty(U,G)
\]
from smooth $U$-parameterized families of $\mathfrak{g}$-valued 1-forms on the interval to 
smooth functions from $U$ to $G$. If we now consider a $\mathfrak{g}$-valued 1-form $\omega$ on the $n$-simplex instead, parallel transport along the sequence of edges $[0,1]$, $[1,2]$,\dots, $[n-1,n]$ defines an element in $G^{n+1}$, and so we have an induced map $\Omega^1_{\mathrm{si}}(U\times\Delta^n,\mathfrak{g})_{\mathrm{vert}}\to \mathbf{B}G(U)_n$. This map, however is \emph{not} in general a map of simplicial sets: the composition of parallel 
transport along $[0,1]$ and $[1,2]$ is in general not the same as the parallel transport along the edge $[0,2]$ so parallel transport is not compatible with face maps. But precisely if the $\mathfrak{g}$-valued 1-form
is flat does its parallel transport (over the contractible simplex) only depend on the endpoint of the 
path along which it is transported. Therefore we have in particular the following
\begin{proposition}\label{prop.tra.G}
Let $G$ be a Lie group with Lie algebra $\mathfrak{g}$. 
\begin{itemize}
\item Parallel transport along the edges of simplices induces a morphism of smooth $\infty$-groupoids
\[
 \mathrm{tra}:\exp_\Delta(\mathfrak{g})\to \mathbf{B}G.
\]

\item
  When $G$ is simply connected, there is a canonical bijection between smooth flat $\mathfrak{g}$-valued
  1-forms $A$ on $\Delta^n$ and smooth maps $\phi : \Delta^n \to G$ that send the 0-vertex to the neutral element.
  This bijection is given by $A = \phi^* \theta$, where $\theta$ is the Maurer-Cartan form of $G$.
\end{itemize}
\end{proposition}
As for every morphism of Kan complexes, we can look at \emph{coskeletal approximations} of parallel transport given by the morphism of coskeleta towers 
\[
\xymatrix{
  \exp(\mathfrak{g})
    \ar[r]
    \ar[d]_{\mathrm{tra}}
   &
   \cdots\ar[r]&\mathbf{cosk}_{n+1}(\exp(\mathfrak{g}))\ar[d]\ar[r]&\mathbf{cosk}_n(\exp(\mathfrak{g}))\ar[d]\ar[r]&\cdots\ar[r]&{*}\ar[d]\\
\mathbf{B}G\ar[r]&\cdots\ar[r]&\mathbf{cosk}_{n+1}(\mathbf{B}G)\ar[r]&\mathbf{cosk}_n(\mathbf{B}G)\ar[r]&\cdots\ar[r]&{*}
}
\]
\begin{proposition}\label{CoskBGAcyclicFibration}
If the Lie group $G$ is $(k-1)$-connected, then the induced maps
\[
\mathbf{cosk}_n(\exp_\Delta(\mathfrak{g}))\to \mathbf{cosk}_n(\mathbf{B}G)
\]
are acyclic fibrations in $[\mathrm{CartSp}^{\mathrm{op}}, \mathrm{sSet}]_{proj}$ for any $n\leq k$. 
\end{proposition}
\proof
Recall that an acyclic fibration in $[\mathrm{CartSp}^{\mathrm{op}}, \mathrm{sSet}]_{\mathrm{proj}}$ is 
a morphism of simplicial presheaves that is objectwise an acyclic Kan fibration of simplicial sets.
By standard simplicial homotopy theory \cite{goerss-jardine}, the latter are
precisely the maps that have the left lifting property against all 
simplex boundary inclusions $\partial \Delta[p] \hookrightarrow \Delta[p]$.

Notice that for $n = 0$ and $n = 1$ the statement is trivial. For $n \geq 2$ we have
an isomorphism $\mathbf{B}G \to \mathbf{cosk}_n \mathbf{B}G$.  Hence we need to prove that
for $2 \leq n \leq k$ we have for all $U \in \mathrm{CartSp}$ 
lifts $\sigma$ in diagrams of the form
\[
  \xymatrix{
    \partial\Delta[n]\ar[r]\ar[d]_{i} & \exp_\Delta(\mathfrak{g})\ar[d]^{\mathrm{tra}}\\
\Delta[n]\ar[r]\ar@{-->}[ru]^{\sigma} & \mathbf{B}G
}
\]
By parallel transport and using the Yoneda lemma, the outer diagram is equivalently given by 
a map $U \times \partial\Delta^p \to G$ that is smooth with sitting instants on each face $\Delta^{p-1}$. 
By proposition \ref{SmoothGluing} this may be thought of as a smooth map
$U \times S^{p-1} \to G$. The lift $\sigma$ then corresponds to a smooth map with sitting instants 
$\sigma : U \times \Delta^n \to G$ extending this, hence to a smooth map $\sigma : U \times D^p \to G$ 
that in a neighbourhood of $S^{p-1}$ is constant in the direction perpendicular to that boundary.

By the connectivity assumption on $G$ there is a continuous map with these properties. By the
Steenrod-Wockel-approximation theorem \cite{wockel}, this delayed homotopy on a smooth function is itself 
continuously homotopic to a smooth such function. This smooth enhancement of the continuous 
extension is a lift $\sigma$.
\endofproof
 For $n = 1$ the Kan complex $\mathbf{cosk}_{1}(\mathbf{B}G)$ is equivalent to the point. For
$n = 2$ we have an isomorphism $\mathbf{B}G \to \mathbf{cosk}_2 \mathbf{B}G$ (since $\mathbf{B}G$ is the nerve of a Lie groupoid) and so the proposition asserts
that for simply connected Lie groups $\mathbf{cosk}_2\exp_\Delta(\mathfrak{g})$ is equivalent to $\mathbf{B}G$.
\begin{corollary}If $G$ is a compact connected and simply connected Lie group with Lie algebra $\mathfrak{g}$, then the natural morphism $\exp_\Delta(\mathfrak{g})\to \mathbf{B}G$ induces an acyclic fibration $\mathbf{cosk}_3(\exp_\Delta(\mathfrak{g}))\to \mathbf{B}G$ in the global model structure.
\end{corollary}
\begin{proof}Since a compact connected and simply connected Lie group is automatically $2$-connected, we have an induced acyclic fibration $\mathbf{cosk}_3(\exp(\mathfrak{g}))\to\mathbf{cosk}_3(\mathbf{B}G)$. Now notice that $\mathbf{B}G$ is 2-coskeletal, i.e, its coskeleta tower stabilizes at $\mathbf{cosk}_2(\mathbf{B}G)=\mathbf{B}G$.
\end{proof}

\subsubsection{Line $n$-group} 

\begin{definition}
  For $n \geq 1$ write $b^{n-1}\mathbb{R}$ for the 
  \emph{line Lie $n$-algebra}: the  $L_\infty$-algebra characterized by the 
  fact that its Chevalley-Eilenberg algebra is generated from a single generator $c$ in 
  degree $n$ and has trivial differentual 
  $\mathrm{CE}(b^{n-1}\mathbb{R}) = (\wedge^\bullet \langle c\rangle, d = 0)$
  \,.
\end{definition}

\begin{proposition}
  Fiber integration over simplices induces an equivalence
\[
\int_{\Delta^\bullet} : \exp_\Delta(b^{n-1}\mathbb{R}) \stackrel{\simeq}{\to} \mathbf{B}^n\mathbb{R}.
\]
\end{proposition}
\proof
By the Dold-Kan correspondence we only need to show that integration along the simplices is a chain map from the normalized chain complex of $\exp_\Delta(b^{n-1}\mathbb{R})$ to $C^\infty(-)[n]$. The normalized chain complex $N_\bullet(\exp_\Delta(b^{n-1}\mathbb{R}))$ has in degree $-k$ the abelian group $C^\infty(-)\hat \otimes \Omega^n_{cl}(\Delta^k)$, and the differential
\[
\partial: N^{-k}(\exp_\Delta(b^{n-1}\mathbb{R}))\to N^{-k+1}(\exp_\Delta(b^{n-1}\mathbb{R}))
\]
maps a differential form $\omega$ to the alternating sum of its restrictions on the faces of the simplex. If $\omega$ is an element in $C^\infty(-)\otimes \Omega^n_{cl}(\Delta^k)$, integration of $\omega$ on $\Delta^k$ is zero unless $k=n$, which shows that integration along the simplex maps $N^{\bullet}(\exp_\Delta(b^{n-1}\mathbb{R}))$ to $C^\infty(-)[n]$. Showing that this map is actually a map of chain complexes is trivial in all degrees but for $k=n+1$; in this degree, checking that itegration
along simplices is a chain map amounts to checking that for a closed $n$-form
$\omega$ on the $(n+1)$-simplex, the integral of $\omega$ on the boundary of
$\Delta^{n+1}$ vanishes, and this is obvious by Stokes theorem.
\endofproof
\begin{remark}
For $n=1$, the morphism $\exp_\Delta(\mathbb{R})\to \mathbf{B}\mathbb{R}$ coincides with the morphism described in Proposition \ref{prop.tra.G}, for $G=\mathbb{R}$.
\end{remark}

\subsubsection{Smooth string 2-group} 
\label{section.string-2_group}

\begin{definition}
  \label{StringLie2Algebra}
  Let $\mathfrak{string} := \mathfrak{so}_{\mu_3}$ 
  be the extension of the Lie algebra $\mathfrak{so}$ classified
  by its 3-cocycle $\mu_3 = \frac{1}{2}\langle-,[-,-]\rangle$ according to definition \ref{LooExtension}.
  This is called the \emph{string Lie 2-algebra}. 
  Let 
  $$
   \mathbf{B}\mathrm{String} := \mathbf{cosk}_3 \exp_\Delta(\mathfrak{so}_{\mu_3})
  $$
  be its Lie integration. We call this the delooping of the smooth \emph{string 2-group}.
\end{definition}
  The Banach-space $\infty$-groupoid version of this Lie integration is discussed in \cite{henriques}.
\begin{remark}
  \label{7CocycleOnStringLie2Algebra}
  The 7-cocycle $\mu_7$ on $\mathfrak{so}$
  is still, in the evident way, a cocycle on $\mathfrak{so}_{\mu_3}$
  $$
    \mu_7 : \mathfrak{so}_3 \to b^6 \mathbb{R}
    \,.
  $$
\end{remark}
\begin{proposition}
  \label{StringByPullback}
  We have an $\infty$-pullback of smooth $\infty$-groupoids
  $$
    \xymatrix{
      \mathbf{B}\mathrm{String} \ar[d]\ar[r] & {*} \ar[d]
      \\
      \mathbf{B}\mathrm{Spin} \ar[r]^{\frac{1}{2}\mathbf{p}_1}
      &
      \mathbf{B}^3 U(1)
    }
  $$
  presented by the ordinary pullback of simplicial presheaves  
  $$
    \xymatrix{
      \mathbf{B}\widetilde{\mathrm{String}}
      \ar[d]
      \ar[r]
      &
      \mathbf{E}\mathbf{B}^2 (\mathbb{Z} \to \mathbb{R})
      \ar[d]
      \\
      \mathbf{cosk}_3 \exp_\Delta(\mathfrak{so})
      \ar[r]^{\exp_\Delta(\mu_3)}
      \ar[d]^{\wr}
      &
      \mathbf{B}^3 (\mathbb{Z} \to \mathbb{R})
      \\
      \mathbf{B}\mathrm{Spin}
    }
    \,,
  $$
  where $\mathbf{B}\mathrm{String} \stackrel{\sim}{\to} \mathbf{B}\widetilde {\mathrm{String}}$ is 
  induced by integrating the 2-form over simplices.
\end{proposition}
\begin{remark} 
  In terms of definition \ref{definition.AbelianNBundles},
  $\mathbf{B}\mathrm{String}$ is the smooth $\mathbf{B}^2 U(1)$-principal 3-bundle over 
  $\mathbf{B}\mathrm{Spin}$ classified by the smooth refinement of the first fractional
  Pontryagin class.
\end{remark}
\proof
  Since all of $\mathbf{cosk}_3 \exp_\Delta(\mathfrak{so})$, $\mathbf{B}^3 (\mathbb{Z} \to \mathbb{R})$ and
  $\mathbf{E} \mathbf{B}^2 (\mathbb{Z} \to \mathbb{R})$ are fibrant in 
  $[\mathrm{CartSp}^{op}, \mathrm{sSet}]_{\mathrm{proj}}$ and since 
  $\mathbf{E}\mathbf{B}^2 U(1) \to \mathbf{B}^3 (\mathbb{Z} \to \mathbb{R})$ is a fibration (being the image
  under DK of a surjection of complexes of sheaves), we have by standard facts about homotopy
  pullbacks that the ordinary pullback is a homotopy pullback in 
  $[\mathrm{CartSp}^{\mathrm{op}}, \mathrm{sSet}]_{\mathrm{proj}}$. By \cite{lurie} this 
  presents the $\infty$-pullback of $\infty$-presheaves on $\mathrm{CartSp}$. And since 
  $\infty$-stackification is left exact, this is also presents the $\infty$-pullback of 
  $\infty$-sheaves.
  
  This ordinary pullback manifestly has 2-cells given by 2-simplices in $G$ labeled by 
  elements in $U(1)$ and 3-cells being 3-simplices in $G$ such that the labels of their faces
  differ by $\int_{\Delta^3 \to G} \mu \; \mathrm{mod} \mathbb{Z}$. This is 
  the definition of   $\mathbf{B}\mathrm{String}$. That $\exp_\Delta(\mu_3)$ indeed presents
  a smooth refinement of the second fractional Pontryagin class as indicated is shown below.
\endofproof
\begin{proposition}
  \label{EquivalentModelsForStringGroup}
  There is a zig-zag of equivalences
  $$
    \mathbf{cosk}_3 \exp(\mathfrak{so}_{\mu_3})
    \simeq \cdots \simeq \mathbf{B}(\hat \Omega \mathrm{Spin} \to P \mathrm{Spin})
  $$
  in $[\mathrm{CartSp}^{\mathrm{op}}, \mathrm{sSet}]_{\mathrm{proj}}$,
  of the Lie integration,
  prop. \ref{StringByPullback}, of $\mathfrak{so}_{\mu_3}$ 
  with the strict 2-group, def. \ref{CrossedModuleAndStrict2Group} 
  coming from the crossed module $(\hat \Omega \mathrm{Spin} \to P \mathrm{Spin})$ 
  of Fr{\'e}chet Lie groups, discussed in \cite{bcss}, consisting of the centrally extended loop group
  and the path group of $\mathrm{Spin}$.
\end{proposition}
This is proven in section 4.2 of \cite{survey}.
\begin{proposition}
  The object $\mathbf{B}\mathrm{String} = \mathbf{cosk}_3(\exp(\mathfrak{so}_{\mu_3}))$ 
  is fibrant in $[\mathrm{CartSp}^{\mathrm{op}}, \mathrm{sSet}]_{\mathrm{proj}, \mathrm{loc}}$.
\end{proposition}
\proof
  Observe first that both object are fibrant in 
  $[\mathrm{CartSp}^{\mathrm{op}}, \mathrm{sSet}]_{\mathrm{proj}}$ (the Lie integration by
  prop. \ref{proposition.exp.is.infty.groupoid}, the delooped strict 2-group by 
  observation \ref{GlobalFibrancyOfBOf2Group}).
  The claim then follows with prop. \ref{EquivalentModelsForStringGroup} and prop. 
  \ref{LocalFibrancyOfBOf2Group}, which imply that for $C(\{U_i\}) \to \mathbb{R}^n$
  the {\v C}ech nerve of a good open cover, hence a cofibrant resultions, there is
  a homotopy equivalence
  $$
    [\mathrm{CartSp}^{\mathrm{op}}, \mathrm{sSet}](\check{C}(\mathcal{U}), \mathbf{cosk}_3 \exp(\mathfrak{so}_{\mu_3}))
    \simeq
    [\mathrm{CartSp}^{\mathrm{op}}, \mathrm{sSet}](\check{C}(\mathcal{U}), \mathbf{B}(\hat \Omega \mathrm{Spin} \to P \mathrm{Spin}))
     \,.
  $$
\endofproof
\begin{corollary} 
A Spin-principal bundle $P\to X$ can be lifted to a String-principal bundle 
precisely if it 
 trivializes $\frac{1}{2}\mathbf{p}_1$, i.e., if the induced mophism
$\mathbf{H}(X,\mathbf{B}Spin) \to \mathbf{H}(X,\mathbf{B}^3U(1))$ is
homotopically trivial. The choice of such a lifting is called a \emph{String structure} on the Spin-bundle.
\end{corollary}
We discuss string structures and their twisted versions further in \ref{HomotopyFibers}.

\subsection{Principal $\infty$-bundles with connection}
\label{section.infinity-stack_of_principal_bundles_with_connection}

For an ordinary Lie group $G$ with Lie algebra $\mathfrak{g}$, we have met 
in section \ref{section.BG_BG-conn_and_principal_G-bundles_with_connection}
the smooth groupoids $\mathbf{B}G$, $\mathbf{B}G_{\mathrm{conn}}$ and $\mathbf{B}G_{\mathrm{diff}}$ arising from $G$, and  in \ref{section.principal_infinity-bundles}  
the smooth $\infty$-groupoid $\exp_\Delta(\mathfrak{g})$ coming from $\mathfrak{g}$, and have shown that they are related by  a diagram
\[
\xymatrix{
& & \exp_\Delta(\mathfrak{g})\ar[d]\\
\mathbf{B}G_{\mathrm{conn}}\ar@{^{ (}->}[r]&\mathbf{B}G_{\mathrm{diff}}\ar@{->>}[r]^{\sim}&\mathbf{B}G
}
\]
and that $\mathbf{B}G_{\mathrm{conn}}$ is the moduli stack of $G$-principal bundles with connection.
Now we dicuss such \emph{differential refinements} 
$\exp_\Delta(\mathfrak{g})_{\mathrm{diff}}$ and $\exp_\Delta(\mathfrak{g})_{\mathrm{conn}}$ that complete
the above diagram for any integrated smooth $\infty$-group $\exp_\Delta(\mathfrak{g})$. 
Where a truncation of $\exp_\Delta(\mathfrak{g})$ is the object that classifies
$G$-principal $\infty$-bundles, the corresponding truncation of $\exp_\Delta(\mathfrak{g})_{\mathrm{conn}}$
classifies \emph{principal $\infty$-bundles with connection}. Between $\exp_\Delta(\mathfrak{g})_{\mathrm{conn}}$ and $\exp_\Delta(\mathfrak{g})_{\mathrm{diff}},$ we will also meet the \emph{Chern-Weil $\infty$-groupoid} $\exp_\Delta(\mathfrak{g})_{\mathrm{CW}}$ which is the natural ambient for $\infty$-Chern-Weil theory to live in.

For the following, let $\mathfrak{g}$ be any Lie $\infty$-algebra. 
\begin{definition}
The differential refinement $\exp_\Delta(\mathfrak{g})_{\mathrm{diff}}$ of  $\exp_\Delta(\mathfrak{g})$ is the simplicial presheaf on the site of Cartesian spaces given by the assignment
\[
  (U,[k])
  \mapsto
  \left\{
    \raisebox{24pt}{
   \xymatrix{\Omega_{\mathrm{si}}^\bullet(U\times\Delta^k)_{\mathrm{vert}} 
     & \mathrm{CE}(\mathfrak{g})\ar[l]_{\phantom{mmm}A_{\mathrm{vert}}}\\
\Omega_{\mathrm{si}}^\bullet(U\times \Delta^k)\ar[u] & W(\mathfrak{g})\ar[l]_{\phantom{mmm}A}\ar[u]
}
  }
  \right\}
  \,,
\]
  where on the right we have the set of commuting diagrams in $\mathrm{dgcAlg}$ as indicated.
\end{definition}
\begin{remark}This means that a $k$-cell in $\exp_\Delta(\mathfrak{g})_{\mathrm{diff}}$ 
over $U \in \mathrm{CartSp}$ is 
a $\mathfrak{g}$-valued form $A$ on $U \times \Delta^k$ that satisfies the condition that
its curvature forms $F_A$ vanish when restricted in all arguments to vectors on the simplex.
This is the analog of the \emph{first Ehresmann condition} on a connection form
on an ordinary principal bundle: the form $A$ on the trivial
simplex bundle $U \times \Delta^k \to U$ is flat along the fibers.
\end{remark}
\begin{proposition}\label{prop.acyclic-exp}
The evident morphism of simplicial presheaves 
$$
  \xymatrix{
    \exp_\Delta(\mathfrak{g})_{\mathrm{diff}}
     \ar@{->>}[r]^{\sim}
      &
   \exp_\Delta(\mathfrak{g})
  }
$$ 
is an acyclic fibration of smooth $\infty$-groupoids in the global model structure.
\end{proposition}
\proof
We need to check that, for all $U \in \mathrm{CartSp}$ and $[k] \in \Delta$ and for all 
diagrams
$$
  \xymatrix{
    \partial \Delta[k] \ar[r]^{A|_\partial} \ar[d] & \exp_\Delta(\mathfrak{g})_{\mathrm{diff}}(U) \ar[d]
    \\
    \Delta[k] \ar[r] \ar[r]^{A_{\mathrm{vert}}} 
    \ar@{-->}[ur]|{A}
    & \exp_\Delta(\mathfrak{g})(U)
  }
$$
 we have a lift as indicated by the dashed morphism. For that we need to extend the composite
\[
\mathrm{W}(\mathfrak{g})\to \mathrm{CE}(\mathfrak{g}) \stackrel{A_{\mathrm{vert}}}{\to} 
\Omega_{\mathrm{si}}^\bullet(U\times\Delta^n)_{\mathrm{vert}}
\]
to an element in 
$\Omega_{\mathrm{si}}^\bullet(U\times\Delta^k)\otimes\mathfrak{g}$ with fixed boudary value $A_\partial$ in $\Omega_{\mathrm{si}}^\bullet(U\times\partial\Delta^k)\otimes\mathfrak{g}$. 
To see that this is indeed possible, use the decomposition
\[
  \Omega_{\mathrm{si}}^\bullet(U\times\Delta^k)
  =
  \Omega_{\mathrm{si}}^\bullet(U\times\Delta^n)_{\mathrm{vert}}
   \oplus
  \left(\Omega^{>0}(U)\hat{\otimes}\Omega_{\mathrm{si}}^\bullet(\Delta^k)\right)
\]
to write $A_\partial=A_{\mathrm{vert}}|_{\partial\Delta^k}+A^{>0}_\partial$. 
Extend $A^{>0}_\partial$ to an element $A^{>0}$ in $\Omega^{>0}(U)\hat{\otimes}\Omega_{\mathrm{si}}^\bullet(\Delta^k)$. 
This is a trivial extension problem: any smooth differential form on the boundary of an $k$-simplex can be extended to a smooth differential form on the whole simplex. 
Then the degree 1 element $A_{\mathrm{vert}}+A^{>0}$ is a solution to our original extension problem.
\endofproof
\begin{remark}This means that $\exp_\Delta(\mathfrak{g})_{\mathrm{diff}}$ is a certain resolution of
$\exp_\Delta(\mathfrak{g})$. In the full abstract theory \cite{survey}, the reason for its existence is that it serves
to model the canonical curvature characteristic map $\mathbf{B}G \to \mathbf{\flat}_{\mathrm{dR}}\mathbf{B}^n U(1)$
in the $\infty$-topos of smooth $\infty$-groupoids by a truncation of the zig-zag
$\exp_\Delta(\mathfrak{g}) \stackrel{\sim}{\leftarrow} \exp_\Delta(\mathfrak{g})_{\mathrm{diff}}
\to \exp_\Delta(b^{n-1}\mathbb{R})$ of simplicial presheaves. By the nature of acyclic fibrations,
we have that for every $\exp_\Delta(\mathfrak{g})$-cocycle $X \xleftarrow{\sim_\mathrm{loc}}\check{C}(\mathcal{U}) 
\stackrel{g}{\to} 
\exp_\Delta(\mathfrak{g})$ there is a lift $g_{\mathrm{diff}}$ to an $\exp_\Delta(\mathfrak{g})_{\mathrm{diff}}$-cocycle
$$
  \xymatrix{
    QX\ar[d]^{\wr} \ar[r]^{g_{\mathrm{diff}}} & \exp_\Delta(\mathfrak{g})_{\mathrm{diff}} \ar[d]^{\wr}
    \\
    \check{C}(\mathcal{U})\ar[d]^{\begin{turn}{270}$\scriptstyle{\sim_\mathrm{loc}}$\end{turn}} \ar[r]^g & \exp_\Delta(\mathfrak{g})
    \\
    X
  }
  \,.
$$
For the abstract machinery of $\infty$-Chern-Weil theory to work, it is only the existence of this
lift that matters. However, in practice it is useful to make certain nice choices of lifts. 
In particular, when $X$ is a paracompact smooth manifold, there is always a choice of lift with the
property that the corresponding curvature characteristic forms are globally defined forms on $X$, instead
of more general (though equivalent) cocycles in total \v{C}ech-de Rham cohomology. 
Moreover, in this case the local connection forms can be chosen to have $\Delta$-horizontal curvature. Lifts with this special property are \emph{genuine $\infty$-connections} 
on the $\infty$-bundles classified by $g$. 
The following definitions formalize this. But it is important to note that genuine $\infty$-connections
are but a certain choice of gauge among all differential lifts. Notably when the base $X$ is not a manifold
but for instance a non-trivial orbifold, then genuine $\infty$-connections will in general not even exist,
whereas the differential lifts always do exist, and always support the $\infty$-Chern-Weil homomorphism.
\end{remark}
\begin{proposition}
If $\mathfrak{g}$ is the Lie algebra of a Lie group $G$, then there is a natural commutative diagram
\[
\xymatrix{
\exp_\Delta(\mathfrak{g})_\mathrm{diff}\ar@{->>}[r]^{\sim}\ar[d]& \exp_\Delta(\mathfrak{g})\ar[d]\\
\mathbf{B}G_\mathrm{diff}\ar@{->>}[r]^{\sim}&\mathbf{B}G
}
\]
In particular, if $G$ is $(k-1)$-connected, with $k\geq 2$, then the induced morphism $\mathrm{cosk}_k(\exp_\Delta(\mathfrak{g})_\mathrm{diff})\to\mathbf{B}G_\mathrm{diff}$ is an acyclic fibration in the global model structure. 
\end{proposition}
\proof
We have seen in Remark \ref{remark.BG_diff} that there is a natural isomorphism 
$\mathbf{B}G_{\mathrm{diff}}\cong \mathbf{B}G\times\mathrm{Codisc}(\Omega^1(-;\mathfrak{g})$, so in order to give the morphism $\exp_\Delta(\mathfrak{g})_\mathrm{diff}\to \mathbf{B}G_\mathrm{diff}$ making the above diagram commute we only need to give a natural morphism $\exp_\Delta(\mathfrak{g})_\mathrm{diff}\to\mathrm{Codisc}(\Omega^1(-;\mathfrak{g})$; this is evaluation  of the connection form
$A$ on the vertices of the simplex.
\par
Assume now $G$ is $(k-1)$-connected, with $k\geq 2$. Then, by Propositions \ref{CoskBGAcyclicFibration} and \ref{prop.acyclic-exp}, both $\mathrm{cosk}_k(\exp_\Delta(\mathfrak{g})_\mathrm{diff}) \to \mathrm{cosk}_k(\exp_\Delta(\mathfrak{g}))$ and $\mathrm{cosk}_k(\exp_\Delta(\mathfrak{g}))\to \mathbf{B}G$ are acyclic fibrations. We have a commutative diagram 
\[
\xymatrix{
\mathrm{cosk}_k(\exp_\Delta(\mathfrak{g})_\mathrm{diff})\ar@{->>}[r]^{\sim}\ar[d]& \mathrm{cosk}_k(\exp_\Delta(\mathfrak{g}))\ar@{->>}[d]^{\wr}\\
\mathbf{B}G_\mathrm{diff}\ar@{->>}[r]^{\sim}&\mathbf{B}G,
}
\]
so, by the ``two out of three'' rule, also $\mathrm{cosk}_k(\exp_\Delta(\mathfrak{g})_\mathrm{diff})\to\mathbf{B}G_\mathrm{diff}$ is an acyclic fibration.
\endofproof

\begin{definition}The simplicial presheaf  
$\exp_\Delta(\mathfrak{g})_{\mathrm{CW}} \subset \exp_\Delta(\mathfrak{g})_{\mathrm{diff}}$ 
is the sub-presheaf of $\exp_\Delta(\mathfrak{g})_{\mathrm{diff}}$ on those $k$-cells 
$\Omega_{\mathrm{si}}^\bullet(U \times \Delta^k) \stackrel{A}{\leftarrow} \mathrm{W}(\mathfrak{g}) $ that make 
also the bottom square in the diagram
\[
\xymatrix{\Omega_{\mathrm{si}}^\bullet(U\times \Delta^k)_{\mathrm{vert}} 
 & CE(\mathfrak{g})\ar[l]_{\phantom{mmm}A_{\mathrm{vert}}}\\
\Omega_{\mathrm{si}}^\bullet(U\times \Delta^k)\ar[u] & W(\mathfrak{g})\ar[l]_{\phantom{mmm}A}\ar[u]\\
\Omega^\bullet(U)\ar[u] & \mathrm{inv}(\mathfrak{g})\ar[u]\ar[l]_{\phantom{mmm} F_A}
}
\]
commute.
\end{definition}
\begin{remark}
A $k$-cell in $\exp(\mathfrak{g})_{\mathrm{CW}}$ parameterised by a Cartesian space $U$ is
a $\mathfrak{g}$-valued differential form $A$ on the total space $U \times \Delta^k$
such that 
\begin{enumerate}
  \item its restriction to the fiber $\Delta^k$ of $U \times \Delta^k \to U$ is flat, and indeed equal to the canonical $\mathfrak{g}$-valued form there as encoded by the cocycle $A_{\mathrm{vert}}$ (which, recall, is the datum in   $\exp(\mathfrak{g})$ that determines the $G$-bundle itself);
this we may think of as the \emph{first Ehresmann condition} on a connection;

  \item all its curvature characteristic forms $\langle F_A \rangle$ descend to the base space $U$ of $U \times \Delta^k \to U$; this we may think of as a slightly weakened version of the \emph{second Ehresmann condition} on a connection: this is the main \emph{consequence} of the second Ehresmann condition. 
\end{enumerate}
These are the structures that have been considered in \cite{SSSI} and \cite{SSSIII}.
\end{remark}
\begin{proposition}
  $\exp_\Delta(\mathfrak{g})_{\mathrm{CW}}$ is 
  fibrant in $[\mathrm{CartSp}^{\mathrm{op}}, \mathrm{sSet}]_{\mathrm{proj}}$ 
  \,.
\end{proposition}
\proof
As in the proof of Propositition \ref{proposition.exp.is.infty.groupoid} 
we find horn fillers $\sigma$ by pullback along the standard retracts, which are smooth
away from the loci where our forms have sitting instants. 
\[
\xymatrix{
  \Omega_{\mathrm{si}}^\bullet(U\times \Delta^k)_{\mathrm{vert}} 
  &
   \Omega_{\mathrm{si}}^\bullet(U\times \Lambda^k_i)\ar[l] & 
  \mathrm{CE}(\mathfrak{g})\ar[l]_{\phantom{mmm}A_{\mathrm{vert}}}
  \\
  \Omega_{\mathrm{si}}^\bullet(U\times \Delta^k)\ar[u] 
   &
   \Omega_{\mathrm{si}}^\bullet(U\times \Lambda^k_i)\ar[l]\ar[u] 
   & 
   \mathrm{W}(\mathfrak{g})\ar[l]_{\phantom{mmm}A}\ar[u]
    & : \sigma
  \\
\Omega^\bullet(U)\ar[u] &\Omega^\bullet(U)\ar[u]\ar[l] & 
  \mathrm{inv}(\mathfrak{g})\ar[u]\ar[l]_{\phantom{mmm}F_A}
}
\]
\endofproof
We say that $\exp(\mathfrak{g}_\mu)_{\mathrm{CW}}$ is the \emph{Chern-Weil $\infty$-groupoid of $\mathfrak{g}$}.
\begin{definition}
  Write $\exp_\Delta(\mathfrak{g})_{\mathrm{conn}}$
  for the simplicial sub-presheaf of $\exp_\Delta(\mathfrak{g})_{\mathrm{diff}}$ given in degree $k$ by those $\mathfrak{g}$-valued forms
  satisfying the following further \emph{horizontality} condition:
  \begin{itemize}
  \item for all vertical (i.e., tangent to the simplex) vector fields $v$ on $U \times \Delta^k$, we have
  $$
    \iota_v F_A = 0
    \,.
  $$
  \end{itemize}
\end{definition}
\begin{remark}
This extra condition is the direct analog of the \emph{second Ehresmann condition}.
For ordinary Lie algebras we have discussed this form of the second Ehresmann condition in 
section \ref{OrdinaryConnections}.
\end{remark}
\begin{remark}
If we decompose differential forms on the products $U \times \Delta^k$ as
$$
  \Omega_{\mathrm{si}}^\bullet(U \times \Delta^k) = \bigoplus_{p,q \in \mathbb{N}} \Omega^p(U) \hat \otimes \Omega_{\mathrm{si}}^q(\Delta^k)
$$
then 
\begin{enumerate}
  \item $k$-simplices in $\exp_\Delta(\mathfrak{g})_{\mathrm{diff}}$ are those connection forms $A:W(\mathfrak{g})\to \Omega_{\mathrm{si}}^\bullet(U\times\Delta^k)$ whose curvature form has only the $(p,q)$-components with $p>0$;
  \item $k$-simplices in $\exp_\Delta(\mathfrak{g})_{\mathrm{conn}}$ are those $k$-simplices in $\exp_\Delta(\mathfrak{g})_{\mathrm{diff}}$ whose curvature is 
   furthermore constrained to have precisely only the $(p,0)$-components, with $p>0$.
\end{enumerate}
\end{remark}
\begin{proposition}
  We have a sequence of inclusions of simplicial presheaves
  $$
    \exp_\Delta(\mathfrak{g})_{\mathrm{conn}}
    \hookrightarrow
    \exp_\Delta(\mathfrak{g})_{\mathrm{CW}}
    \hookrightarrow
    \exp_\Delta(\mathfrak{g})_{\mathrm{diff}}
    \,.
  $$
\end{proposition}
\proof
  Let $\langle-\rangle$ be an invariant polynomial on $\mathfrak{g}$, and $A$ a $k$-cell of $\exp_\Delta(\mathfrak{g})_{\mathrm{conn}}$. Since $d_{W(\mathfrak{g})}\langle-\rangle=0$,  we have 
  $d \langle F_A \rangle = 0$, and since $\iota_vF_A=0$ we also have $\iota_v \langle F_A \rangle = 0$ for $v$ tangent to the
  $k$-simplex.
  Therefore by Cartan's formula also the Lie derivatives $\mathcal{L}_v \langle F_A \rangle $ are zero.
  This implies that the curvature characteristic forms on $\exp_\Delta(\mathfrak{g})_{\mathrm{conn}}$ descend to $U$
  and hence that $A$ defines a $k$-cell in $\exp_\Delta(\mathfrak{g})_{\mathrm{CW}}$.
\endofproof

\subsubsection{Examples}

We consider the special case of the above general construction again for the 
special case that $\mathfrak{g}$ is an ordinary
Lie algebra and for $\mathfrak{g}$ of the form $b^{n-1}\mathbb{R}$.

\begin{proposition} Let $G$ be a Lie group with Lie algebra $\mathfrak{g}$. 
Then, for any $k \in \mathbb{N}$ there is a  pullback diagram
\[
\xymatrix{
  \mathbf{cosk}_k\exp_\Delta(\mathfrak{g})_{\mathrm{conn}}\ar[d]\ar[r] 
  & 
  \mathbf{cosk}_k\exp_\Delta(\mathfrak{g})\ar[d]\\
\mathbf{B}G_{\mathrm{conn}}\ar[r]& \mathbf{B}G
}
\]
in the category of simplicial presheaves.
\end{proposition}
\proof 
  The result is trivial for $n=0$. For $n=1$ we have to show that given two 1-forms $A_0,A_1\in \Omega^1(U,\mathfrak{g})$, a gauge transformation $g:U\to G$ between them, and  any lift $\lambda(u,t)dt$ of $g$ to a 1-form in $\Omega^1(U\times\Delta^1, \mathfrak{g})$, there exists a unique 1-form $A\in\Omega^1(U\times\Delta^1, \mathfrak{g})$ whose vertical part is $\lambda$, whose curvature is of type (2,0), and such that
\[
A\bigl\vert_{U\times\{0\}}=A_0;\qquad A\bigl\vert_{U\times\{1\}}=A_1.
\]
We may decompose such $A$ into its vertical and horizontal components
\[
  A = \lambda \, dt + A_U
  \,,
\]
where $\lambda \in C^\infty(U \times\Delta^1)$ and $A_U$ in the image of $\Omega^1(U, \mathfrak{g})$.
Then the horizontality condition $\iota_{\partial_t} F_A = 0$ on $A$ is the differential equation
\[
  \frac{\partial}{\partial t}A_U 
   =
   d_U \lambda + [A_U,\lambda]
   \,.
\]
For the given initial condition $A_U(t = 0) = A_0$ this has a unique solution, given by
\[
  A_U(t) =g(t)^{-1}A_0 g(t) + g(t)^{-1} d_t g(t),
  \]
where $g(t)\in G$ is the parallel transport 
for the connection $\lambda\, dt$ along the path $[0,t]$ in the 1-simplex $\Delta^1$. 
Evaluating at $t=1$, and using $g(1)=g$, we find 
\[
A(1)=g^{-1}A_0 g +g^{-1}d g = A_1,
\]
as required. 

These computations carry on without substantial modification to higher simplices: using that 
$\lambda\, d t$ is required to be flat along the simplex, it follows the value of $A_U$ at any point in the simplex
is determined by a differential equation as above, for parallel transport along \emph{any} path
from the 0-vertex to that point. Accordingly we find unique lifts $A$, which concludes the proof.
\endofproof
\begin{corollary}
If $G$ is a compact simply connected Lie group, there is a weak equivalence $\mathbf{cosk}_3\exp(\mathfrak{g})_{\mathrm{conn}}\xrightarrow{\sim}\mathbf{B}G_{\mathrm{conn}}$
in $[\mathrm{CartSp}^{\mathrm{op}}, \mathrm{sSet}]_{\mathrm{proj}}$.
\end{corollary}
\proof
  By proposition \ref{CoskBGAcyclicFibration} 
  we have that $\mathbf{cosk}_3 \exp\Delta(\mathfrak{g}) \to \mathbf{B}G$ is an 
  acyclic fibration in the global model structure. Since these are preserved under pullback, the claim follows by the 
  above proposition.
\endofproof

\begin{proposition}
Integration along simplices gives a morphism of smooth $\infty$-groupoids 
\[
\int_{\Delta^\bullet}^{\mathrm{diff}}:\exp_\Delta(b^{n-1}\mathbb{R})_{\mathrm{diff}}\to 
  \mathbf{B}^n\mathbb{R}_{\mathrm{diff}}.
\]
\end{proposition}
\proof
By means of the Dold-Kan correspondence we only need to show that integration along simplices is a morphism of complexes from the normalized chain complex of $\exp_\Delta(b^{n-1}\mathbb{R})$ to the cone
\begin{equation}\label{displayed-cone}
\xymatrix@R=1pt{
       C^\infty(U) 
        \ar[r]^{d} & 
       \Omega^1(U)
       \ar[r]^d & \Omega^2(U) 
       \ar[r]^d & \cdots
        \ar[r]^d
        &
        \Omega^n(U)
       \\
       \oplus & \oplus & \oplus&\oplus
       \\
       \Omega^1(U) \ar[uur]|{\mathrm{Id}} \ar[r]_d
       & \Omega^2(U)
         \ar@{--}[uur]
         \ar[r]
         & \cdots
         \ar[r]
         &
         \Omega^n(U)
         \ar[uur]|{\mathrm{Id}}
         \ar[r]
         &0
    }
\end{equation}
The normalized chain complex $N^\bullet(\exp_\Delta(b^{n-1}\mathbb{R}))$ has in degree $-k$ the subspace of $\Omega^n(U\times \Delta^k)$ consisting of those $n$-forms whose $(0,n)$-component $\omega_{0,n}$ lies in $C^\infty(U)\hat{\otimes}\Omega^n_{\mathrm{cl}}(\Delta^k)$; the differential
\[
\partial:N^{-k}(\exp_\Delta(b^{n-1}\mathbb{R}))\to N^{-k+1}(\exp_\Delta(b^{n-1}\mathbb{R}))
\]
maps an $n$-form $\omega$ on $U\times \Delta^k$ to the alternate sum of its restrictions to the faces of $U\times \partial\Delta^k$. 
For $k\neq 0$, let 
$\int_{\Delta^\bullet}^{diff}$ be the map
\begin{align*}
\int_{\Delta^k}^{diff}:N^{-k}(\exp_\Delta(b^{n-1}\mathbb{R}))&\to \Omega^{n-k}(-)\oplus \Omega^{n-k+1}(-)\\
\omega&\mapsto\left(\int_{\Delta^k}\omega, \int_{ \Delta^k}d_{dR}\omega\right),
\end{align*}
and, for $k=0$ let $\int_{\Delta^0}^{diff}$ be the identity
\[
\int_{\Delta^0}^{diff}=\mathrm{id}:N^{0}(\exp_\Delta(b^{n-1}\mathbb{R}))\to \Omega^{n}(-).
\]
The map $\int_{\Delta^\bullet}^{diff}$ actually takes its values in the cone (\ref{displayed-cone}). Indeed, if $k>n+1,$ then both the integral of $\omega$ and of $d_{dR}\omega$ are zero by dimensional reasons; for $k=n+1$, the only possibly nontrivial contribution to the integral over $\Delta^{n+1}$ comes from $d_{\Delta^{n+1}}\omega_{0,n}$, which is zero by hypothesis (where we have written $d_{dR}=d_{\Delta^k}+d_U$ for the decomposition of the de Rham differential associated with the product structure of $U\times\Delta^k$).
\par
The fact that $\int_{\Delta^\bullet}^{diff}$ is a chain map immediately follows by the Stokes formula:
\[
\int_{\Delta^{k-1}}\partial\omega=\int_{\partial\Delta^k}\omega=\int_{\Delta^k}d_{\Delta^k}\omega
\]
and by the identity $d_{dR}=d_{\Delta^k}+d_U$.
\endofproof

\begin{corollary}
Integration along simplices induces a morphism of smooth $\infty$-groupoids 
\[
\int_{\Delta^\bullet}^{\mathrm{conn}}:\exp_\Delta(b^{n-1}\mathbb{R})_{conn}\to \mathbf{B}^n\mathbb{R}_{conn}.
\]
\end{corollary}
\proof
By Proposition \ref{prop.curv}, we only need to check that the image of the composition $\mathrm{curv}\circ\int_{\Delta^\bullet}^{\mathrm{diff}}$ lies in the subcomplex $(0 \to 0 \to \cdots \to \Omega^{n+1}_{\mathrm{cl}}(-))$ of $\mathbf{\flat}_{\mathrm{dR}} \mathbf{B}^{n+1}\mathbb{R}$, and this is trivial since by definition of $\exp_\Delta(b^{n-1}\mathbb{R})_{conn}$ the curvature of $\omega$, i.e.,  the de Rham differential $d_{dR}\omega$, is 0 along the simplex. 
\endofproof

\section{$\infty$-Chern-Weil homomorphism}
\label{section.infinity_chern_weil_homomorphism}

With the constructions that we have introduced in the previous sections, there is an evident
Lie integration of a cocycle $\mu : \mathfrak{g} \to b^{n-1}\mathbb{R}$ on a
$L_\infty$-algebra $\mathfrak{g}$ to a morphism
$\exp_\Delta(\mathfrak{g}) \to \exp_\Delta(b^{n-1}\mathbb{R})$ that truncates to a
characteristic map $\mathbf{B}G \to \mathbf{B}^n \mathbb{R}/\Gamma$. 
Moreover, this has an evident lift to a morphism 
$\exp_\Delta(\mathfrak{g})_{\mathrm{diff}} \to \exp_\Delta(b^{n-1}\mathbb{R})_{\mathrm{diff}}$ 
between the differential refinements. Truncations of this we shall now identify with the 
Chern-Weil homomorphism and its higher analogs.

\subsection{Characteristic maps by $\infty$-Lie integration}
\label{section.characteristic_maps_by_Lie_integration}

We have seen in section \ref{ooLieIntegration}
how  $L_\infty$-algebras $\mathfrak{g}$, $b^{n-1}\mathbb{R}$ integrate to 
smooth $\infty$-groupoids $\exp_\Delta(\mathfrak{g})$, $\exp_\Delta(b^{n-1}\mathbb{R})$ and their
differential refinements $\exp_\Delta(\mathfrak{g})_{\mathrm{diff}}$, $\exp_\Delta(b^{n-1}\mathbb{R})_{\mathrm{diff}}$
as well as to various truncations and quotients of these.
We remarked at the end of \ref{ooLieAlgebroid} that a degree $n$ cocycle $\mu$ on $\mathfrak{g}$ 
may equivalently be thought of as a morphism $\mu : \mathfrak{g} \to b^{n-1}\mathbb{R}$, i.e.,  as a 
dg-algebra morphism $\mu:\mathrm{CE}(b^{n-1}\mathbb{R})\to \mathrm{CE}(\mathfrak{g})$.
\begin{definition}
  Given an $L_\infty$-algebra cocycle 
  $\mu : \mathfrak{g} \to b^{n-1}\mathbb{R}$ as in section \ref{section.Lie_infinity-algebroids}, 
  define a morphism of simplicial presheaves
  $$
    \exp_{\Delta}(\mu) : \exp_{\Delta}(\mathfrak{g}) \to \exp_\Delta(b^{n-1}\mathbb{R})
  $$
  by componentwise composition with $\mu$:
  $$
    \begin{aligned}
    \exp_\Delta(\mu)_k :     
    &\left(\Omega_{\mathrm{si}}^\bullet(U \times \Delta^k)_{\mathrm{vert}}
    \stackrel{A_{\mathrm{vert}}}{\leftarrow} \mathrm{CE}(\mathfrak{g})\right)
    \mapsto
    \\
    &\left(\Omega_{\mathrm{si}}^\bullet(U \times \Delta^k)_{\mathrm{vert}}
    \stackrel{A_{\mathrm{vert}}}{\leftarrow} \mathrm{CE}(\mathfrak{g})
    \stackrel{\mu}{\leftarrow} \mathrm{CE}(b^{n-1}\mathbb{R}) : \mu(A_{\mathrm{vert}})\right)
   \end{aligned}
    \,.
  $$
  Write $\mathbf{B}^n \mathbb{R}/\mu$ for the pushout
  $$
    \xymatrix{
      \exp_\Delta(\mathfrak{g}) 
        \ar[r]^{\exp_\Delta(\mu)}
        \ar[d]
      &
      \exp_\Delta(b^{n-1}\mathbb{R})
      \ar[r]^>>>>>{\int_{\Delta^\bullet}}_>>>>>\sim
      &
      \mathbf{B}^n \mathbb{R}
      \ar[d]
      \\
      \mathbf{cosk}_n \exp_\Delta(\mathfrak{g})
      \ar[rr]
      &&
      \mathbf{B}^n \mathbb{R}/\mu
    }
    \,.
  $$
  By slight abuse of notation, we shall denote also the bottom morphism by 
  $\exp_\Delta(\mu)$ and refer to it as the \emph{Lie integration of the cocycle} $\mu$.
 \end{definition}
\begin{remark}
  The object $\mathbf{B}^n \mathbb{R}/\mu$ is typically equivalent to the 
  $n$-fold delooping $\mathbf{B}^n (\Lambda_\mu \to \mathbb{R})$ of the reals modulo a lattice 
  $\Lambda_\mu \subset \mathbb{R}$ of periods of $\mu$, as discussed below. 
  Moreover, as discussed in section \ref{ooLieIntegration}, we will be considering weak equivalences 
  $\mathbf{cosk}_n \exp_\Delta(\mathfrak{g}) \stackrel{\sim}{\to} \mathbf{B}G$. 
  Therefore $\exp_\Delta(\mu)$ defines a characteristic morphism of smooth $\infty$-groupoids
  $\mathbf{B}G \to \mathbf{B}^n (\Lambda_\mu \to \mathbb{R})$, presented by the span of 
  morphisms of simplicial presheaves
  $$
    \xymatrix{
       \mathbf{cosk}_n \exp_\Delta(\mathfrak{g}) 
         \ar[r]^>>>>{\exp_\Delta(\mu)}
         \ar[d]^{\wr} 
         & \mathbf{B}^n \mathbb{R}/\mu
       \\
       \mathbf{B}G
    }
    \,.
  $$ 
\end{remark}
\begin{proposition}\label{proposition.periods}
Let $G$ be a Lie group with Lie algebra $\mathfrak{g}$ and  $\mu:\mathfrak{g} \to b^{n-1}\mathbb{R}$ a degree $n$ Lie algebra cocycle. Then there is a smallest subgroup $\Lambda_\mu$
of $(\mathbb{R},+)$ such that we have a commuting diagram
  $$
    \xymatrix{
      \exp_\Delta(\mathfrak{g}) 
        \ar[r]^{\exp_\Delta(\mu)}
        \ar[d]
      &
      \exp_\Delta(b^{n-1}\mathbb{R})
      \ar[r]^>>>>>{\int_{\Delta^\bullet}}_>>>>>\sim
      &
      \mathbf{B}^n \mathbb{R}
      \ar[d]
      \\
      \mathbf{cosk}_n \exp_\Delta(\mathfrak{g})
      \ar[rr]
      &&
      \mathbf{B}^n (\Lambda_\mu \to \mathbb{R})
    }
    \,.
  $$
\end{proposition}
\proof  
We exhibit the commuting diagram naturally over each Cartesian space $U$. 
The vertical map $\mathbf{B}^n \mathbb{R}(U) \to \mathbf{B}^n (\Lambda_\mu \to \mathbb{R})(U)$ is the
obvious quotient map of simplicial abelian groups. 
Since $\mathbf{B}^n (\Lambda_\mu \to \mathbb{R})$ is $(n-1)$-connected and $\mathbf{cosk}_n \exp_\Delta(\mathfrak{g})$ 
is $n$-coskeletal, it is sufficient to define the horizontal map 
$\mathbf{cosk}_n \exp_\Delta(\mathfrak{g}) \to \mathbf{B}^n (\Lambda_\mu \to \mathbb{R})$ on $n$-cells. 
For the diagram to commute, the bottom morphism must send a form 
$A_{\mathrm{vert}} \in \Omega_{\mathrm{si}}^{1}(U \times\Delta^n, \mathfrak{g})$ to the image of 
$\int_{\Delta^n} \mu(A_{\mathrm{vert}}) \in \mathbb{R}$ under the quotient map. For this 
assignment to constitute a morphism of simplicial sets, it must be true that for all 
$A_{\mathrm{vert}} \in \Omega_{\mathrm{si}}^1(U \times \partial \Delta^{n-1}, \mathfrak{g})$ the integral
$\int_{\partial \Delta^{n+1}} A_{\mathrm{vert}} \in \mathbb{R}$ lands in $\Lambda_\mu \subset \mathbb{R}$.

Recall that we may identify flat $\mathfrak{g}$-valued forms on $\partial \Delta^{n+1}$ with based smooth
maps $\partial \Delta^{n+1} \to G$. We observe that  $\int_{\partial \Delta^{n+1}} A_{\mathrm{vert}}$ only
depends on the homotopy class of such a map:  if we have two homotopic $n$-spheres $A_{\mathrm{vert}}$
and $A'_{\mathrm{vert}}$ then by the
arguments as in the proof of proposition \ref{CoskBGAcyclicFibration}, using \cite{wockel}, 
there is a smooth homotopy interpolating between them,
hence a corresponding extension of $\hat A_{\mathrm{vert}}.$  Since this is closed, the fiber integrals
of $A_{\mathrm{vert}}$ and $A'_{\mathrm{vert}}$ coincide.

Therefore we have a group homomorphism
$\int_{\partial \Delta^{n+1}} : \pi_n(G,e_G) \to \mathbb{R}$.
Take $\Lambda_\mu$ to be the subgroup of $\mathbb{R}$ generated by its image. 
This is the minimal subgroup of $\mathbb{R}$ for which we
have a commutative diagram as stated.
\endofproof
\begin{remark}
  If $G$ is compact and simply connected, then its homotopy groups are finitely generated and so is 
$\Lambda_\mu$. 
\end{remark}
\begin{example}
  \label{FirstBMConstruction}
  Let $G$ be a compact, simple and simply connected Lie group and $\mu_3$ the canonical 3-cocycle on
  its semisimple Lie algebra, normalized such that its left-invariant extension
  to a differential 3-form on $G$ represents a generator of $H^3(G,\mathbb{Z}) \simeq \mathbb{Z}$ in 
  de Rham cohomology. In this case we have $\Lambda_{\mu_3} \simeq \mathbb{Z}$ and the 
  diagram of morphisms discussed above is
\[
\xymatrix{
  \exp_\Delta(\mathfrak{g})
   \ar[d]
   \ar[rr]^{\int_{\Delta^\bullet}\exp_\Delta(\mu)}
   &&
   \mathbf{B}^3 \mathbb{R}
   \ar[d]
  \\
  \mathbf{cosk}_3(\exp_\Delta(\mathfrak{g}))
   \ar[d]^{\wr}
   \ar[rr]
  &&
  \mathbf{B}^3(\mathbb{Z} \to \mathbb{R})
   \ar[d]^{\wr}
   \\
  \mathbf{B}G && \mathbf{B}^3U(1)&
}
\]
This presents a morphism of smooth $\infty$-groupoids $\mathbf{B}G \to \mathbf{B}^3 U(1)$.
Let $X \to \mathbf{B}G$ be a morphism of smooth $\infty$-groupoids presented by a 
{\Cech}-cocycle $X \xleftarrow{\sim_{\mathrm{loc}}} \check{C}(\mathcal{U}) \to \mathbf{B}G$ as in 
section \ref{section.Lie_infinity-groupoids}. Then the composite 
$X \to \mathbf{B}G \to \mathbf{B}^3 U(1)$ is a cocycle for a $\mathbf{B}^2 U(1)$-principal 3-bundle
presented by a span of simplicial presheaves
$$
  \xymatrix{
    Q X
    \ar@{->>}[d]^{\wr}
    \ar[r]^>>>>{\hat g}
    &
    \mathbf{cosk}_3 \exp_\Delta(\mathfrak{g})
    \ar[r]^>>>>>{\exp_\Delta(\mu)}
    \ar@{->>}[d]^{\wr}
    &
    \mathbf{B}^3 (\mathbb{Z} \to \mathbb{R})
    \\
    \check{C}(\mathcal{U})
    \ar[d]^{\begin{turn}{270}$\scriptstyle{\sim_\mathrm{loc}}$\end{turn}}
    \ar[r]
    \ar@/^1pc/[u]
    &
    \mathbf{B}G
    \\
    X
  }
  \,.
$$
Here the acyclic fibration $Q X \to \check{C}(\mathcal{U})$ is the pullback of the acyclic fibration
$\mathbf{cosk}_3 \exp_\Delta(\mathfrak{g}) \to \mathbf{B}G$ from proposition \ref{CoskBGAcyclicFibration}
and $\check{C}(\mathcal{U}) \to Q X$ is any choice of section, guaranteed to exist, uniquely up to homotopy, 
since $\check{C}(\mathcal{U})$ is
cofibrant according to proposition \ref{CofibrancyOfGoddCechCover}.

This span composite encodes a morphism of 3-groupoids of {\Cech}-cocycles
$$
  \mathbf{c}_\mu
  :
  \check{C}(\mathcal{U}, \mathbf{B}G)
  \to 
  \check{C}(\mathcal{U}, \mathbf{B}^3(\mathbb{Z} \to \mathbb{R}))
$$
given as follows
\begin{enumerate}
  \item
    it reads in a {\Cech}-cocycle $(g_{i j})$ for a $G$-principal bundle;
  \item
    it forms a lift $\hat g$ of this {\Cech}-cocycle of the following form:
\begin{itemize}
  \item
    over double intersections we have that 
    $\hat g_{i j} : (U_i \cap U_j) \times\Delta^1 \to G$
    is a smooth family of based paths in $G$, with $\hat g_{i j}(1) = g_{i j}$;
  \item
    over triple intersections we have that 
    $\hat g_{i j k} : (U_i \cap U_j \cap U_k) \times\Delta^2 \to G$
    is a smooth family of 2-simplices in $G$ with boundaries labeled by the based paths on double
    overlaps:
    $$
      \xymatrix{
        & g_{i j} \ar[dr]^{g_{i j} \cdot \hat g_{j k}}
        \\
        e 
          \ar@{-}[rr]_{\hat g_{i k}}^{{\ }}="t" 
          \ar@{-}[ur]^{\hat g_{i j}}
          && 
         g_{i k}
        \ar@{}^{\hat g_{i j k }} "t"+(3,2); "t"+(3,2)  
      }
    $$
   \item
     on quadruple intersections we have that 
     $\hat g_{i j k l} : (U_i \cap U_j \cap U_k \cap U_l) : \Delta^3 \to G$ is 
     a smooth family of 3-simplices in $G$, cobounding the union of the 2-simplices corresponding to the triple intersections. 
\end{itemize}
\item
  The morphism $\exp_\Delta(\mu) : \mathbf{cosk}_3 \exp_\Delta(\mathfrak{g}) \to 
  \exp_\Delta(b^2 \mathbb{R})$ takes these smooth families of 3-simplices and integrates
 over them  the 3-form  $\mu_3(\theta \wedge \theta \wedge \theta)$ to obtain the {\Cech}-cocycle
  $$
    (
    \int_{\Delta^3} \hat g_{i j k l}^* \; \mu(\theta \wedge \theta \wedge \theta) 
    \;\; \mod \mathbb{Z})
    \in
    \check{C}(\mathcal{U}, \mathbf{B}^3 U(1))\,.
  $$
 Note that $\mu_3(\theta \wedge \theta \wedge \theta)$ is the canonical 3-form
 representative of a generator of $H^3(G,\mathbb{Z})$.
 \end{enumerate}
  In total the composite of spans therefore encodes a map that takes a {\Cech}-cocycle
  $(g_{i j})$ for a $G$-principal bundle to a degree 3 {\Cech}-cocycle with values in 
  $U(1)$. 
\end{example}
\begin{remark}
  The map of {\Cech} cocycles obtained in the above example from a composite of spans
  of simplicial presheaves is seen to be the special case of the construction considered in
  \cite{brylinski-mclaughlin} that is discussed in section 4 there, where an explicit {\Cech} cocycle for the second Chern class of a principal $SU(n)$-bundle is described. See \cite{bm:pont} for the analogous treatment of the first Pontryagin class of a principal $SO(n)$-bundle and also \cite{bm:geom4-I,bm:geom4-II}.
\end{remark}
\begin{proposition}
  \label{SmoothFirstPontryaginClass}
  For $G = \mathrm{Spin}$ the morphism $\mathbf{c}_{\mu_3}$ from
  example \ref{FirstBMConstruction} 
  is a smooth refinement of the first fractional Pontryagin class
  $$
    \exp_\Delta(\mu_3) = \frac{1}{2}\mathbf{p}_1 : \mathbf{B} \mathrm{Spin} 
     \stackrel{\sim}{\leftarrow} \mathbf{cosk_3} \exp_\Delta(\mathfrak{g}) \to \mathbf{B}^3 U(1)
  $$
  in that postcomposition with this characteruistic map induces the morphism
  $$
    \frac{1}{2}p_1 : H^1(X,\mathrm{Spin}) \to H^4(X,\mathbb{Z})
    \,.
  $$
\end{proposition}
\proof
  Using the identification from example \ref{FirstBMConstruction} 
  of the composite of spans with the construction in
  \cite{brylinski-mclaughlin} this follows from the main theorem there.
  The strategy there is to refine to a secondary characteristic class
  with values in Deligne cocycles that provide the differential refinement of $H^4(X, \mathbb{Z})$.
  The proof is completed by showing that the curvature 4-form of the refining Deligne cocycle is 
  the correct de Rham image of $\frac{1}{2}p_1$.  
\endofproof
Below we shall rederive this theorem as a special case of the more general $\infty$-Chern-Weil homomorphism. 
We now turn to an example that genuinely lives in higher Lie theory and involves higher principal bundles.
\begin{proposition}
  The canonical projection
  $$
    \mathbf{cosk}_7 \exp_\Delta(\mathfrak{so}_{\mu_3}) \to \mathbf{B}\mathrm{String}
  $$
  is an acyclic fibration in the global model structure.
\end{proposition}
\proof
  The 3-cells in $\mathbf{B}\mathrm{String}$ are pairs consisting of 3-cells in $\exp_\Delta(\mathfrak{so})$,
  together with labels on their boundary, subject to a condition that guarantees that 
  the boundary of a 4-cell in $\mathrm{String}$ never wraps a  3-cycle in $\mathrm{Spin}$. Namely,
  a morphism $\partial \Delta^4 \to \mathbf{B}\mathrm{String}$ is naturally identified with 
  a smooth map $\phi : S^3 \to \mathrm{Spin}$ equipped with a 2-form $\omega \in \Omega^2(S^3)$ 
  such that $d \omega = \phi^* \mu_3(\theta \wedge \theta \wedge \theta)$. But since 
  $\mu_3(\theta \wedge \theta \wedge \theta)$ is the image in de Rham cohomology of the generator
  of $H^3(\mathrm{Spin}, \mathbb{Z}) \simeq \mathbb{Z}$ this means that such $\phi$ must 
  represent the trivial element in $\pi_3(\mathrm{Spin})$.
  
  Using this, the proof of the claim proceeds verbatim as that of proposition \ref{CoskBGAcyclicFibration},
  using that the next non-vanishing homotopy group of $\mathrm{Spin}$ after $\pi_3$ is $\pi_7$ and
  that the generator of $H^8(B \mathrm{String}, \mathbb{Z})$ is $\frac{1}{6}p_2$.
\endofproof
\begin{remark}
  Therefore the Lie integration of the 7-cocycle
  $$
    \xymatrix{
      \mathbf{cosk}_7 \exp_\Delta(\mathfrak{so}_{\mu_3})
      \ar[rr]^{\exp_\Delta(\mu_7)}
      \ar[d]^{\wr}
      &&
      \mathbf{B}^7 (\mathbb{Z} \to \mathbb{R})
      \\
      \mathbf{B}\mathrm{String}
    }
  $$
  presents a characteristic map $\mathbf{B}\mathrm{String} \to \mathbf{B}^7 U(1)$.
\end{remark}
\begin{proposition} \label{proposition.integration_of_mu7}
The Lie integration of $\mu_7 : \mathfrak{so}_{\mu_3} \to b^6\mathbb{R}$
is a smooth refinement 
  $$
    \frac{1}{6}\mathbf{p}_2 : \mathbf{B}\mathrm{String} \to \mathbf{B}^7 U(1)
    \,.
  $$
of the second fractional Pontryagin class \cite{SSSII} in that postcomposition with this morphism 
represents the top horizontal morphism in
$$
  \xymatrix{
    H^1(X,\mathrm{String}) \ar[r]^{\frac{1}{6}p_2} 
    \ar[d]
      & H^8(X, \mathbb{Z})\ar[d]^{\cdot 6}
    \\
    H^1(X,\mathrm{Spin}) \ar[r]^{p_2} & H^8(X, \mathbb{Z})
  }
  \,.
$$
\end{proposition}
\proof
  As in the above case, we shall prove this below by refining to a morphism of
  differential cocycles and showing that the corresponding curvature 8-form represents
  the fractional Pontryagin class in de Rham cohomology.
\endofproof

\subsection{Differential characteristic maps by $\infty$-Lie integration}
\label{DiffCharMapsByLieIntegration}

We wish to lift the integration in section \ref{section.characteristic_maps_by_Lie_integration} 
of Lie $\infty$-algebra cocycles from $\exp_\Delta(\mathfrak{g})$
to its differential refinement $\exp_\Delta(\mathfrak{g})_{\mathrm{diff}}$
in order to obtain differential characteristic maps with coefficients in differential cocycles
such that postcomposition with these is the $\infty$-Chern-Weil homomorphism. 
We had obtained $\exp_\Delta(\mu)$ 
essentially by postcomposition of the $k$-cells in $\exp_\Delta(\mathfrak{g})$ with the cocycle 
$\mathfrak{g} \stackrel{\mu}{\to} b^{n-1}\mathbb{R}$. Since the $k$-cells in 
$\exp_\Delta(\mathfrak{g})_{\mathrm{diff}}$ are diagrams, we need to extend the morphism $\mu$ accordingly
to a diagram. We had discussed in section \ref{section.Lie_infinity-algebroids} how transgressive cocycles
extend to a diagram
\[
  \xymatrix{
    \mathrm{CE}(\mathfrak{g}) 
    & \mathrm{CE}(b^{n-1}\mathbb{R})
    \ar[l]_{\mu}
   \\
   \mathrm{W}(\mathfrak{g})
    \ar[u] 
    & 
   \mathrm{W}(b^{n-1}\mathbb{R})\ar[l]_{\mathrm{cs}}\ar[u]
   \\
   \mathrm{inv}(\mathfrak{g})\ar[u] 
    & \mathrm{inv}(b^{n-1}\mathbb{R})
    \ar[u]
   \ar[l]_{\langle -\rangle}
  }
  \,,
\]
where $\langle -  \rangle$ is an invariant polynomial in transgression with $\mu$ and $\mathrm{cs}$
is a Chern-Simons element witnessing that transgression.

\begin{definition} \label{DifferentialLieIntegration}
  Define the morphism of simplicial presheaves
  $$
    \exp_\Delta(\mathrm{cs})_{\mathrm{diff}}:\exp_\Delta(\mathfrak{g})_{\mathrm{diff}}
    \to 
    \exp_\Delta(b^{n-1}\mathbb{R})_{\mathrm{diff}}
  $$
   degreewise by pasting composition with this diagram:
  \[
   \begin{aligned}
    \exp_\Delta(cs)_k 
     &: 
    \left(
      \raisebox{24pt}{
      \xymatrix{
        \Omega_{\mathrm{si}}^\bullet(U\times\Delta^k)_{\mathrm{vert}}
        & 
         \mathrm{CE}(\mathfrak{g})
          \ar[l]_{A_{\mathrm{vert}}} 
        \\
          \Omega_{\mathrm{si}}^\bullet(U\times\Delta^k)
          \ar[u]
          &
          \mathrm{W}(\mathfrak{g})
          \ar[u]\ar[l]_{A} 
     }}
    \right)
    \\
    & \mapsto
    \left(
    \raisebox{24pt}{
    \xymatrix{
      \Omega_{\mathrm{si}}^\bullet(U\times\Delta^k)_{\mathrm{vert}}
      & 
      \mathrm{CE}(\mathfrak{g})
        \ar[l]_{A_{\mathrm{vert}}} 
        & 
        \mathrm{CE}(b^{n-1}\mathbb{R})\ar[l]_{\mu}
        & : \mu(A_{\mathrm{vert}})
       \\
        \Omega_{\mathrm{si}}^\bullet(U\times\Delta^k)
        \ar[u]
        &
        \mathrm{W}(\mathfrak{g})
        \ar[u]\ar[l]_{A} 
        & 
        \mathrm{W}(b^{n-1}\mathbb{R})
        \ar[l]_{\mathrm{cs}}\ar[u]
        &
        : \mathrm{cs}(A)
}}
  \right)
  \end{aligned}
\]
Write $(\mathbf{B}^n \mathbb{R}/\mathrm{cs})_{\mathrm{diff}}$ for the pushout
  $$
    \xymatrix{
      \exp_\Delta(\mathfrak{g})_{\mathrm{diff}}
        \ar[r]^{\exp_\Delta(\mathrm{cs})}
        \ar[d]
      &
      \exp_\Delta(b^{n-1}\mathbb{R})_{\mathrm{diff}}
      \ar[r]^>>>>>{\int_{\Delta^\bullet}}_>>>>>\sim
      &
      \mathbf{B}^n \mathbb{R}_{\mathrm{diff}}
      \ar[d]
      \\
      \mathbf{cosk}_n \exp_\Delta(\mathfrak{g})_{\mathrm{diff}}
      \ar[rr]
      &&
      (\mathbf{B}^n \mathbb{R}/\mathrm{cs})_{\mathrm{diff}}
    }
    \,.
  $$
\end{definition}
\begin{remark}
  This induces a corresponding morphism on the Chern-Weil subobjects
  $$
    \exp_\Delta(\mathrm{cs})_{\mathrm{CW}}:\exp_\Delta(\mathfrak{g})_{\mathrm{CW}}
    \to 
    \exp_\Delta(b^{n-1}\mathbb{R})_{\mathrm{CW}}
  $$
   degreewise by pasting composition with the full transgression diagram
  \[
   \begin{aligned}
    \exp_\Delta(cs)_k 
     &: 
    \left(
      \raisebox{38pt}{
      \xymatrix{
        \Omega_{\mathrm{si}}^\bullet(U\times\Delta^k)_{\mathrm{vert}}
        & 
         \mathrm{CE}(\mathfrak{g})
          \ar[l]_{A_{\mathrm{vert}}} 
        \\
          \Omega_{\mathrm{si}}^\bullet(U\times\Delta^k)
          \ar[u]
          &
          \mathrm{W}(\mathfrak{g})
          \ar[u]\ar[l]_{A} 
         \\
          \Omega^\bullet(U)
           \ar[u]
           & 
           \mathrm{inv}(\mathfrak{g})
          \ar[l]_{F_A}\ar[u] 
     }}
    \right)
    \\
    & \mapsto
    \left(
    \raisebox{38pt}{
    \xymatrix{
      \Omega_{\mathrm{si}}^\bullet(U\times\Delta^k)_{\mathrm{vert}}
      & 
      \mathrm{CE}(\mathfrak{g})
        \ar[l]_{A_{\mathrm{vert}}} 
        & 
        \mathrm{CE}(b^{n-1}\mathbb{R})\ar[l]_{\mu}
        & : \mu(A_{\mathrm{vert}})
       \\
        \Omega_{\mathrm{si}}^\bullet(U\times\Delta^k)
        \ar[u]
        &
        \mathrm{W}(\mathfrak{g})
        \ar[u]\ar[l]_{A} 
        & 
        \mathrm{W}(b^{n-1}\mathbb{R})
        \ar[l]_{\mathrm{cs}}\ar[u]
        &
        : \mathrm{cs}(A)
       \\
      \Omega^\bullet(U)
       \ar[u]
       &
       \mathrm{inv}(\mathfrak{g})
        \ar[l]_{F_A}
        \ar[u] 
        & 
       \mathrm{inv}(b^{n-1}\mathbb{R})
        \ar[u]\ar[l]_{\langle -\rangle}
       &
       : \langle F_A \rangle
}}
  \right)
  \end{aligned}
\]  
Moreover, this restricts further to a morphism of the genuine $\infty$-connection subobjects
$$
  \exp_\Delta(\mathrm{cs})_{\mathrm{conn}}:\exp_\Delta(\mathfrak{g})_{\mathrm{conn}}
  \to 
  \exp_\Delta(b^{n-1}\mathbb{R})_{\mathrm{conn}}
  \,.
$$
Indeed, the commutativity of the lower part of the diagram encodes the classical equation
\[
d \mathrm{cs}(A)=\langle F_A \rangle
\]
stating that the curvature of the connection $\mathrm{cs}(A)$ is the horizontal differential form $\langle F_A\rangle$ in $\Omega(U)$. This shows that the image of   $\exp_\Delta(\mathrm{cs})_{\mathrm{CW}}$ is actually contained in $\exp_\Delta(b^{n-1}\mathbb{R})_{\mathrm{conn}}$, and so the restriction to $\exp_\Delta(\mathfrak{g})_{\mathrm{conn}}$ defines a morphism between the genuine $\infty$-connection subobjects. 
\end{remark}
\begin{remark}
  In the typical application -- see the examples discussed below -- we have that
  $(\mathbf{B}^n \mathbb{R}/\mathrm{cs})_{\mathrm{diff}}$ is 
  $\mathbf{B}^n (\Lambda_\mu \to \mathbb{R})_{\mathrm{diff}}$ and usually
  even $\mathbf{B}^n (\mathbb{Z} \to \mathbb{R})_{\mathrm{diff}}$. The above 
  constructions then yield a sequence of spans in $[\mathrm{CartSp}^{\mathrm{op}}, \mathrm{sSet}]$: 
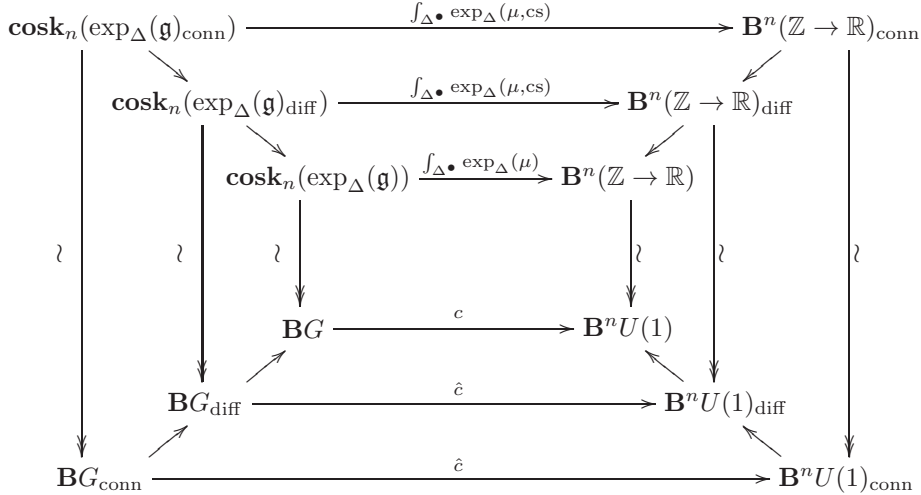
\[
\begin{xy}
,(-45,30)*{\mathbf{cosk}_n(\exp_\Delta(\mathfrak{g})_{\mathrm{conn}})\phantom{i}};(48,30)*{\phantom{i}\mathbf{B}^n(\mathbb{Z}\to\mathbb{R})_{\mathrm{conn}}}**\dir{-}?>*\dir{>}
,(-32,20)*{\mathbf{cosk}_n(\exp_\Delta(\mathfrak{g})_{\mathrm{diff}})\phantom{i}};(32,20)*{\phantom{i}\mathbf{B}^n(\mathbb{Z}\to\mathbb{R})_{\mathrm{diff}}}**\dir{-}?>*\dir{>}
,(-19,10)*{\mathbf{cosk}_n(\exp_\Delta(\mathfrak{g}))\phantom{i}};(21,10)*{\phantom{i}\mathbf{B}^n(\mathbb{Z}\to\mathbb{R})}**\dir{-}?>*\dir{>}
,(-48,-30)*{\mathbf{B}G_{\mathrm{conn}}\phantom{i}};(50,-30)*{\phantom{i}\mathbf{B}^nU(1)_{\mathrm{conn}}}**\dir{-}?>*\dir{>}
,(-34,-20)*{\mathbf{B}G_{\mathrm{diff}}\phantom{i}};(34,-20)*{\phantom{i}\mathbf{B}^nU(1)_{\mathrm{diff}}}**\dir{-}?>*\dir{>}
,(-21,-10)*{\mathbf{B}G\phantom{i}};(21,-10)*{\phantom{i}\mathbf{B}^nU(1)}**\dir{-}?>*\dir{>}
,(-42,27);(-37,23)**\dir{-}?>*\dir{>}
,(42,27);(37,23)**\dir{-}?>*\dir{>}
,(-42,-27);(-37,-23)**\dir{-}?>*\dir{>}
,(42,-27);(37,-23)**\dir{-}?>*\dir{>}
,(-29,17);(-24,13)**\dir{-}?>*\dir{>}
,(29,17);(24,13)**\dir{-}?>*\dir{>}
,(-29,-17);(-24,-13)**\dir{-}?>*\dir{>}
,(29,-17);(24,-13)**\dir{-}?>*\dir{>}
,(-51,27);(-51,-27)**\dir{-}?>*\dir{>>}
,(51,27);(51,-27)**\dir{-}?>*\dir{>>}
,(-35,17);(-35,-17)**\dir{-}?>*\dir{>>}
,(33,17);(33,-17)**\dir{-}?>*\dir{>>}
,(-22,7);(-22,-7)**\dir{-}?>*\dir{>>}
,(22,7);(22,-7)**\dir{-}?>*\dir{>>}
,(-1,-8)*{\scriptstyle{c}}
,(-1,-18)*{\scriptstyle{\hat{c}}}
,(-1,-28)*{\scriptstyle{\hat{c}}}
,(2,12)*{\scriptstyle{\int_{\Delta^\bullet}\exp_\Delta(\mu)}}
,(2,22)*{\scriptstyle{\int_{\Delta^\bullet}\exp_\Delta(\mu,\mathrm{cs})}}
,(2,32)*{\scriptstyle{\int_{\Delta^\bullet}\exp_\Delta(\mu,\mathrm{cs})}}
,(-54,0)*{\wr},(-50,0)
,(-38,0)*{\wr},(-34,0)
,(-25,0)*{\wr},(-21,0)
,(52,0)*{\wr}
,(34,0)*{\wr}
,(23,0)*{\wr}
\end{xy}
\]
Here we have
\begin{itemize}
  \item
    the innermost diagram presents the morphism of smooth $\infty$-groupoids 
    $\mathbf{c}_\mu : \mathbf{B}G \to \mathbf{B}^n U(1)$ that is the characteristic map
    obtained by Lie integration from $\mu$. Postcomposition with this is the morphism
    $$
      \mathbf{c}_\mu : \mathbf{H}(X, \mathbf{B}G) \to \mathbf{H}(X, \mathbf{B}^n U(1))
    $$
    that sends $G$-principal $\infty$-bundles to the corresponding circle $n$-bundles.
    In cohomology/on equivalence classes,  this is the ordinary characteristic class
    $$
      c_\mu : H^1(X, G) \to H^{n+1}(X, \mathbb{Z})
      \,.
    $$
  \item
    The middle diagram is the differential refinement of the innermost diagram. By itself this
    is weakly equivalent to the innermost diagram and hence presents the same characteristic map
    $\mathbf{c}_\mu$. But the middle diagram does support also the projection
    $$
      \mathbf{B}G \stackrel{\sim}{\leftarrow}\mathbf{B}G_\mathrm{diff}
        \to \mathbf{B}^n (\mathbb{Z} \to \mathbb{R})_{\mathrm{diff}} \to 
        \mathbf{\flat}_{\mathrm{dR}}\mathbf{B}^{n+1} \mathbb{R}
    $$
    onto the curvature characteristic classes. This is the simple version of the 
    $\infty$-Chern-Weil homomorphism that takes values in de Rham cohomology
    $H_{dR}^{n+1}(X) = \pi_0 \mathbf{H}(X,\mathbf{\flat}_{\mathrm{dR}}\mathbf{B}^{n+1} \mathbb{R})$
    $$
      \mathbf{H}(X, \mathbf{B}G) \to \mathbf{H}(X,\mathbf{\flat}_{\mathrm{dR}}\mathbf{B}^{n+1} \mathbb{R})
      \,.
    $$
   \item
     The outermost diagram restrict the innermost diagram to differential refinements that
     are genuine $\infty$-connections. These map to genuine $\infty$-connections on circle $n$-bundles
     and hence support the map to secondary characteristic classes
     $$
       \mathbf{H}(X, \mathbf{B}G_{\mathrm{conn}}) \to \mathbf{H}(X,\mathbf{B}^nU(1)_\mathrm{conn})
       \,.
     $$
\end{itemize}
\end{remark}

\subsection{Examples} 

We spell out two classes of examples of the construction of the $\infty$-Chern-Weil homomorphism:

\paragraph{The Chern-Simons circle 3-bundle with connection}

In example \ref{FirstBMConstruction} we had considered the canonical 3-cocycle $\mu_3 \in \mathrm{CE}(\mathfrak{g})$ 
on the semisimple
Lie algebra $\mathfrak{g}$ of a compact, simple and simply connected Lie group $G$ and discussed
how its Lie integration produces a map from {\Cech}-cocycles for $G$-principal bundles to {\Cech}-cocycles
for circle 3-bundles. This map turned out to coincide with that given in
\cite{brylinski-mclaughlin}. We now consider its differential refinement.

From example \ref{ExamplesCSElements} we have a Chern-Simons element $\mathrm{cs}_3$ for $\mu_3$
whose invariant polynomial is the Killing form $\langle -,-\rangle$ on $\mathfrak{g}$.
By definition \ref{DifferentialLieIntegration} this induces a differential Lie integration
$\exp_\Delta(\mathrm{cs})$ of $\mu$.

As a consequence of all the discussion so far, we now simply read off the following corollary.
\begin{corollary}
  \label{TheSmoothFirstFractionalPontryagin}
  Let
  $$
    \xymatrix{
      \mathbf{cosk}_3 \exp_\Delta(\mathfrak{g})_{\mathrm{conn}}
      \ar[rr]^{\int_{\Delta^\bullet}\exp_\Delta(\mathrm{cs}_3)}
      \ar@{->>}[d]^{\wr}
      &&
      \mathbf{B}^3(\mathbb{Z}\to \mathbb{R})_{\mathrm{conn}}
      \\
      \mathbf{B}G_{\mathrm{conn}}
    }
  $$
  be the span of simplicial presheaves obtained from the Lie integration of the differential refinement of the cocycle
  from example \ref{FirstBMConstruction}. Composition with this span 
  $$
    \xymatrix{
      Q X
      \ar[r]^<<<<<<{(\hat g, \hat \nabla)}
      \ar@{->>}[d]^{\wr}
      &
      \mathbf{cosk}_3 \exp_\Delta(\mathfrak{g})_{\mathrm{conn}}
      \ar[rr]^{\int_{\Delta^\bullet}\exp_\Delta(\mathrm{cs}_3)}
      \ar@{->>}[d]^{\wr}
      &&
      \mathbf{B}^3(\mathbb{Z}\to \mathbb{R})_{\mathrm{conn}}
      \\
      \check{C}(\mathcal{U})
       \ar[r]^{(g,\nabla)} 
       \ar[d]^{\begin{turn}{270}$\scriptstyle{\sim_\mathrm{loc}}$\end{turn}}
       \ar@/^1pc/[u]
       & \mathbf{B}G_{\mathrm{conn}} 
      \\
      X
    }
  $$  
  (where $Q X \to \check{C}(\mathcal{U})$ is the pullback acyclic fibration
  and $\check{C}(\mathcal{U}) \to Q X$ any choice of section from the cofibrant
  $\check{C}(\mathcal{U})$ through this acyclic fibration)
  produces a map from {\Cech} cocycles for smooth $G$-principal bundles with connection 
  to  degree 4 {\Cech}-Deligne cocycles 
  $$
    {\hat {\mathbf{c}}}_{\mathrm{cs}}
      : 
     \check{C}(\mathcal{U}, \mathbf{B}G_{\mathrm{conn}})
     \to 
     \check{C}(\mathcal{U}, \mathbf{B}^3 U(1)_{\mathrm{conn}})
  $$
  on a paracompact smooth manifold $X$ as follows:
  \begin{itemize}
    \item
      the input is a set of transition functions and local connection data
      $(g_{i j}, A_i)$ on a differentiably good open cover $\{U_i \to X\}$ 
      as in section \ref{OrdinaryConnections};
      
      (notice that there is a $G$-principal bundle $P \to X$ with Ehresmann connection
      1-form $A \in \Omega^1(P, \mathfrak{g})$ and local sections $\{\sigma_i : U_i \to P|_{U_i}\}$
      such that $\sigma_i|_{U_{i j}}=\sigma_j|_{U_{i j}} g_{i j}$
      and $A_i = \sigma_i^* A$)
    \item
      the span composition produces a lift of this data:
      \begin{itemize}
        \item
         on double intersections a smooth family $\hat g_{i j} : (U_i \cap U_j) \times \Delta^1 \to G$
         of based paths in $G$, together with a 1-form
         $A_{i j} := \hat g_{ij}^* A_i \in \Omega^1(U_{i j} \times \Delta^1, \mathfrak{g})$;
        \item
         on triple intersections a smooth family $\hat g_{i j k} : (U_i \cap U_j \cap U_k) \times \Delta^2 \to G$
         of based 2-simplices in $G$, together with a 1-form
         $A_{i j k} := \hat g_{ijk}^* A_i \in \Omega^1(U_{i j k} \times \Delta^1, \mathfrak{g})$;
        \item
         on quadruple intersections a smooth family 
         $\hat g_{i j k l} : (U_i \cap U_j \cap U_k \cap U_l) \times \Delta^3 \to G$
         of based 2-simplices in $G$, together with a 1-form
         $A_{i j k l} := \hat g_{ijkl}^* A_i \in \Omega^1(U_{i j k l} \times \Delta^1, \mathfrak{g})$;
      \end{itemize}
    \item
      this lifted cocycle data is sent to the {\Cech}-Deligne cocycle
\begin{align*}
        (\mathrm{cs}(A_i),& \int_{\Delta^1} \mathrm{cs}(\hat A_{i j}) , \int_{\Delta^2} \mathrm{cs}(\hat A_{i j k}),
        \int_{\Delta^3} \mu(\hat A_{i j k l}))=\\
             & (\mathrm{cs}(A_i), \int_{\Delta^1} \hat g_{i j}^* \mathrm{cs}(A), 
        \int_{\Delta^2} \hat g_{i jk }^* \mathrm{cs}(A), 
        \int_{\Delta^3} \hat g_{i j k l}^* \mu(A) )
        \,,
\end{align*}
      where $\mathrm{cs}(A)$ is the Chern-Simons 3-form obtained by evaluating a $\mathfrak{g}$-valued
      1-form in the chosen Chern-Simons element $\mathrm{cs}$.
  \end{itemize}
\end{corollary}
\proof
  That we obtain  {\Cech}-Deligne data as indicated is
   a straightforward matter of inserting the definitions of the various morphisms.
  That the data indeed satisifies the {\Cech}-cocycle condition follows from the very
  fact that by construction these are components of a morphism 
  $\check{C}(\mathcal{U}) \to \mathbf{B}^3 (\mathbb{Z} \to \mathbb{R})_{\mathrm{conn}}$,
  as discussed in section \ref{section.Lie_infinity-groupoids}. The curvature 4-form of the 
  resulting {\Cech}-Deligne cocycle is (up to a scalar factor) the Pontryagin form $\langle F_A \wedge F_A\rangle$.
  By the general properties of Deligne cohomology this represents in de Rham cohomology the
  integral class in $H^4(X,\mathbb{Z})$ of the cocycle, so that we find that this is a multiple 
  of the class of the $G$-bundle $P \to X$ corresponding to the Killing form invariant polynomial.
  
  In the case that $G = \mathrm{Spin},$ we have that $H^3(G, \mathbb{Z}) \simeq \mathbb{Z}$.
  By proposition \ref{proposition.periods} it follows that the above construction produces
  a generator of this cohomology group: there cannot be a natural number $\geq 2$ by which this
  $\mathbb{R}/\mathbb{Z}$-cocycle is divisible, since that would mean that 
  $\mu_3(\theta\wedge \theta \wedge \theta)$ had a period greater than 1 around the generator of 
  $\pi_3(G)$, which by construction it does not.  
  But this generator is the fractional Pontryagin class
  $\frac{1}{2}p_1$ (see the review in \cite{SSSII} for instance).
\endofproof
\begin{definition}
  \label{DifferentialFirstPontryagin}
  We write
  $$
    \frac{1}{2}\hat {\mathbf{p}}_1 
     : 
    \mathbf{B}\mathrm{Spin}_{\mathrm{conn}}
     \to
    \mathbf{B}^3 U(1)_{\mathrm{conn}}
  $$
  in $\mathbf{H}$ for the morphism of smooth $\infty$-groupoids given by
  the above corollary and call this the 
  \emph{differential first fractional Pontryagin} map.
\end{definition}
\begin{remark}
  The {\Cech}-Deligne cocycles produced by the span composition in the above corollary 
  are again those considered in section 4 of \cite{brylinski-mclaughlin}.  
  We may regard the
  above corollary as explaining the deeper origin of that construction. But the full impact
  of the construction in the above corollary is that it applies more generally in cases where
  standard Chern-Weil theory is not applicable, as discussed in the introduction. We now 
  turn to the first nontrivial example for the $\infty$-Chern-Weil homomorphism beyond the
  traditional Chern-Weil homomorphism.
\end{remark}

\paragraph{The Chern-Simons circle 7-bundle with connection}

Recall from proposition \ref{proposition.integration_of_mu7} the integration of the 7-cocycle
$\mu_7$ on the String 2-group. We find a Chern-Simons element $\mathrm{cs}_7 \in W(\mathfrak{so}_{\mu_3})$
and use this to obtain the differential refinement of this characteristic map.
\begin{corollary}
  Let
  $$
    \xymatrix{
      \mathbf{cosk}_7 \exp_\Delta(\mathfrak{so}_{\mu_3})_{\mathrm{conn}}
      \ar[rr]^{\int_{\Delta^\bullet}\exp_\Delta(\mathrm{cs}_7)}
      \ar@{->>}[d]^{\wr}
      &&
      \mathbf{B}^7(\mathbb{Z}\to \mathbb{R})_{\mathrm{conn}}
      \\
      \mathbf{B}\mathrm{String}_{\mathrm{conn}}
    }
  $$
  be the span of simplicial presheaves obtained from the Lie integration of the differential refinement of the cocycle
  from proposition \ref{proposition.integration_of_mu7}. Composition with this span 
  $$
    \xymatrix{
      Q X
      \ar[r]^<<<<<<{(\hat g, \hat \nabla)}
      \ar@{->>}[d]^{\wr}
      &
      \mathbf{cosk}_7 \exp_\Delta(\mathfrak{so}_{\mu_3})_{\mathrm{conn}}
      \ar[rr]^{\int_{\Delta^\bullet}\exp_\Delta(\mathrm{cs}_7)}
      \ar@{->>}[d]^{\wr}
      &&
      \mathbf{B}^7(\mathbb{Z}\to \mathbb{R})_{\mathrm{conn}}
      \\
      \check{C}(\mathcal{U})
       \ar[r]^{(g,\nabla)} 
       \ar[d]^{\begin{turn}{270}$\scriptstyle{\sim_\mathrm{loc}}$\end{turn}}
       \ar@/^1pc/[u]
       & \mathbf{B}\mathrm{String}_{\mathrm{conn}} 
      \\
      X
    }
  $$  
  (where $Q X \to \check{C}(\mathcal{U})$ is the pullback acyclic fibration
  and $\check{C}(\mathcal{U}) \to Q X$ any choice of section from the cofibrant
  $\check{C}(\mathcal{U})$ through this acyclic fibration)
  produces a map from {\Cech} cocycles for smooth principal $\mathrm{String}$ 2-bundles with connection 
  to  degree 8 {\Cech}-Deligne cocycles 
  $$
    {\hat {\mathbf{c}}}_{\mathrm{cs}_7}
      : 
     \check{C}(\mathcal{U}, \mathbf{B}\mathrm{String}_{\mathrm{conn}})
     \to 
     \check{C}(\mathcal{U}, \mathbf{B}^7 U(1)_{\mathrm{conn}})
  $$
  on a paracompact smooth manifold $X$.
  For $P \to X$ a principal $\mathrm{Spin}$ bundle with $\mathrm{String}$-structure, i.e. with a 
trivialization of $\frac{1}{2}p_1(P)$, the integral part of 
  ${\hat {\mathbf{c}}}_{\mathrm{cs}_7}(P)$ is the second fractional Pontryagin class
  $\frac{1}{6}p_2(P)$.
\end{corollary}
\proof
  As above.
\endofproof

This completes the proof of theorem \ref{theorem.main}.
\begin{definition}
  \label{DiferentialSecondPontryagin}
  We write
  $$
    \frac{1}{6}\hat {\mathbf{p}_2} : \mathbf{B}\mathrm{String}_{\mathrm{conn}} 
    \to \mathbf{B}^7 U(1)_{\mathrm{conn}}
  $$
  in $\mathbf{H}$ for the morphism of smooth $\infty$-groupoids presented by
the above construction, and speak of the 
 \emph{differential second fractional Pontryagin} map. 
\end{definition}
\begin{remark}
  Notice how the fractional differential class $\frac{1}{6}\hat{\mathbf{p}}_2$ comes out as compared
  to the construction in \cite{brylinski-mclaughlin}, where a \v{C}ech cocycle
representing $-2p_2$ is obtained . There, in order to be able to fill the simplices
  in the 7-coskeleton one works with chains in the Stiefel manifold $\mathrm{SO}(n)/\mathrm{SO}(q)$
  and \emph{multiplies} these with the cardinalities of the torsion homology groups in order to ensure that
  they they become chain boundaries that may be filled.
  
  On the other hand, in the construction above the lift to the {\Cech} cocycle of a $\mathrm{String}$ 2-bundle
  ensures that all the simplices of the cocycle in $\mathrm{Spin}(n)$ can  already be filled 
  genuinely, without passing to multiples. 
  Therefore the cocycle constructed here is a fraction of the cocycle constructed there
  by these integer factors.
\end{remark}

\section{Homotopy fibers of Chern-Weil: twisted differential structures}
\label{HomotopyFibers}

Above we have shown how to construct refined secondary characteristic maps as
morphisms of smooth $\infty$-groupoids of differential cocycles. 
This homotopical refinement of secondary characteristic classes gives access
to their \emph{homotopy fibers}. Here we discuss general properties of these
and indicate how the resulting \emph{twisted differential structures} 
have applications in string physics. 

Some of the computations necessary
for the following go beyond the scope of this article and will not be
spelled out. Details on these are in the followup \cite{FSSII}.
See also section 4.2 of \cite{survey}.

In \ref{section.c_structures} below we consider some basic concepts of obstruction theory in order to 
set the scene for the  its differential refinement further below in \ref{section.Differential_c_structures}.
Before we get to that, it may be worthwhile to note the following subtlety.

There are two different roles
played by topological
spaces in the homotopy theory of higher bundles: 
\begin{enumerate}
  \item they serve as a model for
\emph{discrete} $\infty$-groupoids via the standard Quillen equivalence
$$
  \xymatrix{
   \mathrm{Top}
   \ar@<-3pt>[r]_{\mathrm{Sing}}
   \ar@<+3pt>@{<-}[r]^{\vert - \vert} & \mathrm{sSet}
 }
 \simeq \infty \mathrm{Grpd}
 \,,
$$
where the $\infty$-groupoids on the right are ``discrete'' in direct
generalization to the
sense in which a \emph{discrete group} is discrete, 
\item and they also
serve to model actual geometric structure in the sense of ``continuous cohesion'',
that for instance distinguishes a non-discrete topological group from the underlying
discrete group.
\end{enumerate}
Therefore a \emph{topological group} and more generally a
\emph{simplicial topological group}
is a model for something that pairs these two aspects of topological spaces.

To make this precise, let 
 $\mathrm{Top}$
 be a small
category of suitably nice topological spaces and continuous maps
between them, equipped with the standard Grothendieck topology of open
covers. Then we can consider the the $\infty$-topos of
$\infty$-sheaves over $\mathrm{Top}$, presented by
simplicial presheaves over $\mathrm{Top}$, and this is
the context that contains {\em topological $\infty$-groupoids} in
direct analogy to the smooth $\infty$-groupoids that we considered in
the bulk of the article.

$$
  \infty \mathrm{Grpd} \simeq \mathrm{Sh}_{\infty}(*) \simeq (\mathrm{sSet})^{\mathrm{op}}
$$
$$
  \mathrm{Top}\infty \mathrm{Grpd} \simeq \mathrm{Sh}_{\infty}(\mathrm{Top})
     \simeq ([\mathrm{Top}^{\mathrm{op}}, \mathrm{sSet}]_{\mathrm{loc}})^{\mathrm{op}}
$$
$$
  \mathbf{H} := \mathrm{Smooth}\infty \mathrm{Grpd} \simeq \mathrm{Sh}_{\infty}(\mathrm{CartSp})
     \simeq ([\mathrm{CartSp}^{\mathrm{op}}, \mathrm{sSet}]_{\mathrm{loc}})^{\mathrm{op}}
$$
We have geometric realization functors (see 3.2 and 3.3  in \cite{survey})
$$
 \Pi : \mathrm{Top}\infty \mathrm{Grpd} \to \infty \mathrm{Grpd}
$$
and
$$
 \Pi : \mathrm{Smooth}\infty \mathrm{Grpd} \to \infty \mathrm{Grpd}
 \,,
$$
which on objects represented by simplicial topological spaces are given by the 
traditional geometric realization operation.
For $G$ a topological group or topological $\infty$-group, we write $\mathbf{B}G$ for its
delooping in $\mathrm{Top}\infty \mathrm{Grpd}$. Under geometric realization this
becomes the standard classifying space $B G := \Pi(\mathbf{B}G)$,
which,
while naturally presented by a topological space, is really to be regarded as a presentation
for a discrete $\infty$-groupoid.

\subsection{Topological and smooth $\mathbf{c}$-Structures}
\label{section.c_structures}

An important fact about the geometric realization of topological
$\infty$-groupoids is Milnor's theorem \cite{milnor}:
\begin{theorem}
  \label{MilnorTheorem}
 For every connected $\infty$-groupoid (for instance presented by a connected
homotopy type modeled on a topological space) there is a topological
group such that its topological delooping groupoid $\mathbf{B}G$ has a
geometric realization wealy equivalent to it.
\end{theorem}
This has the following simple, but important consequence.
Let $G$ be a topological group and consider some characteristic map 
$c:BG\to K(\mathbb{Z},n+1)$, representing a characteristic class $[c] \in H^{n+1}(B G, \mathbb{Z})$. 
Then consider the homotopy fiber
$$
\xymatrix{
  BG^c \ar[r] \ar[d]& {*}\ar[d]
  \\
  BG\ar[r]^{c} & K(\mathbb{Z},n+1)
}
$$
formed in $\infty \mathrm{Grpd}$. While this homotopy pullback takes place in discrete
$\infty$-groupoids, Milnor's theorem ensures that there is in fact a topological group
$G^c$ such that $B G^c$ is indeed its classifying space.



%
%
  For $X \in \mathrm{Top}_{\mathrm{sm}}$, the set of homotopy classes $[X, B G]$
  is in natural bijection with equivalence classes of $G$-principal topological bundles
  $P \to X$. One says that $P$ {\it has $c$-structure} if it is in the image of 
  $[X, B G^c ] \to [X, B G]$.

\begin{remark}
By the defining universal property of homotopy fibers, the datum of a (equivalence class of a) principal $G^c$-bundle over $X$ 
is equivalent to the datum of a principal $G$-bundle $P$ over $X$ whose characteristic class $[c(P)]$ vanishes. 
\end{remark}
\begin{example}
Classical examples of this construction are $O^{w_1}=SO$ and $U^{c_1}=SU$. Indeed is well known that the structure group of an $O$-bundle can be reduced to $SO$ if and only if its first Stiefel-Withney class vanishes. More precisely, an principal $SO$-bundle can be seen as a principal $O$-bundle with a trivializiation of the associated orientation $\mathbb{Z}/2\mathbb{Z}$-bundle. Similarly, an $SU$-bundle is a $U$ bundles with a trivialization of the associated determinant bundle, and such a trivialization exists if and only if the first Chern class of the given $U(n)$-bundle vanishes.
\par
A more advanced example is the one described in Section \ref{section.infinity_chern_weil_homomorphism}: $\mathrm{Spin}^{\frac{1}{2}p_1}=\mathrm{String}$, i.e., $\mathrm{String}$-bundles are $\mathrm{Spin}$-bundles with a trivialization of the associated 2-gerbe.
\end{example}
For a more refined description of $c$-structures, we need to consider not just the \emph{set}
of equivalence classes of bundles, but the full cocycle $\infty$-groupoids: whose objects
are such bundles, whose morphisms are equivalences between such bundles, whose 2-morphisms
are equivalaneces between such equivalences, and so on. But for this purposes it matters 
whether we form homotopy fibers in \emph{discrete} or in \emph{topological} $\infty$-groupoids.
We shall be interested in homotopy fibers of topological $\infty$-groupoids.
\begin{definition}
  \label{TopologicalcStructures}
  Let $G$ be a simplicial topological group and $\mathbf{c} : \mathbf{B}G \to \mathbf{B}^n U(1)$
  a classifying map. Write $\mathbf{B}G^c$ for the homotopy fiber
  $$
    \xymatrix{
      \mathbf{B}G^c \ar[r] \ar[d] & {*} \ar[d]
      \\
      \mathbf{B}G \ar[r]^{\mathbf{c}} & \mathbf{B}^n U(1)
    }
  $$ 
  of topological $\infty$-groupoids. Then for $X$ a paracompact topological space, we say that the
  $\infty$-groupoid
  $$
    \mathbf{c}\mathrm{Struc}(X) := \mathrm{Top}\infty\mathrm{Grpd}(X, \mathbf{B}G^c)
  $$
  is the {\it $\infty$-groupoid of topological $\mathbf{c}$-structures} on $X$.
  
  Analogously, for $G$ a  smooth $\infty$-group 
  and $\mathbf{c} : \mathbf{B}G \to \mathbf{B}^n U(1)$ a morphism of smooth $\infty$-groupoids
  as in \ref{section.Lie_infinity-groupoids}, we write $\mathbf{B}G^c$ for its homotopy fiber in 
  $\mathbf{H} = \mathrm{Smooth}\infty\mathrm{Grpd}$ and says that
  $$
    \mathbf{c}\mathrm{Struc}(X) := \mathbf{H}(X, \mathbf{B}G^c)
  $$
  is the {\it $\infty$-groupoid of smooth $\mathbf{c}$-structures} on $X$.
\end{definition}
Among the first nontrivial examples for these notions is the following
\begin{definition}
  Let 
  $$
    \frac{1}{2}\mathbf{p}_1 : \mathbf{B}\mathrm{Spin} \to \mathbf{B}^3 U(1)
  $$
  be the smooth refinement of the first fractional Pontryagin class, from 
  corollary \ref{TheSmoothFirstFractionalPontryagin}. We write
  $$
    \mathbf{B}\mathrm{String} := \mathbf{B}\mathrm{Spin}^{\frac{1}{2}\mathbf{p}_1}
  $$
  and call $\mathrm{String}$ the \emph{smooth String 2-group}.
\end{definition}
By prop. \ref{StringByPullback} the smooth 2-groupoid $\mathbf{B}\mathrm{String}$ is
presented by the simplicial presheaf $\mathbf{cosk}_3 \exp(\mathfrak{so}_{\mu_3})$.
\begin{proposition}
  Under geometric realization the delooping of the smooth String 2-group 
  yields the classifying space of the topological string group
  $$
    \Pi \mathbf{B}\mathrm{String} \simeq B \mathrm{String}
    \,.
  $$
  Moreover, in cohomology smooth $\frac{1}{2}\mathbf{p}_1$-structures on a manifold $X$ are
  equivalent to ordinary String-structures, hence $\frac{1}{2}p_1$-structures.
  \end{proposition}
\proof
  The first statement is proven in section 4.2 of \cite{survey}. The second statement
  follows with proposition \ref{EquivalentModelsForStringGroup} from 
  proposition 4.1 in \cite{NikolausWaldorf}.
\endofproof

\subsection{Twisted differential $\mathbf{c}$-Structures}
\label{section.twisted_c_structures}
\label{section.Differential_c_structures}
By the universal property of the homotopy pullback, the {\it $\infty$-groupoid of topological $\mathbf{c}$-structures} on $X$, def. \ref{TopologicalcStructures}, 
can be equivalently described as the homotopy pullback
$$
    \xymatrix{
      {{\mathbf{c}}}\mathrm{Struc}(X)
      \ar[r]
      \ar[d]
      & {*} \ar[d]
      \\
      \mathrm{Top}\infty\mathrm{Grpd}(X, \mathbf{B}G)
      \ar[r]^{{{\mathbf{c}}}}
      &
      \mathrm{Top}\infty\mathrm{Grpd}(X, \mathbf{B}^{n}U(1))
    }
  $$
  of $\infty$-groupoids of cocycles over $X$, where the right vertical morphism picks 
  any cocycle representing the trivial class.
 From this point of view, there is no reason to restrict one's attention to the fiber of 
 \[
  \mathrm{Top}\infty\mathrm{Grpd}(X, \mathbf{B}G)\xrightarrow{\mathbf{c}}\mathrm{Top}\infty\mathrm{Grpd}(X, \mathbf{B}^{n}U(1))
 \]
over the distinguished point in $\mathrm{Top}\infty\mathrm{Grpd}(X, \mathbf{B}^{n}U(1)$ corresponding to the trivial $\mathbf{B}^{n}U(1)$-bundle over $X$. Rather, it is more natural and convenient to look at all homotopy fibers at once, i.e. to consider all possible (isomorphism classes of) $\mathbf{B}^{n}U(1)$-bundles over $X$.
\begin{definition}\label{def.twisted_c-structures}
  \label{TwistedTopologicalStructures}
  For $\mathbf{c} : \mathbf{B}G \to \mathbf{B}^n U(1)$ a characteristic map in 
  either $\mathbf{H} = \mathrm{Top}\infty \mathrm{Grpd}$ or $\mathbf{H} = \mathrm{Smooth}\infty\mathrm{Grpd}$,
  and for $X$ a paracompact topological space or paracompact smooth manifold, respectively, let 
    ${{\mathbf{c}}} \mathrm{Struc}_{\mathrm{tw}}(X)$ be the $\infty$-groupoid defined by the
  homotopy pullback
  $$
    \xymatrix{
      {{\mathbf{c}}}\mathrm{Struc}_{\mathrm{tw}}(X)
      \ar[r]^{\mathrm{tw}}
      \ar[d]_{\chi} 
      & H^{n+1}(X;\mathbb{Z}) \ar[d]
      \\
      \mathbf{H}(X, \mathbf{B}G)
      \ar[r]^{{{\mathbf{c}}}}
      &
     \mathbf{H}(X, \mathbf{B}^{n}U(1))
    }
    \,,
  $$
  where the right vertical morphism from the cohomology set into the cocycle $n$-groupoid
  picks one basepoint in each connected component, i.e., picks a representative $U(1)$-$(n-1)$-gerbe for each degree $n+1$ integral cohomology class.
\par    
  We call ${{\mathbf{c}}}\mathrm{Struc}_{\mathrm{tw}}(X)$ the 
  \emph{$\infty$-groupoid of (topological or smooth) \emph{twisted ${{\mathbf{c}}}$-structures}}.
  For $\tau \in {\mathbf{c}}\mathrm{Struc}_{\mathrm{tw}}(X)$ we say 
  $\mathrm{tw}(\tau) \in H^{n+1}(X;\mathbb{Z})$ is its \emph{twist} and
  $\chi(\tau) \in \mathrm{Top}\infty\mathrm{Grpd}(X, \mathbf{B}G)$ is the
  (topological or smooth) \emph{underlying $G$-principal $\infty$-bundle} of $\tau$, or that $\tau$ is a
  \emph{$\mathrm{tw}(\tau)$-twisted lift of $\chi(\tau)$}.
\par
For $[\omega]\in H^{n+1}(X;\mathbb{Z})$ a cohomology class,  ${{\mathbf{c}}}\mathrm{Struc}_{\mathrm{tw}=[\omega]}(X)$ is the full sub-$\infty$-groupoid of ${{\mathbf{c}}}\mathrm{Struc}_{\mathrm{tw}}(X)$ on those twisted structures with twist $[\omega]$. 
\end{definition}
The following list basic properties of $\mathbf{c}\mathrm{Struc}_{\mathrm{tw}}(X)$ 
that follow directly on general abstract grounds.
\begin{proposition}\label{prop.features.c-structures}
\begin{enumerate}
\item
The definition of ${{\mathbf{c}}}\mathrm{Struc}_{\mathrm{tw}}(X)$ is independent, up to equivalence, of the choice of the right vertical
      morphism. Indeed, all choices of such are (non canonically) equivalent as $\infty$-functors.
\item For $\mathbf{B}G$ a topological $k$-groupoid for $k \leq n-1$, the $\infty$-groupoid ${{\mathbf{c}}}\mathrm{Struc}_{\mathrm{tw}}(X)$ is an $(n-1)$-groupoid.
\item
The following pasting diagram of homotopy pullbacks shows how ${{\mathbf{c}}}\mathrm{Struc}_{\mathrm{tw}=[\omega]}(X)$ can be equivalently seen as the homotopy fiber of $ \mathrm{Top}\infty\mathrm{Grpd}(X, \mathbf{B}G)\xrightarrow{\mathbf{c}}\mathrm{Top}\infty\mathrm{Grpd}(X, \mathbf{B}^{n}U(1))$ over a representative $U(1)$-$(n-1)$-gerbe for the cohomology class $[\omega]$:
$$
    \xymatrix{
      {{\mathbf{c}}}\mathrm{Struc}_{\mathrm{tw}=[\omega]}(X)
      \ar[r]
      \ar[d] 
      & {*}\ar[d]^{[\omega]}\\
      {{\mathbf{c}}}\mathrm{Struc}_{\mathrm{tw}}(X)
      \ar[r]^{\mathrm{tw}}
      \ar[d]_{\chi} 
      & H^{n+1}(X;\mathbb{Z}) \ar[d]
      \\
      \mathrm{Top}\infty\mathrm{Grpd}(X, \mathbf{B}G)
      \ar[r]^{{{\mathbf{c}}}}
      &
     \mathrm{Top}\infty\mathrm{Grpd}(X, \mathbf{B}^{n}U(1)
    }
    \,,
  $$
In particular one has
\[
\mathbf{c}\mathrm{Struc}_{\mathrm{tw}=0}(X)\cong\mathbf{c}\mathrm{Struc}(X).
\]
\end{enumerate}
\end{proposition}
  \label{SmoothTwistedStringAndFivebraneStructures}
We consider the following two examples, being the direct differential refinement of those
of def. \ref{SmoothTwistedStringAndFivebraneStructures}:
\begin{definition}
  For $\frac{1}{2}\mathbf{p}_1 : \mathbf{B} \mathrm{Spin} \to \mathbf{B}^3 U(1)$
  the smooth first fractional Pontryagin class from prop. \ref{SmoothFirstPontryaginClass}, we call
  $$
    \frac{1}{2}\mathbf{p}_1 \mathrm{Struc}_{\mathrm{tw}}(X)
  $$
  the 2-groupoid of \emph{smoth twisted String-structures} on $X$. For
  $\frac{1}{6}\mathbf{p}_2 : \mathbf{B} \mathrm{Spin} \to \mathbf{B}^7 U(1)$
  the smooth second fractional Pontryagin class from prop. \ref{proposition.integration_of_mu7},
  we call
  $$
    \frac{1}{6}\mathbf{p}_2 \mathrm{Struc}_{\mathrm{tw}}(X)
  $$
  the 6-groupoid of \emph{twisted differential Fivebrane-structures} on $X$.
\end{definition}
The terminology here arises from the applications in string theory that originally motivated
these constructions, as described in \cite{SSSII}.

In order to explicity compute simplicial sets modelling $\infty$-groupoids of 
smooth twisted $\mathbf{c}$-structures,
the usual recipe for computing homotopy fibers applies: it is sufficient to present 
the smooth cocycle $\mathbf{c}$ by a fibration of simplicial presheaves
and then form an ordinary pullback
of simplicial presheaves. We shall discuss now how to obtain such fibrations by Lie integration
of factorizations of the $L_\infty$-cocycles $\mu_3 : \mathfrak{so} \to b^2 \mathbb{R}$
and $\mu_7 : \mathfrak{so}_{\mu_3} \to b^6 \mathbb{R}$. These factorizations
at the $L_\infty$-algebra level are due to \cite{SSSIII}. The full proofs that their 
Lie integration produces the desired fibration is due to \cite{FSSII} and can be found in 
section 4.2 of \cite{survey}.

\begin{definition} Let $\mathfrak{string}:=\mathfrak{so}_{\mu_3}$ be the
string Lie 2-algebra from Definition \ref{StringLie2Algebra}, and
  let $(b \mathbb{R} \to \mathfrak{string})$ be the Lie 3-algebra
  defined by the fact that its Chevalley-Eilenberg algebra is that of $\mathfrak{so}$
  with two additional generators, $b$ in degree 2 and $c$ in degree 3, 
  and with the differential extended to these as
  $$
    d_{\mathrm{CE}} b = c - \mu_3
  $$
  $$
    d_{\mathrm{CE}} c = 0
    \,.
  $$
\end{definition}
  There is an evident sequence of morphisms of $L_\infty$-algbras
  $$
    \mathfrak{so} \to (b \mathbb{R} \to \mathfrak{string}) \to b^2 \mathbb{R}
    \,
  $$
  factoring the 3-cocycle $\mu_3:\mathfrak{so}  \to b^2 \mathbb{R}$.
\begin{proposition}
  \label{FactorizationOfexp(mu3)}
  The Lie integration, according to definition \ref{LieIntegrationOfLInfinityAlgebra}, of this sequence 
  of $L_\infty$-algebra morphisms is a factorization
  $$
    \frac{1}{2}{\mathbf{p}}_1 
    : \mathbf{cosk}_3 \exp(\mathfrak{so})
     \stackrel{\sim}{\to}
    \mathbf{cosk}_3 \exp(b \mathbb{R} \to \mathfrak{string})
      \xymatrix{\ar@{->>} [r] & }
    \mathbf{B}^3 U(1)
  $$
  of the smooth refinement of the first fractional Pontryagin class
  from proposition \ref{SmoothFirstPontryaginClass} into 
  a weak equivalence followed by a fibration in $[\mathrm{CartSp}^{\mathrm{op}}, \mathrm{sSet}]_{\mathrm{proj}}$.
\end{proposition}
\begin{corollary}
  The 2-groupoid of 
  twisted string structures on a smooth manifold $X$
  is
  presented by the ordinary fibers of
  $$
     [\mathrm{CartSp}^{\mathrm{op}}, \mathrm{sSet}](\check{C}(\mathcal{U}), 
        \mathbf{cosk}_3 \exp(b \mathbb{R} \to \mathfrak{so}_{\mu_3}))
     \to 
     [\mathrm{CartSp}^{\mathrm{op}}, \mathrm{sSet}](\check{C}(\mathcal{U})), \mathbf{B}^3 U(1))
     \,.
  $$
\end{corollary}
We spell out the explicit presentation for $\frac{1}{2}\mathbf{p}_1 \mathrm{Struc}_{\mathrm{tw}}(X)$
further below, after passing to the following differential refinement.

Recall that when an $L_\infty$-algebra cocycle $\mu : \mathfrak{g} \to b^n \mathbb{R}$ 
can be transgressed to an invariant polynomial by a Chern-Simons element, as in section 
\ref{DiffCharMapsByLieIntegration}, then the smooth characteristic map 
$\mathbf{c} = \exp(\mu)$ refines to a \emph{differential characteristic map}
  $$
    \hat {\mathbf{c}} : \mathbf{B}G_{\mathrm{conn}} \to \mathbf{B}^n U(1)_{\mathrm{conn}}
    \,,
  $$
  where
  $$
    \mathbf{B}G_{\mathrm{conn}}:= \mathbf{cosk}_{n+1}\exp_\Delta(\mathfrak{g})_{\mathrm{conn}}. 
  $$
In terms of this there is a straightforward refinement of \ref{def.twisted_c-structures}:
\begin{definition}
  \label{TwistedDifferentialStructures}
  For $X$ a smooth manifold,  let ${\hat {\mathbf{c}}} \mathrm{Struc}_{\mathrm{tw}}(X)$ be the $\infty$-groupoid defined by the
  homotopy pullback
  $$
    \xymatrix{
      {\hat {\mathbf{c}}}\mathrm{Struc}_{\mathrm{tw}}(X)
      \ar[r]^{\mathrm{tw}}
      \ar[d]_{\chi} 
      & \hat{H}^{n+1}_{\mathrm{diff}}(X;\mathbb{Z}) \ar[d]
      \\
     \mathbf{H}(X, \mathbf{B}G_{\mathrm{conn}})
      \ar[r]^{{\hat {\mathbf{c}}}}
      &
      \mathbf{H}(X, \mathbf{B}^{n}U(1)_{\mathrm{conn}})
    }
    \,,
  $$
  where the right vertical morphism from the cohomology set into the cocycle $n$-groupoid
  picks one basepoint in each connected component.
  
	We call ${\hat {\mathbf{c}}}\mathrm{Struc}_{\mathrm{tw}}(X)$ the $\infty$-groupoid
	of \emph{twisted differential $\hat {\mathbf{c}}$-structures} on $X$.
   \end{definition}
Such twisted differential structures enjoy the analogous properties
listed in prop. \ref{prop.features.c-structures}.  
In particular, also for differential refinements one has a natural interpretation of untwisted 
$\hat{\mathbf{c}}$-structures:  the component of ${\hat {\mathbf{c}}}\mathrm{Struc}(X)$ over the 0-twist
  is the $\infty$-groupoid of $\hat G$-$\infty$-connections
  $$
    {\hat {\mathbf{c}}}\mathrm{Struc}_{\mathrm{tw = 0}}(X)
    \simeq
    \mathrm{Smooth}\infty\mathrm{Grpd}(X, \mathbf{B}\hat G_{\mathrm{conn}})
    \,,
  $$
  where $\mathbf{B}^{n-2}U(1) \to \hat G \to G$ is the extension of $\infty$-groups 
  classified by $\mathbf{c} : \mathbf{B}G \to \mathbf{B}^n U(1)$. This is shown in detail in \cite{FSSII}, see also section 4.2 of \cite{survey}. 

\subsection{Examples}

We consider the following two examples:
\begin{definition}
  For $\frac{1}{2}\hat {\mathbf{p}}_1 : \mathbf{B} \mathrm{Spin}_{\mathrm{conn}} \to 
  \mathbf{B}^3 U(1)_{\mathrm{conn}}$
  the differential first fractional Pontryagin class from definition \ref{DifferentialFirstPontryagin}
  and $\frac{1}{6}\hat {\mathbf{p}}_2 : \mathbf{B} \mathrm{String}_{\mathrm{conn}} \to 
  \mathbf{B}^7 U(1)_{\mathrm{conn}}$ the differential second fractional 
  Pontryagin class from definition \ref{DiferentialSecondPontryagin}, we call
  $$
    \frac{1}{2}\hat {\mathbf{p}}_1 \mathrm{Struc}_{\mathrm{tw}}(X)
  $$
  the 2-groupoid of \emph{twisted differential String-structures} on $X$
  and
  $$
    \frac{1}{6}\hat {\mathbf{p}}_2 \mathrm{Struc}_{\mathrm{tw}}(X)
  $$
  the 6-groupoid of \emph{twisted differential Fivebrane-structures} on $X$.
\end{definition}
We indicate now explicit constructions of these higher groupoids of 
twisted structures. 

\paragraph{\bf Twisted differential String-structures.}
The factorization 
$$
    \frac{1}{2}{\mathbf{p}}_1 
    : \mathbf{cosk}_3 \exp(\mathfrak{so})
     \stackrel{\sim}{\to}
    \mathbf{cosk}_3 \exp(b \mathbb{R} \to \mathfrak{string})
      \xymatrix{\ar@{->>} [r] & }
    \mathbf{B}^3 U(1)
  $$
of the smooth first fractional Pontryagin class from prop. \ref{FactorizationOfexp(mu3)}
has a differential refinement, from which we can 
compute the 2-groupoid of twisted differential string structures by an ordinary
pullback of simplicial sets. This is achieved by factoring the commutative diagram
$$
\xymatrix{
\mathrm{CE}(\mathfrak{so})&\mathrm{CE}(b^2\mathbb{R})\ar[l]_\mu\\
\mathrm{W}(\mathfrak{so})\ar[u]&\mathrm{W}(b^2\mathbb{R})\ar[u]\ar[l]_{\mathrm{cs}}\\
\mathrm{inv}(\mathfrak{so})\ar[u]&\mathrm{inv}(b^2\mathbb{R})\ar[u]\ar[l]_{\langle\,-\,\rangle}
}
$$
as a commutative diagram
$$
\xymatrix{
  \mathrm{CE}(\mathfrak{so})
   &
  \mathrm{CE}(b\mathbb{R}\to  \mathfrak{string})
   \ar[l]_<<<<<{\sim} 
   &
   \mathrm{CE}(b^2\mathbb{R})\ar[l]
   \\
    \mathrm{W}(\mathfrak{so})\ar[u]
    &
    \tilde{\mathrm{W}}(b\mathbb{R}\to  \mathfrak{string})
    \ar[u]\ar[l]_<<<<<{\sim}
    &
    \mathrm{W}(b^2\mathbb{R})
    \ar[u]\ar[l]
    \\
    \mathrm{inv}(\mathfrak{so})\ar[u]
    &
    \mathrm{inv}(b \mathbb{R} \to \mathfrak{string})
    \ar[u]\ar[l]_<<<<<{=} 
    &
    {\mathrm{inv}}(b^2\mathbb{R})
    \ar[u]\ar[l]
}
 \,
$$
as in \cite{SSSIII}. In the above diagram the Weil algebra $\mathrm{W}(b \mathbb{R} \to \mathfrak{string})$ is replaced by the modified Weil algebra $\tilde{\mathrm{W}}(b \mathbb{R} \to \mathfrak{string})$ presented by
\begin{align*}
&d t^a   =  - \frac{1}{2}C^a{}_{b c} t^b \wedge t^c  +r^a
       \\
   & d b   =  c- \mathrm{cs}_3 + h     
       \\
       &dc = g
       \\
       &d r^a  =  - C^a{}_{b c} t^b \wedge r^c
       \\
       &d h = \langle -,-\rangle - g
       \\
       &d g = 0.
\end{align*}
Here $\{t^a\}$ are the coordinates on $\mathfrak{so}$ relative to a basis $\{e_a\}$, 
$C^a{}_{bc}$ are the structure constants of the Lie brackets of $\mathfrak{so}$ with respect to this basis, 
$b$ and $c$ are the additional generators of the Chevalley-Eilenberg algebra $\mathrm{CE}(b\mathbb{R}\to \mathfrak{string})$, the generators  $r^a,h,g$ are the images of of $t^a,b,c$ via the shift isomorphism, and $\mathrm{cs}_3$ is a Chern-Simons element transgressing the cocycle $\mu_3$ to the Killing form $\langle-,-\rangle$.  The modified Weil algebra $\tilde{\mathrm{W}}(b \mathbb{R} \to \mathfrak{string})$ is ismorphic (via a distinguished isomorphism) to the Weil algebra  $\mathrm{W}(b \mathbb{R} \to \mathfrak{string})$ as a dgca, but the isomorphism between the two does not preserves the graded subspaces of polynomials in the shifted generators. In particuar the modified algebra takes care of realizing the horizontal homotopy between $\langle -,-\rangle$ and $g$ as a polynomial in the shifted generators, see the third item in example \ref{ExamplesCSElements}. Since the notion of curvature forms depends on the splitting of the generators of the Weil algebra into shifted and unshifted generators (see Remark \ref{remark.curvature-forms}), the modified Weil algebra will lead to a modified version of $\exp(b \mathbb{R} \to \mathfrak{string})_{\mathrm{conn}}$ which we will denote by $\exp(b \mathbb{R} \to \mathfrak{string})_{\widetilde {\mathrm{conn}}}$. This is a resolution of $\exp(\mathfrak{so})_{\mathrm{conn}}$ that is naturally adapted to the computation of the homotopy fiber of $\frac{1}{2}\mathbf{p}_1$.   As we will show below, it is precisely this resolution that is the relevant one for applications to the Green-Schwarz mechanism.
\begin{proposition}
  \label{FactorizationOfexp(mu3)Conn}
Lie integration of the above diagram of differential $L_\infty$-algebra cocycles
  provides a factorization
  $$
    \frac{1}{2}{\hat {\mathbf{p}}}_1: 
    \mathbf{cosk}_3 \exp(\mathfrak{so})_{\mathrm{conn}}
     \stackrel{\sim}{\to}
    \mathbf{cosk}_3 \exp(b \mathbb{R} \to \mathfrak{string})_{\widetilde {\mathrm{conn}}}
      \xymatrix{\ar@{->>} [r] & }
    \mathbf{B}^3 U(1)_{\mathrm{conn}}
  $$
  of the  differential first fractional Pontryagin class 
  from definition \ref{DifferentialFirstPontryagin} into 
  a weak equivalence followed by a fibration in $[\mathrm{CartSp}^{\mathrm{op}}, \mathrm{sSet}]_{\mathrm{proj}}$.
\end{proposition}
This is due to \cite{FSSII}. Details can be found in section 4.2 of \cite{survey}. 
\begin{corollary}
  \label{DiffFormDataForTwistedDiffStringStruc}
  The 2-groupoid of 
  twisted differential string structures on a smooth manifold $X$ 
  with respect to a differentiably good open cover $\mathcal{U} = \{U_i \to X\}$ is presented by
  the ordinary fibers of the morphism of simplicial sets
  $$
     [\mathrm{CartSp}^{\mathrm{op}}, \mathrm{sSet}](\check{C}(\mathcal{U}), 
        \mathbf{cosk}_3 \exp_\Delta(b \mathbb{R} \to \mathfrak{string})_{\widetilde{\mathrm{conn}}})
     \to 
     [\mathrm{CartSp}^{\mathrm{op}}, \mathrm{sSet}](\check{C}(\mathcal{U}),\mathbf{B}^3 U(1)_{\mathrm{conn}})
     \,.
  $$
  A $k$-simplex for $k \leq 3$ 
in the simplicial set of local differential forms data describing a differential twisted string structure consists, for any $k$-fold intersection $U_I := U_{i_0, \cdots, i_k}$ 
in the cover $\mathcal{U}$, of a triple $(\omega, B,C)_{I}$ of connection data  such the corresponding curvature data 
$(F_\omega,H,\mathcal{G})_I$ are horizontal. Here
  \[
    \omega_I\in\Omega^1_{\mathrm{si}}(U_I\times\Delta^k;\mathfrak{so}),
   \qquad 
   B_I\in \Omega^2_{\mathrm{si}}(U_I\times\Delta^k;\mathbb{R}), 
   \qquad C_I\in \Omega^3_{\mathrm{si}}(U_I\times\Delta^k;\mathbb{R})
 \]
 and
\[
 F_{\omega_I}=d\omega_I+\frac{1}{2}[\omega_I,\omega_I], 
  \qquad
 H_I =d B_I+\mathrm{cs}(\omega_I)-C_I, \qquad
\mathcal{G}_I=d C_I
  \,.
\]
\end{corollary}
\begin{remark}
The curvature forms of a twisted string structure obey the \emph{Bianchi identities}
\[
 dF_{\omega_I}=-[\omega_I,F_{\omega_I}], 
  \qquad
 d H_I=\langle F_{\omega_I}\wedge F_{\omega_I}\rangle-\mathcal{G}_I, \qquad
 d\mathcal{G}_I=0.
\]
\end{remark}
\paragraph{\bf Twisted differential $\mathrm{String}$-structures and the Green-Schwarz mechanism.}
The above is the local differential form data governing what in string theory 
is called the \emph{Green-Schwarz mechanism}. We briefly indicate what this means and
how it is formalized by the notion of twisted differential $\mathrm{String}$-structures
(for background and references on the physics story see for instance \cite{SSSII}).

The standard action functionals of higher dimensional supergravity theories are
generically \emph{anomalous} in that instead of being functions
on the space of field configurations, they are just sections of a line bundle over these spaces.
In order to get a well defined action principle as input for a path-integral quantization
to obtain the corresponding quantum field theories, one needs to prescribe in addition the
data of a \emph{quantum integrand}. This is a choice of trivialization of these line
bundles, together with a choice of flat connection. For this to be possible, the 
line bundle has to be trivializable and flat in the first place. Its failure to be tivializable
-- its Chern class -- is called the \emph{global anomaly}, and its failure to be flat -- its
curvature 2-form -- is called its local anomaly.

But moreover, the line bundle in question is the tensor product of two different line bundles
with connection. One is a Pfaffian line bundle induced from the fermionic degrees of freedom
of the theory, the other is a line bundle induced from the higher form fields of the theory
in the presence of higher \emph{electric and magnetic charge}. The Pfaffian line bundle is fixed by 
the requirement of supersymmetry, but there is freedom in choosing the background higher 
electric and magnetic charge. Choosing these appropriately such as to ensure that 
the tensor product of the two anomaly line bundles produces a flat trivializable line bundle
is called an \emph{anomaly cancellation} by a \emph{Green-Schwarz mechanism}.

Concretely, the higher gauge background field of 10-dimensional heterotic supergravity 
is the Kalb-Ramond field, which in the absence of \emph{fivebrane magnetic charge} is
modeled by a circle 2-bundle (a bundle gerbe) with connection and curvature 3-form 
$H \in \Omega^3(X)$, satisfying the higher \emph{Maxwell equation}
$$
  d H = 0
  \,.
$$

In order to cancel the relevant quantum anomaly it turns out that a magnetic background charge
density is to be added to the system, whose differential form representative is the difference
$j_{\mathrm{mag}} := 
 \langle F_{\nabla_{\mathrm{Spin}}} \wedge F_{\nabla_{\mathrm{Spin}}}\rangle
 -
\langle F_{\nabla_{\mathrm{SU}}} \wedge F_{\nabla_{\mathrm{SU}}} \rangle$
 between the Pontryagin forms of the Spin-tangent bundle and of a given $\mathrm{SU}$-gauge bundle
(here we leave normalization constants implict in the definition of the 
invariant polynomials $\langle-,- \rangle$).
This modifies the above Maxwell equation locally, on a patch $U_i \subseteq X$ to 
$$
  d H_i = 
  \langle F_{\omega_i} \wedge F_{\omega_i}\rangle
   -
  \langle F_{A_i} \wedge F_{A_i} \rangle
   \,.
$$ 
Comparing with proposition \ref{DiffFormDataForTwistedDiffStringStruc}  
we see that, while such $H_i$ is no longer be the local curvature 3-forms of a circle 2-bundle (2-gerbe),
they are that of a  \emph{twisted circle 3-bundle}
--
a {\v C}ech-Deligne 2-cochain that trivializes the difference of the two 
Chern-Simons {\v C}ech-Deligne 3-cocycles 
-- that is part of the data 
of a twisted differential string-structure with $\mathcal{G}_i = \langle F_{A_i} \wedge F_{A_i} \rangle$. Note that the above differential form
equation exhibits a de Rham homotopy between the two Pontryagin forms. This is
the local differential aspect of the very definition of a twisted differential string-structure:
a homotopy from the Chern-Simons circle 3-bundle of the Spin-tangent bundle to a given
twisting circle 3-bundle, which here is itself a Chern-Simons 3-bundle, coming from an 
$\mathrm{SU}$-bundle.

    This anomaly cancellation has been known in the physics literture since the
  seminal article \cite{Killingback}.
      Recently \cite{Bunke} has given a rigorous proof in the special 
      case that underlying topological class of the
   twisting gauge bundle is trivial.  This proof used the 
   model of  twisted differential string structures with topologically tivial twist 
   given in \cite{waldorf}.
   This model is constructed in terms of bundle 2-gerbes and does not
   exhibit the homotopy pullback property of definition \ref{TwistedDifferentialStructures} explicitly. 
   However, the author shows that his model satisfies the
   properties \ref{prop.features.c-structures} satisfied by the abstract homotopy pullback.

\paragraph{\bf Twisted differential fivebrane structures.}

The construction of an explicit Kan complex model for the 6-groupoid of twisted differential
fivebrane structures proceeds in close analogy to the above discussion for twisted differential
string structures, by adding throughout one more layer of generators in the CE-algebra.

\begin{definition}
Write
$$
  \mathfrak{fivebrane} := (\mathfrak{so}_{\mu_3})_{\mu_7}
$$
for the $L_\infty$-algebra extension of the $\mathfrak{string}$ Lie 2-algebra
(def. \ref{StringLie2Algebra}) by the 7-cocycle $\mu_7 : \mathfrak{so}_{\mu_3} \to b^6 \mathbb{R}$
(remark \ref{7CocycleOnStringLie2Algebra}) according to prop. \ref{LooExtension}.
Following \cite{SSSI} we call this the \emph{fivebrane Lie 6-algebra}.
\end{definition}
\begin{remark}
  The Chevelley-Eilenberg algebra $\mathrm{CE}(\mathfrak{fivebrane})$ is
  given by
  $$
    d t^a = - \frac{1}{2}C^a{}_{b c} t^b \wedge t^c
  $$
  $$
    d b_2 = -\mu_3 := - \frac{1}{2}\langle -,[-,-]\rangle
  $$
  $$
    d b_6 = -\mu_7 := - \frac{1}{8}\langle -,[-,-], [-,-], [-,-]\rangle
  $$
  for $\{t^a\}$ and $b_2$ generators of degree 1 and 2, respectively, as for the
  string Lie 2-algebra, and $b_6$ a new generator in degree 6.
\end{remark}
\begin{definition}
  Let $(b^5 \mathbb{R} \to \mathfrak{fivebrane})$ be the Lie 7-algebra 
  defined by having CE-algebra given by
  $$
    d t^a = - \frac{1}{2}C^a{}_{b c} t^b \wedge t^c
  $$
  $$
    d b_2 = c_3 - \mu_3 
  $$
  $$
    d b_6 = c_7 - \mu_7 
  $$  
  $$
    d c_3 = 0
  $$
  $$
    d c_7 = 0
    \,.
  $$
\end{definition}
\begin{proposition}
  In the evident factorization
  $$
    \mu_7 : 
    \xymatrix{
       \mathrm{CE}(\mathfrak{string}) \ar@{<-}[r]^<<<<<<<\sim &  
      \mathrm{CE}(b^5 \to \mathfrak{fivebrane}) \ar@{<-}@{<-}[r] 
      &
      \mathrm{CE}(b^6 \mathbb{R})
    }
  $$
  of the 7-cocycle $\mu_7$, the first morphism is a quasi-isomorphism.
\end{proposition}
As before, it is convenient to lift this factorization to the differential refinement
by using a slightly modified Weil algebra to collect horizontal generators 
$$
  \tilde W(b^5 \to \mathfrak{fivebrane}) \simeq W(b^5 \to \mathfrak{fivebrane})
$$
given by
  $$
    d t^a = - \frac{1}{2}C^a{}_{b c} t^b \wedge t^c + r^a
  $$
  $$
    d b_2 = c_3 - \mathrm{cs}_3  + h_3
  $$
  $$
    d b_6 = c_7 - \mathrm{cs}_7  + h_7
  $$  
  $$
    d c_3 = g_4
  $$
  $$
    d c_7 = g_8
  $$
  $$
    d h_3 = \langle -,- \rangle - g_4
  $$
  $$
    d h_7 = \langle -,-,-,- \rangle - g_8,
  $$
  where $\langle -,-,-,- \rangle$ is the second Pontryagin polynomial for $\mathfrak{so}$, to obtain a factorization
$$
\xymatrix{
  \mathrm{CE}(\mathfrak{string})
   &
  \mathrm{CE}(b^5\mathbb{R}\to  \mathfrak{fivebrane})
   \ar[l]_<<<<<{\sim} 
   &
   \mathrm{CE}(b^6\mathbb{R})\ar[l]
   \\
    \mathrm{W}(\mathfrak{string})\ar[u]
    &
    \tilde{\mathrm{W}}(b^5\mathbb{R}\to  \mathfrak{fivebrane})
    \ar[u]\ar[l]_<<<<<{\sim}
    &
    \mathrm{W}(b^6\mathbb{R})
    \ar[u]\ar[l]
    \\
    \mathrm{inv}(\mathfrak{string})\ar[u]
    &
    \mathrm{inv}(b^5 \mathbb{R} \to \mathfrak{fivebrane})
    \ar[u]\ar[l]_<<<<<{=} 
    &
    {\mathrm{inv}}(b^6\mathbb{R})
    \ar[u]\ar[l]
  }
 \,.
$$
This is the second of the big diagrams in \cite{SSSIII}. Using this and following
through the same steps as for twisted differential string-structures above, one finds that 
the 6-groupoid of twisted differential fivebrane structures over some $X$ with respect to
a diffrentiably good open cover $\mathcal{U}$ has $k$-cells for $k \leq 7$ given by
differential form data 

$$
   \omega_I\in\Omega^1_{\mathrm{si}}(U_I\times\Delta^k;\mathfrak{so}),
   \qquad 
   (B_2)_I\in \Omega^2_{\mathrm{si}}(U_I\times\Delta^k;\mathbb{R}), 
   \qquad 
   (B_6)_I\in \Omega^6_{\mathrm{si}}(U_I\times\Delta^k;\mathbb{R}), 
$$
$$
    (C_3)_I \in \Omega^3_{\mathrm{si}}(U_I\times\Delta^k;\mathbb{R}), 
    \qquad
    (C_7)_I \in \Omega^7_{\mathrm{si}}(U_I\times\Delta^k;\mathbb{R})
$$
with horizontal curvature forms
$$
 F_{\omega_I}=d\omega_I+\frac{1}{2}[\omega_I,\omega_I], 
$$
$$
 (H_3)_I =d (B_2)_I+\mathrm{cs}_3(\omega_I)-(C_3)_I, 
  \qquad
 (H_7)_I =d (B_6)_I+\mathrm{cs}_7(\omega_I)-(C_7)_I, 
$$
$$
  (\mathcal{G}_4)_I=d (C_4)_I , \qquad (\mathcal{G}_8)_I=d (C_8)_I
  \,.
$$
And \emph{Bianchi identities}
$$
 d F_{\omega_I}=-[\omega_I,F_{\omega_I}], 
$$
$$
 d (H_3)_I=\langle F_{\omega_I}\wedge F_{\omega_I}\rangle- (\mathcal{G}_4)_I, 
  \qquad
 d (H_7)_I=\langle F_{\omega_I}\wedge F_{\omega_I}\wedge F_{\omega_I}\wedge F_{\omega_I}\rangle
   -(\mathcal{G}_8)_I, \qquad
$$
$$
 d(\mathcal{G}_4)_I=0, \qquad d(\mathcal{G}_8)_I=0
$$

\paragraph{\bf Twisted differential fivebrane structures and the dual Green-Schwarz mechanism.}
On a 10-dimensional smooth manifold $X$ a (twisted) circle 2-bundle with local connection form $\{(B_2)_I\}$
and (local) curvature forms $\{(H_3)_I\}$ is the electric/magnetic dual of a
(twisted) circle 6-bundle with local connection 6-forms $\{(B_2)_I\}$ and
(local) curvature forms $\{(H_7)_I\}$. It is expected 
(see the references in \cite{SSSII}) that there is a magnetic dual
quantum heterotic string theory where the string -- electrically charged under $B_2$ --
is replaced by the fundamental fivebrane -- magnetically charged under $B_6$. 
While the understanding of the 6-dimensional fivebrane sigma-model is rudimentary,
its fermionic worldvolume quantum anomaly can and has been computed and the corresponding
anomaly cancelling Green-Schwarz mechanism has been written down (all reviewed in \cite{SSSII}).
If $X$ does have differential string structure then its local differential expression is the relation
$$
 d (H_7)_I=\langle F_{\omega_I}\wedge F_{\omega_I}\wedge F_{\omega_I}\wedge F_{\omega_I}\rangle
   -
   \langle F_{A_I}\wedge F_{A_I}\wedge F_{A_I}\wedge F_{A_I}\rangle
$$
for some normalization of invariant polynomials,
where the second term is the curvature characteristic form of the next higher Chern class of 
the background $\mathrm{SU}$-principal gauge bundle. Comparing with the above formula,
we find that this is indeed modeled by twisted differential fivebrane structures.

\newpage

\section*{Appendix: $\infty$-Stacks over the site of Cartesian spaces}
\label{section.Cartesian-spaces}
\addcontentsline{toc}{section}{Appendix: $\infty$-Stacks over the site of Cartesian spaces}

Here we give a formal description of simplicial presheaves over the site of Cartesian spaces 
and prove several statements mentioned in Section \ref{section.Lie_infinity-groupoids}.

\begin{definition} \label{DifferentiablyGoodOpenCover}
  For $X$ a $d$-dimensional paracompact smooth manifold, a \emph{differentiably good open cover} is an open cover 
  $\mathcal{U} = \{U_i \to X\}_{i \in I}$
  such that for all $n \in \mathbb{N}$ every $n$-fold intersection $U_{i_1} \cap \cdots \cap U_{i_n}$
  is either empty or \emph{diffeomorphic} to $\mathbb{R}^d$.
\end{definition}
Notice that this is asking a little more than that the intersections are contractible, as for ordinary
good open covers. 
\begin{propositionapp} 
  Differentiably good open covers always exist.
\end{propositionapp}
\proof
  By \cite{greene} every paracompact manifold admits a Riemannian metric with positive 
  convexity radius $r_{\mathrm{conv}} \in \mathbb{R}$. Choose such a metric and choose an open cover consisting 
  for each point $p\in X$ of the geodesically convex open subset $U_p := B_p(r_{\mathrm{conv}})$ given by the
  geodesic $r_{\mathrm{conv}}$-ball at $p$. Since the injectivity radius of any metric is at least 
  $2r_{\mathrm{conv}}$
  \cite{berger} it follows from the minimality of the geodesics in a geodesically convex region 
  that inside every finite nonempty intersection $U_{p_1} \cap \cdots \cap U_{p_n}$ 
  the geodesic flow around any point $u$ is of radius less than or equal the injectivity radius
  and is therefore a diffeomorphism onto its image. Moreover, the preimage of the intersection region
  under the geometric flow is a star-shaped region in the tangent space $T_u X$:  
  the intersection of geodesically convex regions is itself geodesically convex, so that for
  any $v \in T_u X$ with $\exp(v) \in U_{p_1} \cap \cdots \cap U_{p_n}$ the whole geodesic segment
  $t \mapsto \exp(t v)$ for $t \in [0,1]$ is also in the region.
So we have that every finite non-empty intersection of the $U_p$ is diffeomorphic to a star-shaped
  region in a Euclidean space. It is then a folk theorem that every star-shaped region is diffeomorphic to 
  an $\mathbb{R}^n$; an explicit proof of this fact is in theorem 237 of \cite{ferus}.
\endofproof
\medskip 

\noindent Recall the following notions \cite{johnstone}.
\begin{definitionapp}
  A \emph{coverage} on a small category $\mathcal{C}$ is for each object $U \in \mathcal{C}$ a choice of
  collections of morphisms $\mathcal{U} = \{U_i \to U\}$ -- called \emph{covering families} -- such that 
  whenever $\mathcal{U}$ is a covering family and $V \to U$ any morphism in $\mathcal{C}$ there exists
  a covering family $\{V_j \to V\}$ such that all diagrams
  $$
    \xymatrix{
      V_j \ar[d]\ar[r]^{\exists} & U_i \ar[d]
      \\
      V \ar[r] & U
    }
  $$
  exist as indicated. The \emph{covering sieve} corresponding to a covering family $\mathcal{U}$ is the colimit
  $$
    S(\mathcal{U}) = {\lim_{\to}}_{[k] \in \Delta} \check{C}(\mathcal{U})_k \in [\mathcal{C}^{\mathrm{op}},\mathrm{Set}]
  $$
  of the {\Cech}-nerve, formed after Yoneda embedding in the category of presheaves on $\mathcal{C}$.
 \end{definitionapp}
 \begin{definitionapp} 
  A \emph{site} is a small category $\mathcal{C}$ equipped with a coverage. A \emph{sheaf} on a site
  is a presheaf $A : \mathcal{C}^{\mathrm{op}} \to \mathrm{Set}$ such that for each covering sieve
  $S(\mathcal{U}) \to U$ the morphism
  $$
    A(U) \simeq [\mathcal{C}^{\mathrm{op}}, \mathrm{Set}](U,A) \to [\mathcal{C}^{\mathrm{op}}, \mathrm{Set}](S(\mathcal{U}),A)
  $$ 
  is an isomorphism.
\end{definitionapp}
\begin{remarkapp}
  Often this is formulated in terms of Grothendieck topologies instead of coverages. But
  every coverage induces a unique Gorthendieck topology such that the corresponding notions
  of sheaf coincide. An advantage of using coverages is that there are fewer morphisms to check
  the sheaf condition against. 

  In the language of left exact reflective localizations: the coverage sieve projections of a covering 
  family form a small set such that localizing the presheaf category at this set produces the
  category of sheaves. This localization however inverts more morphisms than just the coverage sieves.
  This saturated class of inverted morphisms contains also the sieve projections of the corresponding
  Grothendieck topology.
   
  Below we use this for obtaining the $\infty$-stacks/$\infty$-sheaves by left Bousfield localization
  just at a coverage. 
\end{remarkapp}
\begin{corollaryapp}
  Differentiably good open covers form a coverage
  on the category $\mathrm{CartSp}$.
\end{corollaryapp}
\proof
  The pullback of a differentiably good open cover always exists in the category of manifolds, where
  it is an open cover. By the above, this may always be refined again by a differentiably good 
  open cover. 
\endofproof
\begin{definitionapp}
We consider $\mathrm{CartSp}$ as a site 
by equipping it with this differentiably-good-open-cover coverage.
\end{definitionapp}
\begin{definitionapp}
Write $[\mathrm{CartSp}^{\mathrm{op}}, \mathrm{sSet}]_{\mathrm{proj}}$ for the global
projective model category structure on simplicial presheaves whose weak equivalences and
fibrations are objectwise those of simplicial sets. Write $[\mathrm{CartSp}^{\mathrm{op}}, \mathrm{sSet}]_{\mathrm{proj},\mathrm{loc}}$ for the  
left Bousfield localization 
of $[\mathrm{CartSp}^{\mathrm{op}}, \mathrm{sSet}]_{\mathrm{proj}}$ at at the set of coverage {\Cech}-nerve projections $\check{C}(\mathcal{U}) \to U$.
This is a simplicial model category with respect to the canonical simplicial enrichment of 
simplicial presheaves, see \cite{dugger}. For $X, A$ two objects, we write $[\mathrm{CartSp}^{\mathrm{op}}, \mathrm{sSet}](X,A) \in \mathrm{sSet}$ for the simplicial hom-complex of morphisms between them.
\end{definitionapp}
\begin{propositionapp} \label{CofibrancyOfGoddCechCover}
  \label{CechResolution}
  In $[\mathrm{CartSp}^{\mathrm{op}}, \mathrm{sSet}]_{\mathrm{proj},\mathrm{loc}}$ the {\Cech}-nerve
  $\check{C}(\mathcal{U}) \to X$ of a differentiably good open cover over a paracompact smooth manifold
  $X$ is a cofibrant resolution of $X$.
\end{propositionapp}
\proof
  By assumption $\check{C}(\mathcal{U})$ is degreewise a coproduct of representables 
  (this is what the definition of \emph{differentiably good open cover} formulates).
  Clearly its degeneracies split off as a direct summand
  in each degree (the summand of intersections $U_{i_0} \cap \cdots U_{i_n}$ where at least one
  index repeats). With this it follows from corollary 9.4
  in \cite{dugger} that $\check{C}(\mathcal{U})$ is cofibrant in the global projective  
  model structure. Since left Bousfield localization keeps the cofibrations unchanged, it follows
  that it is also cofibrant in the local structure. That the projection $\check{C}(U) \to X$ is a 
  weak equivalence in the local structure
  follows by using our theorem A.\ref{NatureOfTheLocalization} below in proposition A.4  
  of \cite{dugger-hollander-isaksen}. 
\endofproof
\begin{corollaryapp} \label{CharacterizationOfLocalFibrancy}
  The fibrant objects of $[\mathrm{CartSp}^{\mathrm{op}}, \mathrm{sSet}]_{\mathrm{proj}, \mathrm{loc}}$
  are precisely those simplicial presheaves $A$ that are objectwise Kan complexes and such that 
  for all differentiably good open covers $\mathcal{U}$ of a Cartesian space $U$ the induced morphism
  $$
    A(U) \stackrel{\simeq}{\to} [\mathrm{CartSp}^{\mathrm{op}}, \mathrm{sSet}](U,A) \to 
   [\mathrm{CartSp}^{\mathrm{op}}, \mathrm{sSet}](\check{C}(\mathcal{U}),A)
  $$ is a weak equivalence of Kan complexes.

This is the \emph{descent condition} or \emph{$\infty$-sheaf/$\infty$-stack condition} on $A$.
\end{corollaryapp}
\proof
  By standard facts about left Bousfield localizations we have that the fibrant objects are
  the degreewise fibrant object such that the morphisms
  $$
    \mathbb{R}\mathrm{Hom}(U,A) \to \mathbb{R}\mathrm{Hom}(\check{C}(\mathcal{U}), A)
  $$
  are weak equivalences of Kan complexes, where $\mathbb{R}\mathrm{Hom}$ denotes the 
  right derived simplicial hom-complex in the global projective model structure. 
  Since every representable $U$ is cofibrant and since
  $\check{C}(\mathcal{U})$ is cofibrant by the above proposition, these hom-complexes are 
  equivalent to the hom-complexes in $[\mathrm{CartSp}^{\mathrm{op} }, \mathrm{sSet}]$ as
  indicated.
\endofproof
Finally we establish the equivalence of the localization at a coverage that we are using to the
localization at the corresponding Grothendieck topology, which is the one commonly found
discussed in the literature. 
\begin{theoremapp} \label{NatureOfTheLocalization}
  Let $\mathcal{C}$ be any small category equipped with a coverage given by covering families $\{U_i \to U\}$.
  
  Then the $\infty$-topos presented by the left Bousfield localization of 
  $[\mathcal{C}^{\mathrm{op}}, \mathrm{CartSp}]_{\mathrm{proj}}$ at the coverage covering families is equivalent
  to that presented by the left Bousfield localization at the covers for the corresponding 
  Grothendieck topology.
\end{theoremapp}
We prove this for 
the \emph{injective} model structure on simplicial presheaves. 
The result then follows since that is Quillen equivalent to 
the projective one and so presents  the same $\infty$-topos.

Write $S(\mathcal{U}) \to j(U)$ for the sieve corresponding to a covering family, 
regarded as a subfunctor of the representable functor $j(U)$ (the Yoneda embedding of $U$), 
which we both regard as simplicially discrete objects in $[\mathcal{C}^{\mathrm{op}}, \mathrm{sSet}]$.
Write $[\mathcal{C}^{\mathrm{op}}, \mathrm{sSet}]_{\mathrm{inj}, \mathrm{cov}}$ for the left 
Bousfield localization of the injective model structure at the morphisms $S(\mathcal{U}) \to j(U)$ 
corresponding to covering families.
\begin{lemmaapp}
A subfunctor inclusion $\tilde S \hookrightarrow j(U)$ corresponding to a sieve that contains a covering sieve $S(\mathcal{U})$ is a weak equivalence in $[\mathcal{C}^{\mathrm{op}}, \mathrm{sSet}]_{\mathrm{inj},\mathrm{cov}}$
\end{lemmaapp}
\proof
Let $J$ be the set of morphisms in the bigger sieve that are not in the smaller sieve.
By assumption we can find for each $j \in J$ a covering family
$\{V_{j,k} \to V_j\}$ such that for all $j,i$ the diagrams
$$
  \xymatrix{
    V_{j,k} \ar[r] \ar[d]& U_{i} \ar[d]
    \\
    V_j \ar[r]^f & U
  }
$$ 
commute. Consider then the commuting diagram
$$
  \xymatrix{
    \coprod_{j} S(\{V_{j,k}\}) \ar@{^{(}->}[r] \ar[d]^{\wr} & S(\{U_i\} \cup \{V_{j,i}\}) \ar[d]
    \\
    \coprod_j j(V_j) \ar[r] & S(\{U_i\} \cup \{V_j\}) = \tilde S
  }
  \,.
$$
Observe that this is a pushout in $[\mathcal{C}^{\mathrm{op}}, \mathrm{sSet}]$, 
that the top morphism is a cofibration 
in $[\mathcal{C}^{\mathrm{op}}, \mathrm{sSet}]_{\mathrm{inj}}$ 
and hence in $[\mathcal{C}^{\mathrm{op}}, \mathrm{sSet}]_{\mathrm{inj},\mathrm{cov}}$, 
that the left morphism a weak equivalence in the local structure 
and that by general properties of left Bousfield localization 
the localization is left proper. 
Therefore the pushout morphism $S(\{U_i\} \cup \{V_{j,k}\}) \to S(\{U_i\} \cup \{V_j\}) = \tilde S$ 
is a weak equivalence.

Then observe that from the horizontal morphisms of the above commuting diagrams that defined the 
covers $\{V_{j,k} \to V_j\}$ we have an induced morphism $S(\{U_i\} \cup \{V_{j,k}\}) \to S(\{U_i\})$
that exhibit $S(\{U_i\})$ as a retract
$$
  \xymatrix{
    S(\{U_i\}) \ar[d]\ar[r]& S(\{U_i\} \cup \{V_{j,k}\}) \ar[d]\ar[r]& S(\{U_i\}) \ar[d]
    \\
    \tilde S \ar[r]^=& \tilde S \ar[r]^=& \tilde S 
  }
  \,.
$$
By closure of weak equivalences under retracts, this shows that the inclusion $S(\{U_i\}) \to \tilde S$ 
is a weak equivalence.  By 2-out-of-3 this finally means that $\tilde S \hookrightarrow j(U)$ is a weak equivalence.
\endofproof
\begin{corollaryapp}
For $S(\{U_i\}) \to j(U)$ a covering sieve, 
its pullback $f^*S(\{U_i\}) \to j(V)$ in $[\mathcal{C}, \mathrm{sSet}]$ along any morphism $j(f) : j(V) \to j(U)$ 
$$
  \xymatrix{
    f^* S(\{U_i\}) \ar[r] \ar[d] & S(\{U_i\}) \ar[d]
    \\
    j(V) \ar[r]^{j(f)}& j(U)
  }
$$
is also a weak equivalence.
\end{corollaryapp} 
\begin{lemmaapp}
If $S(\{U_i\}) \to j(U)$ is the sieve of a covering family and $\tilde S \hookrightarrow j(U)$ is any sieve such that for every $f_i : U_i \to U$ the pullback $f^* \tilde S$ is a weak equivalence, then $\tilde S \to j(U)$ becomes an isomorphism in the homotopy category.
\end{lemmaapp}
\proof
First notice that if $f_i^* \tilde S$ is a weak equivalence for every $i$, 
then the pullback of $\tilde S$ to any element of the sieve $S(\{U_i\})$ is a weak equivalence. 
Use the Yoneda lemma to write
$$
  S(\{U_i\}) \simeq \lim_{\underset{V \to U_i \to U}{\to}} j(V)
  \,.
$$
Then consider these objects in the $\infty$-category of $\infty$-presheaves 
that is presented by $[\mathcal{C}^{\mathrm{op}}, \mathrm{sSet}]_{\mathrm{inj}}$ \cite{lurie}.
Since that has universal colimits we have the pullback square
$$
  \xymatrix{
    i^* \lim\limits_\to j(V)
    \ar[r]^\sim & \lim\limits_{\to} f_V^* \tilde S
    \ar[r] \ar[d] & \tilde S \ar[d]^i
    \\
    S(\{U_i\}) \ar[r]^\sim & 
    \lim\limits_{\underset{f_V : V \to U_i \to U}{\to}} j(V)
    \ar[r]^>>>{(f_V)}&
    j(U)
  }
$$
and the left vertical morphism is a colimit over morphisms that are weak equivalences in 
$[\mathcal{C}^{\mathrm{op}},\mathrm{sSet}]_{\mathrm{inj},\mathrm{loc}}$. By the general properties of 
reflective sub-$\infty$-categories this means that the total left vertical morphism becomes an 
isomorphism in the homotopy category of $[\mathcal{C}^{\mathrm{op}}, \mathrm{sSet}]_{\mathrm{inj},\mathrm{cov}}$. Also the bottom morphism is an isomorphism there, and hence the right vertical one is.
\endofproof
{\it Proof of the theorem}.
The two lemmas show that all morphisms $S(\{V_j\}) \to j(V)$ for covering sieves of the
Grothendieck topology that is generated by the coverage are also weak equivalences in the 
left Bousfield localization just at the coverage sieves. It follows that this coincides with 
the localization at the full Grothendieck topology.
\endofproof

\addcontentsline{toc}{section}{References}

\end{document}